\documentclass[11pt,leqno]{article}

\usepackage{graphicx}
\usepackage{latexsym} 
\usepackage{amsmath}
\usepackage{amssymb}
\usepackage{amsfonts}
\usepackage{enumerate}
\usepackage{theorem}

\theoremstyle{plain}
\newtheorem{teo}{Theorem}[section]
\newtheorem{lem}[teo]{Lemma}

\newtheorem{prop}[teo]{Proposition}

{\theorembodyfont{\rmfamily}
\newtheorem{defin}[teo]{Definition}
\newtheorem{oss}[teo]{Remark}

\newtheorem{claim}{Claim}}

\renewcommand{\eqref}[1]{\textnormal{(\ref{#1})}}

\numberwithin{equation}{section}

\newcommand{\mcd}{\mathcal{D}}
\newcommand{\cvd}{\hfill$\square$}
\newcommand{\proof}[1]{\noindent\textsc{Proof#1}}
\newcommand{\rmi}{\mathrm{i}}
\newcommand{\rme}{\mathrm{e}}
\newcommand{\rmd}{\mathrm{d}}

\graphicspath{{./PaperFigs/}}

\title{A variational approach to the inverse photolithography problem}

\author{Luca \textsc{Rondi}\thanks{Dipartimento di Matematica e Geoscienze,
Universit\`a degli Studi di Trieste, via Valerio, 12/1, 34127
Trieste, Italy. E-mail: \texttt{rondi@units.it}} \and
Fadil \textsc{Santosa}\thanks{School of Mathematics, University of Minnesota,
127 Vincent Hall, 206 Church St. SE,
Minneapolis, MN 55455, USA. E-mail: \texttt{santosa@math.umn.edu} }\and Zhu \textsc{Wang}\thanks{Department of Mathematics, University of South Carolina, 
426 LeConte College, 1523 Greene Street,
Columbia, SC 29208, USA. E-mail: \texttt{wangzhu@math.sc.edu}} }

\date{}

\begin{document}

\maketitle

\setcounter{section}{0}
\setcounter{secnumdepth}{2}

\begin{abstract}
Photolithography is a process in the production of integrated circuits
in which a mask is used to create an exposed pattern with a desired
geometric shape. In the inverse problem of photolithography, a desired
pattern is given and the mask that produces an exposed pattern which
is close to the desired one is sought.  We propose a variational
approach formulation of this shape design problem and introduce a
regularization strategy. The main novelty in this work is the
regularization term that makes the thresholding operation involved
in photolithography stable.  The
potential of the method is demonstrated in numerical experiments.

\end{abstract}

\section{Introduction}\label{intro}

Photolithography is a key step in the production of integrated
circuits.  We provide a brief description of the process and refer the
reader to a more detailed readable account in \cite{Schell03}.

Integrated circuits are created in layers.  The circuit layout in
each layer is made by first treating the substrate with a
photo-resist.  A pattern in transferred to the photo-resist using
ultraviolet (UV) light and a mask.  The UV light, diffracted by
the mask, selects a pattern on
the photo-resist that is to be removed.  Once the pattern is removed,
the substrate without the photo-resist is then etched.

The mask can be viewed as an opaque screen with cut-outs. UV light from a source
goes through a system of lenses and is diffracted by the mask.
The diffracted light creates an image on the photo-resist which is placed at
the focal distance from the lenses.  The photo-resist is light sensitive.
Parts that are exposed to image intensity greater than some threshold can
be removed.  

For the purpose of this work, we call the pattern we wish to remove the `target
pattern'.  For a given mask, the exposed pattern is the set of points on the 
photo-resist where the UV light intensity is greater than some threshold.
The inverse problem in photolithography is the problem of finding the mask
that produces an exposed pattern that is as close to the target pattern as possible.
Such an inverse problem can be thought of as a shape design problem.

The mask is a set, and it may be represented
by its characteristic function $m$, that is by a binary function.
Let $I=I(m)$ be the light intensity on
the photo-resist plane for a given mask $m$.  The exposed region is given by
$$\Omega(m)=\{x:\ I(x)>h\},$$
where $h$ is the threshold.  Thus $\Omega(m)$
is the suplevel set of the real valued function $I$ at level $h$. 
Such a thresholding operation, besides being highly nonlinear, is not stable 
for instance with respect to variations of the threshold value $h$, 
in particular from the topological point of view, whenever $h$ is close to a 
critical value for $I$. 
Notice that in order to describe $\Omega$ we can use again its characteristic function, namely
$$\chi_{\Omega}=\mathcal{H}(I-h),$$
where $\mathcal{H}$ is the Heaviside function. The fact that $\mathcal{H}$ is not 
differentiable is another issue that has to be taken into account for the numerics.

Finally, the operator that maps the mask to the corresponding light
intensity on the photo-resist plane is smoothing, therefore a perfect
agreement with the target pattern might be impossible, especially if
it has some corners. This is the reason why we set the problem as an
optimal design problem.

Cobb \cite{Cobb98} was the first to tackle this problem from the
point of view of optimal design, using a physically-based
model. This approach was further developed first by using a level set method, 
\cite{Tuzel09}, and then by a variational method, \cite{Ron-San}.
In \cite{Poona07} a different computational method, where 
the mask is modelled as a pixelated binary image, is proposed.

Our starting point is the variational approach developed in \cite{Ron-San} by two of the authors. Given a desired circuit $\Omega_0$ we wish to find a mask $m$ minimizing the distance, in a suitable sense, of $\Omega(m)$ from $\Omega_0$. In order for the mask to be constructed in a relatively easy way, we require it to be not too irregular, therefore we add a perimeter penalization on the mask $m$. In \cite{Ron-San} a suitable approximation, in the sense of $\Gamma$-convergence, of the resulting functional was proposed. Such an approximation was amenable to computation using, for example, finite difference
approximations on structured grids and steepest descent for minimization, and was based on approximating binary functions $m$ by so-called phase-field functions $u$ taking values in $[0,1]$ and extending the intensity functional $I$ to be defined on phase-field and not only binary functions. The approximation of the perimeter penalization used there was the one developed in this phase-fields framework by Modica and Mortola, \cite{Mod e Mor}. We recall here that the same idea lies in the approximation of the Mumford-Shah functional due to Ambrosio and Tortorelli. Furthermore, also the Heaviside function was replaced by a smooth approximation. 

In order to apply the analysis in \cite{Ron-San}, a crucial point is that the threshold $h$ is not a critical value of the intensity. In \cite{Ron-San} this was obtained by imposing suitable technical restrictions to the model used. Instead in this paper we greatly improve the results in \cite{Ron-San} because we allow an extremely general model, that includes the one usually used in the industry which is based on the so-called Hopkins
aerial intensity representation. In fact we are able to carry over the analysis by a adding a further penalization term, the main theoretical novelty of this work. Such a regularization term, which we call $\mathcal{R}$ and that it is applied to the intensity $I$, has the aim to penalize critical points at values close to the threshold and has two important effects. From the theoretical point of view, it allows the development of the analysis, and, from the practical point of view, it allows the reconstruction to be more stable, especially from the topological point of view, with respect to variations of the threshold, that is with respect to errors in the evaluation of the threshold value.

Using the approximation developed in \cite{Ron-San} for the distance, the perimeter penalization and the Heaviside function, and devising a suitable approximation for the regularization term $\mathcal{R}$, we construct an approximated functional which is still amenable to computation. We compute its gradient, at least in a discretised version of it, and we test it by numerical experiments. The tests show that the method performs rather well, leading to reconstructed circuits that are good approximations of the desired ones.

The plan of the paper is the following. After a brief discussion of the mathematical preliminaries, Section~\ref{prelsec}, we introduce the inverse photolithography problem
and develop our variational approach, Section~\ref{model}. In Section~\ref{numsec} we present our numerical experiments. Final comments and conclusions are in Section~\ref{conclusionsec}.

\section{Mathematical preliminaries}\label{prelsec}

The following notation will be used. For every $x\in\mathbb{R}^2$, we shall set $x=(x_1,x_2)$, where $x_1$ and $x_2\in\mathbb{R}$. For every $x\in\mathbb{R}^2$ and $r>0$, we shall denote by $B_r(x)$ the open ball
in $\mathbb{R}^2$ centered at $x$ of radius $r$.
Usually we shall write $B_r$ instead of $B_r(0)$.
For any set $E\subset \mathbb{R}^2$, we denote by $\chi_E$
its characteristic function, and for any $r>0$,
$B_r(E)=\bigcup_{x\in E}B_r(x)$.

For any $f\in\mathcal{S}'(\mathbb{R}^2)$, the space of tempered distributions,
we denote by $\hat{f}$ its \emph{Fourier transform}, which, if
$f\in L^1(\mathbb{R}^2)$, may be written as
$$\hat{f}(\xi)=\int_{\mathbb{R}^2}f(x)\rme^{-\rmi \xi\cdot x}\rmd x,\quad \xi\in\mathbb{R}^2.$$
We recall that $f(x)=(2\pi)^{-2}\hat{\hat{f}}(-x)$, that is, when also $\hat{f}\in L^1(\mathbb{R}^2)$,
$$f(x)=\frac{1}{(2\pi)^2}\int_{\mathbb{R}^2}\hat{f}(\xi)\rme^{\rmi \xi\cdot x}\rmd \xi,\quad x\in\mathbb{R}^2.$$



For any function $f$ defined on $\mathbb{R}^2$ and any positive constant $s$, we denote $f_s(x)=s^{-2}f(x/s)$,
$x\in \mathbb{R}^2$.
We note that $\|f_s\|_{L^1(\mathbb{R}^2)}=\|f\|_{L^1(\mathbb{R}^2)}$ and $\widehat{f_s}(\xi)=\hat{f}(s\xi)$, $\xi\in\mathbb{R}^2$.

By $\mathcal{H}^1$ we denote the $1$-dimensional Hausdorff measure and by
$\mathcal{L}^2$ we denote the $2$-dimensio\-nal Lebesgue measure. We recall that,
if $\gamma\subset\mathbb{R}^2$ is a smooth curve,
then $\mathcal{H}^1$ restricted to
$\gamma$ coincides with its arclength.
For any Borel $E\subset\mathbb{R}^2$ we denote $|E|=\mathcal{L}^2(E)$.

Let $\mcd$ be a bounded open set contained in $\mathbb{R}^2$, with boundary $\partial \mcd$.
We say that $\mcd$ has a \emph{Lipschitz boundary}
if for every $x=(x_1,x_2)\in\partial\mcd$ there exist a Lipschitz
function $\varphi:\mathbb{R}\to\mathbb{R}$ and a positive constant $r$
such that for any $y\in B_r(x)$ we have, up to a rigid transformation,
\begin{equation}\label{Lipdomain}
y=(y_1,y_2)\in\mcd\quad \text{if and only if}\quad  y_2<\varphi(y_1).
\end{equation}

We note that $\mcd$ has a finite number of connected components, whereas
$\partial \mcd$ is formed by a finite number of rectifiable Jordan curves, therefore
$\mathcal{H}^1(\partial\mcd)=\mathrm{length}(\partial\mcd)<+\infty$.

For any integer $k=0,1,2,\ldots$, any $\alpha$, $0< \alpha\leq 1$, and any positive constants $r$ and $L$,
we say that a bounded open set $\mcd\subset\mathbb{R}^2$
is $C^{k,\alpha}$ \emph{with constants} $r$ \emph{and} $L$ if 
for every $x\in\partial\mcd$ there exists a $C^{k,\alpha}$
function $\varphi:\mathbb{R}\to\mathbb{R}$, with $C^{k,\alpha}$ norm bounded by $L$, such that for any
$y\in B_r(x)$, and up to a rigid transformation, \eqref{Lipdomain} holds. We note that we shall often use the notation \emph{Lipschitz} instead of $C^{0,1}$.

Let us fix three positive constants $r$, $L$ and $R$.
For any integer $k=0,1,2,\ldots$ and any $\alpha$, $0< \alpha\leq 1$, we denote with
$\mathcal{A}^{k,\alpha}(r,L,R)$ the class of all bounded open sets, contained in $B_R\subset \mathbb{R}^2$,
which are $C^{k,\alpha}$ with constants $r$ and $L$.

We recall some basic properties
of functions of bounded variation and sets of finite perimeter. For a more comprehensive treatment of
these subjects see, for instance, \cite{Amb e Fus e Pal}.

Given a bounded open set $\mcd\subset \mathbb{R}^2$,
we denote by $BV(\mcd)$ the Banach space of \emph{functions of bounded
variation}. We recall that $u\in BV(\mcd)$ if and only if  
$u\in L^1(\mcd)$ and its distributional derivative $Du$ is a bounded
vector
measure. We endow $BV(\mcd)$ with the standard norm as follows. Given
$u\in BV(\mcd)$, we denote by $|Du|$ the total variation of its
distributional derivative and
we set $\|u\|_{BV(\mcd)}=\|u\|_{L^1(\mcd)}+|Du|(\mcd)$.

We say that a sequence of $BV(\mcd)$ functions $\{u_h\}_{h=1}^{\infty}$
\emph{weakly}$^*$ \emph{converges} in $BV(\mcd)$ to $u\in BV(\mcd)$ if and only if
$u_h$ converges to $u$ in $L^1(\mcd)$ and $Du_h$
weakly$^*$ converges to $Du$ in $\mcd$, that is
\begin{equation}\label{weakstarconv}
\lim_{h}\int_{\mcd}v \rmd Du_h=
\int_{\mcd}v \rmd Du,\quad\text{for any }v\in C_0(\mcd).
\end{equation}
We recall that if
a sequence of $BV(\mcd)$ functions $\{u_h\}_{h=1}^{\infty}$ is bounded in $BV(\mcd)$ and converges to $u$ in $L^1(\mcd)$, then
$u\in BV(\mcd)$ and $u_h$ converges to $u$
weakly$^*$ in $BV(\mcd)$.

We say that a sequence of $BV(\mcd)$ functions $\{u_h\}_{h=1}^{\infty}$
\emph{strictly converges} in $BV(\mcd)$ to $u\in BV(\mcd)$ if and only if
$u_h$ converges to $u$ in $L^1(\mcd)$ and $|Du_h|(\mcd)$
converges to $|Du|(\mcd)$. Indeed, for any $a>0$,
\begin{equation}\label{d_st}
d_{st}(u,v)=\int_{\mcd}|u-v|+a\big||Du|(\mcd)-|Dv|(\mcd) \big|
\end{equation}
is a distance on $BV(\mcd)$ inducing the strict convergence.
We also note that strict convergence implies weak$^*$ convergence.

We recall that if $\mcd$ is a bounded open set with Lipschitz boundary,
then for any $C>0$ the set
$\{u\in BV(\mcd):\ \|u\|_{BV(\mcd)}\leq C\}$
is a compact subset of $L^1(\mcd)$.

Let $E$ be a bounded 
Borel set contained in $B_R\subset \mathbb{R}^2$. We shall denote by $\chi_E$ its characteristic function. We notice that 
$E$ is compactly
contained in $B_{R+1}$, which we shall denote by $E\Subset B_{R+1}$.
We say that $E$ is a \emph{set of finite perimeter} if
$\chi_E$ belongs to $BV(B_{R+1})$ and we call the number
$P(E)=|D\chi_E|(B_{R+1})$ its \emph{perimeter}.

Let us finally remark that the intersection of two sets of finite perimeter is
still a set of finite perimeter. Moreover,
whenever $E$ is open and $\mathcal{H}^1(\partial E)$ is finite, then $E$ is a set of finite
perimeter. In particular
a bounded open set
$\mcd$ with Lipschitz boundary
is a set of finite perimeter and its perimeter $P(\mcd)$ coincides with
$\mathcal{H}^1(\partial\mcd)$.

We conclude this preliminary section by describing a classical  $\Gamma$-convergence approximation of the perimeter functional
due to Modica and Mortola, \cite{Mod e Mor}. For the definition and properties of $\Gamma$-convergence we refer to \cite{DaM}.
%
Throughout the paper, for any $p$, $1\leq p\leq +\infty$, we shall denote its conjugate exponent by $p'$, that is $p^{-1}+(p')^{-1}=1$.

\begin{teo}\label{Mod-Morteo}
Let us fix $R>0$. Let $1<p<+\infty$ and $W:\mathbb{R}\to[0,+\infty)$ be a continuous function such that
$W(t)=0$ if and only if $t\in\{0,1\}$. Let $c_p=(\int_0^1(W(s))^{1/p'}\mathrm{d}s)^{-1}$.

For any $\varepsilon>0$ we define the functional
$\mathcal{P}_{\varepsilon}:L^1(\mathbb{R}^2)\to [0,+\infty]$ as follows
\begin{equation}\label{modmordef}
\mathcal{P}_{\varepsilon}(u)=\left\{
\begin{array}{ll}
\displaystyle{\frac{c_p}{p'\varepsilon}\int_{B_R}W(u)+\frac{c_p\varepsilon^{p-1}}{p}\int_{B_R}|\nabla u|^p}&\text{if }u\in W^{1,p}_0(B_R),\\ &\quad u=0\text{ a.e. outside }B_R,\\
\vphantom{\displaystyle{\int}}+\infty&\text{otherwise}.
\end{array}
\right.
\end{equation}

Let $\mathcal{P}:L^1(\mathbb{R}^2)\to [0,+\infty]$ be such that
\begin{equation}\label{Pdef}
\mathcal{P}(u)=
\left\{
\begin{array}{ll}
\vphantom{\displaystyle{\int}}|Du|(B_{R+1})
&\text{if }u\in BV(B_{R+1}),\ u\in\{0,1\}\text{ a.e.},\\ &\quad
u=0\text{ a.e. outside }B_R,\\
\vphantom{\displaystyle{\int}}+\infty&\text{otherwise}.
\end{array}
\right.
\end{equation}

Then $\mathcal{P}=\Gamma\textrm{-}\!\lim_{\varepsilon\to 0^+}\mathcal{P}_{\varepsilon}$ with respect to the $L^1(\mathbb{R}^2)$ norm.
\end{teo}

\begin{oss}
We observe that $\mathcal{P}(u)=P(E)$ if $u=\chi_E$ where $E$ is a set of finite perimeter contained in $B_R$ and $\mathcal{P}(u)=+\infty$ otherwise.

Furthermore, we note that the result does not change if in the definition of $\mathcal{P}_{\varepsilon}$ we set $\mathcal{P}_{\varepsilon}(u)=+\infty$ whenever $u$ does not satisfy the
constraint
\begin{equation}
0\leq u\leq 1\text{ a.e. in }B_R.
\end{equation}
Actually, in the numerics we shall always implicitly impose such a constraint.
\end{oss}


Also the following result, due to Modica, \cite{Mod}, will be useful.

\begin{prop}\label{compactnessoss}
Let us consider
any family $\{u_{\varepsilon}\}_{0<\varepsilon\leq \tilde{\varepsilon}}$ such that,
for some positive constant $C$ and
for any $\varepsilon$, $0<\varepsilon\leq\tilde{\varepsilon}$,
we have $0\leq u_{\varepsilon}\leq 1$ almost everywhere and
$\mathcal{P}_{\varepsilon}(u_{\varepsilon})\leq C$.
Then $\{u_{\varepsilon}\}_{0<\varepsilon\leq \tilde{\varepsilon}}$
is precompact in $L^1(\mathbb{R}^2)$.
\end{prop}

\section{The inverse problem and its approximation}\label{model}

Kirchhoff approximation is presently favored as a modeling tool for
the optical phenomena in photolithography.  This is due to the fact
that Kirchhoff approximation can be very efficiently computed and it
is relatively accurate.  It is true however that more accurate
optical modeling may be needed in the future.  
Under this approximation, the open portions of the mask acts as
light sources; the amplitude of light at the mask opening is that
of the incident field from the light source.
Propagation through the lenses can be
calculated using Fourier optics.  It is further assumed that the image
plane, in this case the plane of the photo-resist, is at the focal
distance of the optical system.  If there were no
diffraction, a perfect image of the mask would be formed on the image
plane. Diffraction, together with partial coherence of the light
source, acts to distort the formed image.

The mask, which consists of cut-outs, is represented as a 
binary function, the
characteristic function of the cut-outs $D$. Namely the mask is given by
\[
m(x) = \chi_D(x).
\]
The light intensity on the image plane is given by
\cite{Pati97}
\begin{equation}\label{Hopkins}
I(x) = \int_{\mathbb{R}^2} \int_{\mathbb{R}^2} m(\xi)K(x-\xi) J(\xi-\eta) \overline{K}(x-\eta) m(\eta) \rmd\xi \rmd\eta,\quad x\in\mathbb{R}^2. 
\end{equation}
In the above expression the kernel $K(\cdot)$ is called the \emph{coherent point spread function}
and describes the optical system. The function $J(\cdot)$ is called the \emph{mutual
intensity function}.  If the illumination is fully coherent then $J\equiv 1$ but in practice illumination is never fully coherent.
The equation (\ref{Hopkins}) is often referred to as the Hopkins aerial
intensity representation.

\subsubsection{Assumptions on $K$ and $J$.}
We assume that $K$ is a complex valued function such that for a constant $\alpha$, $0<\alpha\leq 1$, we have $K\in C^{1,\alpha}(\mathbb{R}^2)$.
Furthermore, we assume that $|K|$ converges to $0$ uniformly as $\|x\|\to +\infty$,  that is for any $\varepsilon>0$ there exists $r>0$ such that for any $x\in \mathbb{R}^2$ with $\|x\|\geq r$ we have $|K(x)|\leq \varepsilon$.

We assume that $J$ is the Fourier transform of a function $\hat{J}$ such that
$\hat{J}\in L^1(\mathbb{R}^2)$ and $\hat{J}\geq 0$ almost everywhere in $\mathbb{R}^2$.
In particular $J$ is a continuous complex valued function.

\bigskip

A typical model for $K$ and $J$ is the following.
For an optical system with a
circular aperture, once the wavenumber of the light
used, $k>0$, has
been chosen, the kernel depends on a single parameter called the Numerical
Aperture, $\text{NA}$. Notice that the wavelength
is $\lambda=2\pi/k$.
Let us recall that the so-called Jinc function is defined as
\[
\textrm{Jinc}(x)=\frac{J_1(|x|)}{2\pi|x|},\quad x\in\mathbb{R}^2,
\]
where $J_1$ is the Bessel function of order 1. We notice that in the 
Fourier space
$$\widehat{\textrm{Jinc}}(\xi)=\chi_{B_1}(\xi),\quad \xi\in\mathbb{R}^2.$$

If we denote by $s=(k\text{NA})^{-1}$, then
the kernel
is usually modeled as follows
\begin{equation}\label{Kdefin}
K(x)=\textrm{Jinc}_s(x)=\frac{k\text{NA}}{2\pi}\frac{J_1(k \text{NA}|x|)}{|x|},\quad x\in\mathbb{R}^2,
\end{equation}
therefore
$$\hat{K}(\xi)=\chi_{B_1}(s\xi)=\chi_{B_{1/s}}(\xi)=\chi_{B_{k\text{NA}}}(\xi),\quad \xi\in\mathbb{R}^2.$$
If $\text{NA}$ goes to $+\infty$, that is $s\to 0^+$, then $\hat{K}$ converges pointwise to $1$, thus $K$ approximates in a suitable sense the Dirac delta.

The mutual
intensity function $J(\cdot)$ is parametrized
by a \emph{coherency coefficient} $\sigma$.  A typical model for $J$ is
\begin{equation}\label{Jdef}
J(x)=2\frac{J_1(k \sigma \text{NA}|x|)}{k\sigma\text{NA}|x|}=\pi \textrm{Jinc}(k\sigma\text{NA}|x|),\quad x\in\mathbb{R}^2.
\end{equation}
Thus,
\begin{equation}\label{Jdefstar}
\frac{1}{(2\pi)^2}\hat{J}(\xi)=\frac{1}{\pi(k\sigma\text{NA})^2}\chi_{B_{k\sigma\text{NA}}}(\xi),\quad\xi\in\mathbb{R}^2,
\end{equation}
that, as $\sigma\to 0^+$, converges, in a suitable sense, to the Dirac delta. Therefore
full coherence is achieved for $\sigma\to 0^+$. In fact, if $\sigma\to 0^+$, $J$ converges to $1$ uniformly on any compact subset of $\mathbb{R}^2$.

\bigskip

The photo-resist material responds to the intensity of the image.  When
intensity at the photo-resist goes over a certain threshold, it is
then considered exposed and can be removed.  Therefore, the
exposed pattern, given a mask $m(x)$, is
\begin{equation}\label{exposed}
\Omega = \{ x\in\mathbb{R}^2\,:\ I(x) > h \},
\end{equation}
where $h>0$ is the exposure threshold.  Clearly, $\Omega$ depends on the
mask function $m(x)$, which we recall is given by the characteristic
function of $D$ representing the cut-outs, that is $\Omega=\Omega(D)$.  In photolithography, we have a
desired exposed pattern $\Omega_0$ which we wish to achieve.  The inverse problem
is to find a mask that achieves this desired exposed pattern, that is to find $D$ such that $\Omega(D)=\Omega_0$.
Mathematically, this cannot, in general, be done.  Therefore, the
inverse problem must be posed as an optimal design problem.

\subsubsection{Assumptions on the target pattern $\Omega_0$.}
Let us fix $R>0$.
We assume that $\Omega_0$ is a bounded open set compactly contained in $B_R$ such that $\Omega_0$ is a set of finite perimeter.

\bigskip

Suppose the desired pattern is given by $\Omega_0$.  We pose the
minimization problem
\begin{equation}\label{design_prob}
\displaystyle{\min_{D\in{\mathcal A}}} \; d(\Omega(D), \Omega_0).
\end{equation}
For what concerns the distance function $d(\cdot,\cdot)$ and the admissible set
$\mathcal A$, we shall choose the following.
We set
$$\mathcal{A}=\{E\subset B_R:\ E\text{ is a set of finite perimeter}\}$$
and for any $E\in\mathcal{A}$ we denote by $P(E)$ its perimeter and we notice that
$$P(E)=\mathcal{P}(\chi_E),$$
where, for any function $u\in L^1(\mathbb{R}^2)$, $\mathcal{P}(u)$ is defined in \eqref{Pdef}.
With a slight abuse of notation we shall identify sets with their characteristic functions, so that
$\mathcal{A}$ may also denote
$$\mathcal{A}=\{u\in L^1(\mathbb{R}^2):\ \mathcal{P}(u)<+\infty\}.$$

About the distance we shall choose
\begin{equation}
d(\Omega_1,\Omega_2)=d_{st}(\chi_{\Omega_1},\chi_{\Omega_2})=\int|\chi_{\Omega_1}-\chi_{\Omega_2}|+a|P(\Omega_1)-P(\Omega_2)|,\label{d3def}
\end{equation}
where $a$ is a positive tuning parameter.
We recall that in \cite[Section~3.3]{Ron-San} the choice of the distance has been thoroughly discussed.

We shall add to \eqref{design_prob} two regularization terms.
The first one is on the independent variable $D$, that is on the mask. To ensure manufacturing of the mask, the optimal mask may not be too irregular, therefore we shall add a perimeter penalization on the mask. The second regularization term allows us to stabilize the optimization procedure. In fact, the thresholding operation that, given the intensity, determines the target domain is not stable. For instance, if $h$ is a critical value of the intensity $I$, a very small modification of the mask might lead to a change in the topology of the reconstructed circuit. In order to avoid this, we shall
discard masks such that $h$ is close to a critical value of the corresponding intensity $I$. We shall achieve this aim by adding a second penalization term $\mathcal{R}$ which we describe later in this section.

Let us set up the regularized minimization problem. 
We denote
$$A=\{u\in L^1(\mathbb{R}^2):\ 0\leq u\leq 1\text{ a.e. in }\mathbb{R}^2\text{ and }u=0 \text{ a.e outside }B_R\}.$$
Let us define the following operator $\mathcal{I}:A\to C^0(\mathbb{R}^2)$
such that for any $u\in A$ we have
\begin{equation}\label{mathcalPdef}
\mathcal{I}(u)(x) = \int_{\mathbb{R}^2} \int_{\mathbb{R}^2} u(\xi)K(x-\xi) J(\xi-\eta) \overline{K}(x-\eta) u(\eta) \rmd\xi \rmd\eta,\quad x\in\mathbb{R}^2. 
\end{equation}
Let us note that if $u=\chi_D$ for a mask $D\subset B_R$, then $\mathcal{I}(u)$ coincides with the intensity
$I$ as defined in \eqref{Hopkins}.

In the next proposition we describe some of the properties of $\mathcal{I}$.

\begin{prop}\label{regularityprop}
Under our assumptions on $K$ and $J$, the following holds.
\begin{enumerate}[\textnormal{(}i\textnormal{)}]
\item\label{part1} For any $u\in A$, $\mathcal{I}(u)$ is a real valued function such that $\mathcal{I}(u)\geq 0$
in $\mathbb{R}^2$. Obviously,
if $u$ is identically equal to zero, then also $\mathcal{I}(u)$ is identically equal to zero.
\item\label{part2} For any $u\in A$, $\mathcal{I}(u)\in C^{1,\alpha}(\mathbb{R}^2)$ and, for any $R_1>0$, there exists $C>0$ such that
$$\|\mathcal{I}(u)\|_{C^{1,\alpha}(B_{R_1})}\leq C,\quad\text{for any }u\in A.$$
\item\label{part3} For any $R_1>0$, $\mathcal{I}$ is uniformly continuous with respect to the $L^1$ norm on $A$ and the $C^{1,\alpha}$ norm on $C^0(B_{R_1})$.
\item\label{part4} $\mathcal{I}(u)$ converges to $0$ uniformly as $\|x\|\to +\infty$, uniformly with respect to $u\in A$, that is for any $\varepsilon>0$ there exists $r>0$ such that for any $x\in \mathbb{R}^2$ with $\|x\|\geq r$ and any $u\in A$ we have $|\mathcal{I}(u)(x)|\leq \varepsilon$.
\end{enumerate}

\end{prop}

\proof{.}
For any $u\in A$, we define $U\in L^1(\mathbb{R}^4)$ as follows
$$U(x,y)=u(x)u(y)J(x-y),\quad\text{for any }x,y\in\mathbb{R}^2.$$
Then we define 
$H\in C^{1,\alpha}(\mathbb{R}^4)$ in the following way
$$H(x,y)=K(x)\overline{K}(y),\quad\text{for any }x,y\in\mathbb{R}^2.$$
We notice that
$$\mathcal{I}(u)(x)=(H\ast U)(x,x),\quad \text{for any }x\in\mathbb{R}^2,$$
where $\ast$ denotes convolution, in this particular case in $\mathbb{R}^4$.

Therefore, parts (\ref{part2}), (\ref{part3}) and (\ref{part4}) follow immediately from standard properties of convolutions. For what concerns (\ref{part1}), this requires a little more care.
We call $T(u)=K\ast u$, where again $\ast$ denotes convolution, in this case in $\mathbb{R}^2$. Then, for any $k\in\mathbb{R}^2$, we denote
$\tilde{u}_k(x)=u(x)\rme^{\rmi k\cdot x}$ for any $x\in\mathbb{R}^2$.
Then, fixed $x\in \mathbb{R}^2$, we define the function 
$$f(k)=|T(\tilde{u}_k)(x)|^2=\int_{\mathbb{R}^2} \int_{\mathbb{R}^2} u(\xi)\rme^{\rmi k\cdot \xi}K(x-\xi) \overline{K}(x-\eta) u(\eta)\rme^{-\rmi k\cdot \eta} \rmd\xi \rmd\eta,\quad x\in\mathbb{R}^2.$$
Clearly $f(k)$ is nonnegative for any $k\in \mathbb{R}^2$, therefore it would be enough to show that
$$\mathcal{I}(u)(x)=\frac{1}{(2\pi)^2}\int_{\mathbb{R}^2} \hat{J}(k)f(k)\rmd k$$ and this follows simply by Fubini theorem.\cvd

\begin{oss}\label{other_assumptions}
We may replace the assumptions on $J$ by the following one and still the previous proposition and all the next results hold. We may assume that
$J$ is the Fourier transform of a tempered distribution $\hat{J}$ such that
$\hat{J}$ has compact support and it is semipositive definite,
that is $\langle \hat{J},f\rangle\geq 0$ for any $f\in \mathcal{S}$ such that $f\geq 0$ everywhere in $\mathbb{R}^2$. We notice that again $J$ is a continuous complex valued function, therefore the proof of parts (\ref{part2}), (\ref{part3}) and (\ref{part4}) is exactly the same. The argument for proving (\ref{part1}) is slightly more involved in this case and we leave the details to the reader.
\end{oss}

We denote by $\mathcal{H}:\mathbb{R}\to\mathbb{R}$ the Heaviside function
such that $\mathcal{H}(t)=0$ for any $t\leq 0$ and $\mathcal{H}(t)=1$ for any $t> 0$.
For any positive constant $h$ we set $\mathcal{H}_h(t)=\mathcal{H}(t-h)$
for any $t\in\mathbb{R}$.
Then, fixed the threshold $h>0$, we define
$\mathcal{W}:A\to L^{\infty}(\mathbb{R}^2)$ as follows
\begin{equation}\label{Wdefinition}
\mathcal{W}(u)=\mathcal{H}_h(\mathcal{I}(u)),\quad \text{for any }u\in A.
\end{equation}
Clearly, for any $u\in A$, $\mathcal{W}(u)$ is the characteristic function of an open set,
which we shall call $\Omega(u)$. That is
\begin{equation}\label{Omegedefinition}
\Omega(u)=\{x\in\mathbb{R}^2\,:\ \mathcal{I}(u)(x)>h\},\quad \text{for any }u\in A.
\end{equation}
In other
words, $\chi_{\Omega(u)}=\mathcal{W}(u)=\mathcal{H}_h(\mathcal{I}(u))$.
Moreover, whenever $u=\chi_E$, where $E$ is a measurable subset of $B_R$,
we shall denote $\Omega(E)=\Omega(\chi_E)$.

In order to define the regularization term $\mathcal{R}$,
we need a few auxiliary functions. Let us fix a positive constant $\delta_0$.
Let
$f:\mathbb{R}\to [0,+\infty]$ be a continuous function satisfying the following properties
\begin{enumerate}[(i)]
\item $f$ is identically equal to $+\infty$ on $(-\infty,0]$;
\item $f$ is decreasing;
\item $f$ is identically equal to zero on 
$[\delta_0,+\infty)$;
\item the following behaviour at $0$ holds
$$\lim_{s\to 0^+}f(s)s^{2/\alpha}\geq C>0.$$
\end{enumerate}

Let $\varphi:\mathbb{R}\to\mathbb{R}$ be a $C^1$
function such that
\begin{enumerate}[(i)]
\item $\varphi(h)>0$;
\item $\varphi$ is increasing before $h$ and decreasing after $h$;
\item $\varphi(0)<-\delta_0$.
\end{enumerate}

We notice that we can find positive constants $\delta$ and $c_1$ such that
$\varphi$ is greater than or equal to $c_1$ on $[h-\delta,h+\delta]$ and $\varphi(0)\leq -(\delta_0+\delta)$.

For example, we may choose
\begin{equation}\label{fdefinition}
f(s)=\rme^{-\delta_0^2/(\delta_0^2-s^2)}s^{-2/\alpha},\quad \text{for any }0<s<\delta_0,
\end{equation}
with $f\equiv +\infty$ on $(-\infty,0]$ and $f\equiv 0$ on $[\delta_0,+\infty)$,
and
\begin{equation}\label{varphidefin}
\varphi(s)=-a(s-h)^2+b,\quad\text{for any }s\in\mathbb{R},
\end{equation}
for suitable positive constants $a$ and $b$. For instance, if we pick $a>0$ and $\delta_0$, $0<\delta_0\leq h/2$, such that $a\delta_0\geq 1$, then we may take $\delta=c_1=\delta_0$ provided $b$ satisfies
$\delta_0+a\delta_0^2  \leq b \leq ah^2-2\delta_0$.

\begin{defin}\label{regularizationdefin}
Let us define $\mathcal{R}:A\to [0,+\infty]$ as follows
$$\mathcal{R}(u)=\int_{\mathbb{R}^2}f\left(\|\nabla (\mathcal{I}(u))\|^2-\varphi(\mathcal{I}(u))\right),\quad\text{for any }u\in A.$$
\end{defin}
To get a sense of the behavior of the regularization term $\mathcal{R}$, consider
a point $x$ for which $\mathcal{I}(u(x)) \approx h$. This means that $\varphi(I(u))>0$. If
$\| \nabla (\mathcal{I}(u(x))) \|$ is small, then $\mathcal{R}=+\infty$. The penalty
term is zero if  $\mathcal{I}(u(x))$ is away from $h$. Therefore, the term $\mathcal{R}$
does not allow the critical values of $\mathcal{I}(u)$ to be close to $h$.
 
\begin{oss}\label{tildeRremark}
All the theoretical results we are going to prove remain valid if we replace $\mathcal{R}$
with $\tilde{\mathcal{R}}:A\to [0,+\infty]$ defined as follows
$$\tilde{\mathcal{R}}(u)=\mathcal{R}(u)+\frac{1}{|\{h-\delta\leq \mathcal{I}(u)\leq h+\delta\}|},
\quad\text{for any }u\in A.$$
\end{oss}
 
\begin{prop}\label{reglemma}
Under the previous notation and assumptions, we have that the functional $\mathcal{R}:A\to [0,+\infty]$
is continuous, with respect to the $L^1$ convergence in $A$.
\end{prop}

Before proving Proposition~\ref{reglemma}, we need the following.

\begin{lem}\label{domainlemma}
Under the previous notation and assumptions, there exist positive constants $r$, $L$ and $R_1\geq R$ such that , for any $u\in A$ satisfying $\mathcal{R}(u)<+\infty$, and
for any $t\in [h-\delta,h+\delta]$, we have that
$$\Omega(t)=\{\mathcal{I}(u)>t\}$$
is either empty or it belongs to $\mathcal{A}^{1,\alpha}(r,L,R_1)$.
\end{lem}

\proof{.} Since $\mathcal{I}(u)$
decays to zero at infinity, uniformly with respect to $u\in A$, we observe that there exists $R_1\geq R$ such that the following properties hold.
First
$\{\mathcal{I}(u)\geq h-\delta\}\subset B_{R_1}$ for any $u\in A$. Moreover, if we denote
$$g(x)=\|\nabla \mathcal{I}(u)(x)\|^2-\varphi(\mathcal{I}(u)(x)),\quad \text{for any }x \in\mathbb{R}^2,$$
we have that $g(x)\geq \delta_0$, hence
$f(g(x))=0$, for any $x\in\mathbb{R}^2$ with $\|x\|\geq R_1$ and any $u\in A$.

Fixed $u\in A$, let us define $g$ as before.
By the continuity of $g$ and the properties of $f$, if $\mathcal{R}(u)$ is finite
then $g(x)$ must be nonnegative for every $x\in\mathbb{R}^2$. This would be enough for the proof of this lemma, but actually we have that there exists a positive constant $\varepsilon$, depending on $u$, such that $g(x)\geq \varepsilon$ for every $x\in\mathbb{R}^2$. Since this property is crucial in the proof of Proposition~\ref{reglemma}, we sketch its proof here.

Since $g(x)\geq \delta_0$ for any
$x$ outside $B_{R_1}$, it would be enough to prove that $g(x)>0$ for any $x\in\mathbb{R}^2$.
We argue by contradiction and we assume that
$g(x)=0$ for some $x\in B_{R_1}$.

Then, since $g$ is H\"older continuous with exponent $\alpha$, $0<\alpha\leq 1$,  on the closure of $B_{R_1+1}$, we infer that 
for any $y$ in
a neighbourhood of $x$ we have that $g(y)\leq C|y-x|^{\alpha}$,
therefore by using polar coordinates centered at $x$ we obtain, for some $r_0>0$,
$$\mathcal{R}(u)\geq 2\pi\int_{0}^{r_0}sf(Cs^{\alpha})\rmd s.$$
Since the right-hand side is $+\infty$ by our assumptions on $f$, the claim is proved.

If $c_1>0$ is the minimum of $\varphi$ on $[h-\delta,h+\delta]$, since $g(x)\geq 0$, we infer that
$$\|\nabla (\mathcal{I}(u))(x)\|^2\geq c_1\quad\text{if }\mathcal{I}(u)(x)\in [h-\delta,h+\delta].$$

Then the conclusion immediately follows by the uniform $C^{1,\alpha}$ regularity of $\mathcal{I}(u)$ proved in Proposition~\ref{regularityprop} and the implicit function theorem.\cvd

\bigskip

\proof{ of Proposition~\ref{reglemma}.}
Let $u_n\in A$, $n\in\mathbb{N}$, converge to $u$ in $L^1$ as $n\to\infty$. Clearly $u\in A$ as well.

We begin by considering the case in which $\mathcal{R}(u)$ is finite.
In the previous lemma we proved that
there exists a positive constant $\varepsilon$ such that
$g(x)\geq \varepsilon$ for any $x\in\mathbb{R}^2$.
Since $\mathcal{I}(u_n)$ converges to $\mathcal{I}(u)$ locally in $C^1$ as $n\to\infty$, 
by the dominated convergence theorem we easily deduce that $\mathcal{R}(u_n)$ converges to $\mathcal{R}(u)$ as $n\to\infty$.

By the previous lemma, for any $u\in A$ such that $\mathcal{R}(u)<+\infty$ we have that
$|\{h-\delta< \mathcal{I}(u)< h+\delta\}|=
|\{h-\delta\leq \mathcal{I}(u)\leq h+\delta\}|$.
By the uniform convergence of $\mathcal{I}(u_n)$ to $\mathcal{I}(u)$
as $n\to\infty$,
and the dominated convergence theorem, we obtain that
$|\{h-\delta\leq \mathcal{I}(u_n)\leq h+\delta\}|$ converges to  $|\{h-\delta\leq \mathcal{I}(u)\leq h+\delta\}|$ as $n\to\infty$.
Therefore we immediately conclude that also $\tilde{\mathcal{R}}$ is continuous at any $u\in A$ such that $\mathcal{R}(u)<+\infty$.

If $\mathcal{R}(u)=+\infty$, then there exists $x$ such that $g(x)\leq 0$. 
Consequently, if $g_n$ is the corresponding function related to $u_n$, we conclude that
$g_n(x)$ goes to zero as $n\to\infty$. Therefore
$$\mathcal{R}(u_n)\geq 2\pi\int_{0}^{r_0}sf(g_n(x)+Cs^{\alpha})\rmd s\to +\infty\quad\text{as }n\to\infty$$
and the proposition is proved.\cvd

\bigskip

Let $\tilde{R}=R_1+1$ with $R_1$ as in Lemma~\ref{domainlemma}. We notice that the functional $\mathcal{R}$ may be equivalently defined as
\begin{equation}\label{newR}
\mathcal{R}(u)=\int_{B_{\tilde{R}}}f\left(\|\nabla (\mathcal{I}(u))\|^2-\varphi(\mathcal{I}(u))\right),\quad\text{for any }u\in A.
\end{equation}

For any positive constant $C$, let us denote
$$\tilde{A}_C=\{u\in A:\ \mathcal{R}(u)\leq C\}.$$

\begin{lem}\label{firstlemma}
For any $C>0$,
the map $\mathcal{W}:\tilde{A}_C\to BV(B_{\tilde{R}})$ is uniformly continuous with respect to the $L^1$ norm on $\tilde{A}_C$
and the distance $d_{st}$ on $BV(B_{\tilde{R}})$.
\end{lem}

\proof{.} By Lemma~\ref{domainlemma}, we recall that, for any $u\in\tilde{A}_C$,
$\mathcal{W}(u)= \chi_{\Omega(u)}$ where $\Omega(u)$ is either empty or it belongs to
$\mathcal{A}^{1,\alpha}(r,L,R_1)$.
Then, there exists a constant $C_1$ such that
$d_{st}(\mathcal{W}(u),\mathcal{W}(v))\leq C_1$ for any $u$ and $v\in \tilde{A}_C$ .

There exists a function $\omega:[0,+\infty)\to[0,+\infty)$, which is nondecreasing and such that
$\lim_{t\to 0^+}\omega(t)=0$, such that for any $u$ and $v$ belonging to $A$ we have
$$\|\mathcal{I}(u)-\mathcal{I}(v)\|_{L^{\infty}(B_{\tilde{R}})}\leq \omega(\|u-v\|_{L^1(\mathbb{R}^2)}).$$

We need the following claim.

\begin{claim}\label{claim1}
There exists a function $\tilde{g}:[0,+\infty)\to[0,+\infty)$, which is continuous, increasing and such that $\tilde{g}(0)=0$,
satisfying the following property.
For any $u\in\tilde{A}_C$ such that $\Omega(u)$ is not empty,  
for any $\varepsilon>0$ and any $x\in\mathbb{R}^2$ we have
\begin{equation}\label{claim}
\text{if }x\not\in B_{\varepsilon}(\partial\Omega(u)),\text{ then }|\mathcal{I}(u)(x)-h|>\tilde{g}(\varepsilon).
\end{equation}
\end{claim}

To prove this claim, we recall that for a positive constant $c_1$ 
we have $\|\nabla(\mathcal{I}(u))(x)\|=
-\nabla(\mathcal{I}(u))(x)\cdot\nu\geq \sqrt{c_1}$
for any $x\in \partial\Omega(u)$, where as usual $\nu$ is the outer normal to $\Omega(u)$. The $C^{1,\alpha}$ regularity of $\mathcal{I}(u)$ on $B_{\tilde{R}}$, which is uniform with respect to $u\in A$, allows us to conclude the proof of the claim.

We conclude that there exists a positive constant $\eta_0$ such that if $u$ and $v\in \tilde{A}_C$ satisfy $\|u-v\|_{L^1(\mathbb{R}^2)}\leq \eta_0$, then $\Omega(u)$ is empty if and only if 
$\Omega(v)$ is.
If they are both empty, then $d_{st}(\mathcal{W}(u),\mathcal{W}(v))=0$. Therefore, we are interested only in the case in which they are both not empty.

We follow some of the arguments developed in \cite[Theorem~4.2]{Ron-San} which we briefly sketch for the convenience of the reader.

Let us now assume that $u$ and $v\in \tilde{A}_C$ satisfy $\|u-v\|_{L^1(\mathbb{R}^2)}\leq \eta_0$, and $\Omega(u)$ and
$\Omega(v)$ are not empty. Fixed $\varepsilon>0$, we can find $\eta>0$
such that if $\|u-v\|_{L^1(\mathbb{R}^2)}\leq \eta$, then $\|\mathcal{I}(u)-\mathcal{I}(v)\|_{L^{\infty}(B_{\tilde{R}})}\leq \tilde{g}(\varepsilon)$.

Let us now take $x\in\partial \Omega(u)$, that is $x\in\mathbb{R}^2$ such that
$\mathcal{I}(u)(x)=h$. We infer that $|\mathcal{I}(v)(x)-h|\leq \tilde{g}(\varepsilon)$, therefore by the claim
we deduce that $x\in B_{\varepsilon}(\partial\Omega(v))$. That is
$\partial \Omega(u)\subset B_{\varepsilon}(\partial\Omega(v))$. By symmetry,
we conclude that the Hausdorff distance $d_H$ between $\partial\Omega(u)$ and $\partial\Omega(v)$
is bounded by $\varepsilon$. It has been shown in Section~3.3 in \cite{Ron-San} that 
there exist positive constants $C_2$ and $\beta$ such that for any $\Omega(u)$ and $\Omega(v)$ belonging to
$\mathcal{A}^{1,\alpha}(r,L,R_1)$ we have
$$d_{st}(\Omega(u),\Omega(v))\leq C_2\left(d_H(\partial\Omega(u),\partial\Omega(v))\right)^{\beta}.$$
Therefore the thesis immediately follows.\cvd

\bigskip

We are now in the position to set up our optimization problem. Under the previous definitions and assumptions, let us define the functional
$F_0:A\to [0,+\infty]$ such that
\begin{equation}\label{F0}
F_0(u)= \int|\chi_{\Omega(u)}-\chi_{\Omega_0}|+a|P(\Omega(u))-P(\Omega_0)|
+b\mathcal{P}(u)+c\mathcal{R}(u),\quad\text{for any }u\in A,
\end{equation}
where $\mathcal{P}$ is the functional defined in \eqref{Pdef}, and $a$, $b$ and $c$ are positive tuning parameters. We notice that 
$$\int|\chi_{\Omega(u)}-\chi_{\Omega_0}|+a|P(\Omega(u))-P(\Omega_0)|=
d_{st}(\mathcal{W}(u),\chi_{\Omega_0})$$
where
$d_{st}$ is the strict convergence distance in $BV(B_{\tilde{R}})$ given in \eqref{d_st}.

We look for the solution to the following minimization problem
\begin{equation}\label{min0}
\min\{F_0(u)\,:\ u\in A\}.
\end{equation}

By the direct method, we have that $F_0$ admits a minimum on $A$. However,
in order to make the minimization problem meaningful we shall need the following.

\subsubsection{A priori assumptions on minimizers of $F_0$}
We assume that there exists $\tilde{u}\in A$ such that $F_0(\tilde{u})$ is finite and 
$F_0(\tilde{u})<|\Omega_0|+aP(\Omega_0)$.

\bigskip

By these assumptions, we exclude that the function $u_0\equiv 0$ is a minimizer of $F_0$ and we guarantee that for
any minimizer $u$ of $F_0$ we have that $\Omega(u)$ is not empty.
In fact, if 
$u\in A$ is such that $\Omega(u)$ is empty, we have that
$$F_0(u)\geq F_0(u_0)=|\Omega_0|+aP(\Omega_0)>F_0(\tilde{u}).$$

Let us notice that if instead we replace $\mathcal{R}$ with $\tilde{\mathcal{R}}$ as in Remark~\ref{tildeRremark}, we just need to assume that there exists $\tilde{u}\in A$ such that $F_0(\tilde{u})$ is finite. In fact, if $u\in A$ is such that $\Omega(u)$ is empty, we have that
$\tilde{\mathcal{R}}(u)=+\infty$ and consequently $F_0(u)=+\infty$. This follows by this simple argument. If $\mathcal{R}(u)$ is finite, then the maximum of $\mathcal{I}(u)$ is either strictly greater that
$h+\delta$, thus $\Omega(u)$ is not empty, or strictly smaller than $h-\delta$, thus $\tilde{\mathcal{R}}(u)=+\infty$,
see the proof of Lemma~\ref{domainlemma}.

 Moreover, the function $\tilde{u}\in A$ satisfying the previous assumption is the characteristic function of a set of finite perimeter and might be considered a natural choice as an initial guess for any iterative method. Hopefully, the target set $\Omega_0$ may provide such an initial guess, that is it would be desirable that the previous assumption be satisfied by $\tilde{u}=\chi_{\Omega_0}$. As we shall show in the numerical tests, actually in practice it is not always convenient to use as an initial guess the target itself or a small perturbation of it.

We conclude that, under this assumption,
if $u$ is a minimizer of $F_0$, then $u=\chi_E$ where $E$ is a set of finite perimeter and $\Omega(\chi_E)$ is not empty.
Such a set $E$ should be chosen as the optimal mask and $\Omega(\chi_E)$ would be the optimal reconstructed circuit.


The minimization of $F_0$ presents several challenges from a numerical point of view. 
Therefore  we approximate, in the sense of $\Gamma$-convergence, the functional 
$F_0$
with a family of functional $\{F_{\varepsilon}\}_{\varepsilon>0}$ which are easier
to compute with.

We recall that $h>0$ is the fixed threshold.
We take a $C^{\infty}$ function $\phi:\mathbb{R}\to\mathbb{R}$ such that $\phi$
is nondecreasing, $\phi(t)=0$ for any $t\leq -1/2$ and $\phi(t)=1$ for any $t\geq 1/2$.
For any $\eta>0$ let
$$\phi_{\eta}(t)=\phi\left(\frac{t-h}{\eta}\right),\quad\text{for any }t\in\mathbb{R}.$$

For any $\eta>0$, let
$\Phi_{\eta}:A\to C^{1,\alpha}(\mathbb{R}^2)$
be defined as
\begin{equation}\label{Phidefin}
\Phi_{\eta}(u)=\phi_{\eta}(\mathcal{I}(u)),\quad\text{for any }u\in A.
\end{equation}

Let us summarize the properties of such a function $\Phi_{\eta}$.

\begin{prop}\label{uniformcontprop}
For any $\eta>0$, let
$\Phi_{\eta}:A\to C^{1,\alpha}(\mathbb{R}^2)$
be defined as in \eqref{Phidefin}.
We have that
$\Phi_{\eta}$ is uniformly continuous with respect to the $L^1$ norm on $A$ and the $C^{1,\alpha}$ norm on 
$C^{1,\alpha}(B_{\tilde{R}})$.

Furthermore, for any $C>0$, $\Phi_{\eta}$ converges,
as $\eta\to 0^+$,  uniformly to $\mathcal{W}$ on $\tilde{A}_C$
with respect to the distance $d_{st}$ on $BV(B_{\tilde{R}})$, that is for any $\varepsilon>0$ there exists $\eta_0>0$ such that, for any $\eta$, $0<\eta\leq\eta_0$, we have
$$d_{st}(\Phi_{\eta}(u),\mathcal{W}(u))\leq\varepsilon,\quad\text{for any }u\in \tilde{A}_C.$$
\end{prop}

\proof{.} The fist part is an easy consequence of Proposition~\ref{regularityprop}.
The second and more important part may be proved in an analogous way as 
Proposition~5.1 in \cite{Ron-San}.\cvd

\bigskip

Furthermore, for any $\gamma>0$, let us define $\mathcal{R}_{\gamma}:A\to[0,+\infty]$ satisfying the following properties. First,  $\mathcal{R}_{\gamma}$ is lower semicontinuous, with respect to the $L^1$ convergence in $A$. Second, for any $u\in A$ and any $0<\gamma_1<\gamma_2$, we have
$$\mathcal{R}_{\gamma_2}(u)\leq \mathcal{R}_{\gamma_1}(u)\leq \mathcal{R}(u)\quad\text{and}\quad \lim_{\gamma\to 0^+}\mathcal{R}_{\gamma}(u)= \mathcal{R}(u).$$
Let us denote, for consistency, $\mathcal{R}_0=\mathcal{R}$.

We shall use the following result.

\begin{lem}\label{R_epsilon}
For any $n\in\mathbb{N}$, let $u_n\in A$ and $\gamma_n\geq 0$ be such that, as $n\to\infty$, $u_n$ converges to $u\in A$ in $L^1$ and $\gamma_n$ is a nonincreasing sequence converging to $0$.
Then we have that
$$\lim_n\mathcal{R}_{\gamma_n}(u_n)=\mathcal{R}(u).$$
\end{lem}

\proof{.} Let us fix $k\in\mathbb{N}$. Then for any $n\geq k$ we have
$$\mathcal{R}_{\gamma_k}(u_n)\leq \mathcal{R}_{\gamma_n}(u_n)\leq \mathcal{R}(u_n).$$
By the semicontinuity of $\mathcal{R}_{\gamma_k}$ and the continuity of $\mathcal{R}$, we infer that for any $k\in\mathbb{N}$
$$\mathcal{R}_{\gamma_k}(u)\leq\liminf_n \mathcal{R}_{\gamma_k}(u_n)\leq
\liminf_n \mathcal{R}_{\gamma_n}(u_n)\leq
\limsup_n \mathcal{R}_{\gamma_n}(u_n)\leq
 \mathcal{R}(u).$$
Letting $k$ go to infinity we easily conclude the proof.\cvd

\bigskip

The main example of such a family of operators $\mathcal{R}_{\gamma}$ is given by substituting in the definition of  $\mathcal{R}$ in \eqref{newR}
 the function $f$ with a function $f_{\gamma}:\mathbb{R}\to [0,+\infty]$, that is to define
$$\mathcal{R}_{\gamma}(u)=\int_{B_{\tilde{R}}}f_{\gamma}\left(\|\nabla (\mathcal{I}(u))\|^2-\varphi(\mathcal{I}(u))\right),\quad\text{for any }u\in A.$$
In order to have the required properties on the family of operators $\mathcal{R}_{\gamma}$ it is enough that
$f_{\gamma}$ is continuous on $\mathbb{R}$ and, 
for any $x\in \mathbb{R}$ and any $0<\gamma_1<\gamma_2$, we have
$$f_{\gamma_2}(x)\leq f_{\gamma_1}(x)\leq f(x)\quad\text{and}\quad \lim_{\gamma\to 0^+}f_{\gamma}(x)= f(x).$$
The main advantage is that in this way we may choose $f_{\gamma}$ which is smooth and real valued everywhere.

We are now in the position of describing the approximating functionals and
proving the $\Gamma$-convergence result. Let us a fix a constant $p_1$, $1<p_1<+\infty$, and a continuous function $W:\mathbb{R}\to[0,+\infty)$ such that
$W(t)=0$ if and only if $t\in\{0,1\}$.
Let us denote by $\mathcal{P}_{\varepsilon}$, $\varepsilon>0$, the functional defined in
\eqref{modmordef} with $p=p_1$ and the double-well potential $W$. We recall that the functional $\mathcal{P}$
is defined in \eqref{Pdef}.

Then, for any $\varepsilon>0$, let us define
$F_{\varepsilon}:A\to [0,+\infty]$
such that
\begin{equation}\label{Fepsilondefin}
F_{\varepsilon}(u)=d_{st}(\Phi_{\eta(\varepsilon)}(u),\chi_{\Omega_0})+
b\mathcal{P}_{\varepsilon}(u)+c\mathcal{R}_{\gamma(\varepsilon)}(u),\quad\text{for any }u\in A,
\end{equation}
where $\eta:[0,+\infty)\to[0,+\infty)$ is a continuous, increasing function such that
$\eta(0)=0$ and $\gamma:[0,+\infty)\to[0,+\infty)$ is a continuous, nondecreasing function such that
$\gamma(0)=0$. Notice that here we may also assume $\gamma$ identically equal to zero.
Let us observe that, for any $u\in A$,
$$d_{st}(\Phi_{\eta(\varepsilon)}(u),\chi_{\Omega_0})=
\int_{B_{\tilde{R}}}|\Phi_{\eta(\varepsilon)}(u)-\chi_{\Omega_0}|+a\left|\int_{B_{\tilde{R}}}\|\nabla(\Phi_{\eta(\varepsilon)}(u))\|   -P(\Omega_0)\right|.$$

By the direct method, each of the functionals
$F_{\varepsilon}$, $\varepsilon>0$, admits a minimum over $A$.
We now state the $\Gamma$-convergence result.

\begin{teo}\label{gammaconvteo}
As $\varepsilon\to 0^+$,
$F_{\varepsilon}$ $\Gamma$-converges to $F_0$
on $A$ with respect to the $L^1$ norm.
\end{teo}

\proof{.}
Let us fix $\varepsilon_n>0$, $n\in\mathbb{N}$, such that $\varepsilon_n$ converges to $0$ as $n\to\infty$. Let $F_n=F_{\varepsilon_n}$ and let $u\in A$. Without loss of generality we may assume that $\varepsilon_n$ is decreasing with respect to $n\in\mathbb{N}$.


We begin with the $\Gamma$-$\liminf$ inequality. Let $u_n\in A$, $n\in\mathbb{N}$, be such that
$u_n$ converges to $u$ in $L^1$. Without loss of generality, we may assume that
$F_n(u_n)$, $n\in\mathbb{N}$, is uniformly bounded by a constant $C$. Then by Lemma~\ref{R_epsilon}, and by the continuity of $\mathcal{R}$,
we infer that $u$ and $u_n$, for any $n\in\mathbb{N}$ large enough, belong to $\tilde{A}_C$.
Then
$$d_{st}(\Phi_{\eta(\varepsilon_n)}(u_n),\mathcal{W}(u))\leq
d_{st}(\Phi_{\eta(\varepsilon_n)}(u_n),\mathcal{W}(u_n))+
d_{st}(\mathcal{W}(u_n),\mathcal{W}(u)).$$
As $n\to\infty$, the first term of the right-hand side converges to $0$ by Proposition~\ref{uniformcontprop}, whereas the second converges to $0$ by Lemma~\ref{firstlemma}.
Therefore the $\Gamma$-$\liminf$ inequality immediately follows  by the $\Gamma$-convergence of $\mathcal{P}_{\varepsilon}$ to $\mathcal{P}$ and Lemma~\ref{R_epsilon}.

For what concerns the recovery sequence,
without loss of generality we may restrict ourselves to the case in which $F_0(u)$ is finite.
The $\Gamma$-convergence of $\mathcal{P}_{\varepsilon}$ to $\mathcal{P}$ allows us to find $u_n\in A$, $n\in\mathbb{N}$, such that $u_n$ converges to $u$ in $L^1$ and $\mathcal{P}_{\varepsilon_n}(u_n)$ converges to $\mathcal{P}(u)$. By Lemma~\ref{R_epsilon}. and
a reasoning analogous to the one developed for the $\Gamma$-$\liminf$ inequality,
we immediately conclude that $F_n(u_n)$ converges to $F_0(u)$.\cvd

\begin{prop}\label{equicoerciveprop}
There exist $\tilde{\varepsilon}>0$ and a compact subset $\mathcal{K}$ of $A$ such that for any $\varepsilon$, $0<\varepsilon\leq\tilde{\varepsilon}$, we have
$$\min_{\mathcal{K}}F_{\varepsilon}=\min_{A}F_{\varepsilon}.$$
\end{prop}

\proof{.}
By the $\Gamma$-convergence result, in particular by the construction of the recovery sequence applied to a minimizer of $F_0$,
we infer that there exist $\tilde{\varepsilon}>0$ and a positive constant $C_1$ such that for any $\varepsilon$, $0<\varepsilon\leq\tilde{\varepsilon}$, we have
$\inf_AF_{\varepsilon}=\min_AF_{\varepsilon}\leq C_1$.
Let $u_{\varepsilon}\in A$, $0<\varepsilon\leq \tilde{\varepsilon}$,
be such that $F_{\varepsilon}(u_{\varepsilon})=\min_AF_{\varepsilon}$.
Then we observe that the set $\{u_{\varepsilon}\}_{0<\varepsilon\leq\tilde{\varepsilon}}$
satisfies the properties of Proposition~\ref{compactnessoss} for some constant $C$.
Therefore $\{u_{\varepsilon}\}_{0<\varepsilon\leq\tilde{\varepsilon}}$ is precompact in 
$L^1(\mathbb{R}^2)$ and the proof is concluded.\cvd

\bigskip

Using Theorem~\ref{gammaconvteo} and Proposition~\ref{equicoerciveprop}, we apply the Fundamental Theorem of $\Gamma$-convergence to conclude with the following result.

\begin{teo}\label{finalteo}
We have that $F_0$ admits a minimum over $A$ and 
$$
\min_{A} F_0=
\lim_{\varepsilon\to 0^+}\inf_{A}
F_{\varepsilon}=\lim_{\varepsilon\to 0^+}\min_{A}
F_{\varepsilon}.
$$

Let $\varepsilon_n$, $n\in \mathbb{N}$, be a sequence of positive numbers converging to $0$. For any $n\in\mathbb{N}$, let $F_n=F_{\varepsilon_n}$.
If
$\{u_n\}_{n=1}^{\infty}$ is a sequence contained in $A$ which
converges, as $n\to\infty$, to $u\in A$ in $L^1$ and
satisfies $\lim_n F_n(u_n)=\lim_n\inf_{A} F_n$, then $u$ is a minimizer
of $F_0$ on $A$, that is $u$
solves the minimization problem \eqref{min0}.
\end{teo}

We notice that
if
$\{u_n\}_{n=1}^{\infty}$ is a sequence contained in $A$ 
such that $\lim_n F_n(u_n)=\lim_n\inf_{A} F_n$, then, again by Proposition~\ref{compactnessoss},
we have that, up to a subsequence, 
$u_n$
converges, as $n\to\infty$, to a function $u\in A$ in $L^1$, where 
$u$ is a minimizer
of $F_0$ on $A$.

Finally, we point out the following remark that may be of use in the numerical computation of the minimizers.

\begin{oss}\label{distoss} All the results remain valid also with 
the following modifications.
Since we are dealing only with $BV(B_{\tilde{R}})$ functions whose values are between $0$ and $1$, we may replace, in the definition of $F_0$ and $F_{\varepsilon}$, $\varepsilon>0$, the distance $d_{st}$ with the following distance-like function $d_{st}^p$, for any $p$, $1\leq p<+\infty$, 
defined as follows
$$d_{st}^p(u,v)=\int_{B_{\tilde{R}}}|u-v|^p+a\big||Du|(B_{\tilde{R}})-|Dv|(B_{\tilde{R}}) \big|\quad\text{ for any }u,\, v\in BV(B_{\tilde{R}},[0,1]).$$
Moreover, we may allow in all cases the parameter $a$ to be equal to $0$.
\end{oss}

\section{The numerical experiments}\label{numsec}

For some positive $\varepsilon$ we shall minimize
the functional $F_{\varepsilon}$ defined in \eqref{Fepsilondefin}, with $d_{st}$ replaced by
$d_{st}^2$ (see Remark~\ref{distoss}) and where the tuning parameter $a$ is allowed to be $0$.

Let $\Omega$ be the computational domain, that for simplicity we choose as a square centered in the origin. We assume that 
$u$ is always a real valued function on $\Omega$ and is extended to zero outside $\Omega$.

Following \cite{Tuzel09}, we use the kernel
$K$ defined in \eqref{Kdefin} and approximate the mutual intensity function $J$ defined in \eqref{Jdef} and \eqref{Jdefstar} by
$$J_{approx}(x)=
\int_{\mathbb{R}^2}\frac{1}{\pi(k\sigma\text{NA})^2}\rme^{-\beta|\xi|^2}\rme^{\rmi \xi\cdot x}\rmd \xi,\quad x\in\mathbb{R}^2$$
where
$$\beta=\frac{\log 2}{(k\sigma\text{NA})^2}.$$

Correspondingly, the intensity function can be approximated as follows
$$\mathcal{I}_{approx}(u)(x)=\int_{\mathbb{R}^2}\int_{\mathbb{R}^2}u(x-\xi)H(\xi,\eta)u(x-\eta)
\rmd\xi\rmd\eta,\quad x\in\mathbb{R}^2$$
where
$$H(\xi,\eta)= K(\xi)J_{approx}(\eta-\xi)\overline{K}(\eta).$$
Here $H$ is the Hopkins transmission cross coefficients (TCC) function. 
The advantage of using the Hopkins model is that the TCC function is independent of the mask function $m$, therefore, given an optical system, the TCC function only needs to be computed once.

In the discrete formulation we subdivide $\Omega$ into $N\times N$ squares
with sides of length $
\Delta x=\Delta y$.
Setting $i=(i_1,i_2)$, $j=(j_1,j_2)$ and
$k=(k_1,k_2)$, with $i$, $j$, $k\in \{1,\ldots,N \}\times \{1,\ldots,N \}$, we have
$$\mathcal{I}_{discr}(u)(i)=\sum_{k}\sum_{j}\left(u(i-k)H(k,j)u(i-j)\right)\Delta x^2\Delta y^2$$
where $u$ is an $N\times N$ real valued matrix, whereas $H$ is here a $N^2\times N^2$ matrix. We notice that for any $N\times N$ matrix $M$ we shall always assume that
$M(i)=0$ for any $i\in \mathbb{Z}^2$ such that $i\not\in \{1,\ldots,N \}\times \{1,\ldots,N \}$. Analogously we are assuming that $H(k,j)=0$ for any $(k,j)\in \mathbb{Z}^2\times \mathbb{Z}^2$ such that either $k\not\in \{1,\ldots,N \}\times \{1,\ldots,N \}$
or $j\not\in \{1,\ldots,N \}\times \{1,\ldots,N \}$, that is 
we are introducing here another approximation since we are dropping $H$ to zero outside $\Omega\times \Omega$.

Following the arguments in \cite{Tuzel09}, we notice that $H$ is a positive, semi-definite Hermitian matrix and 
we decompose $H$ through a Singular Value Decomposition, actually through the following eigenfunction expansion
$$H(k,j)=\sum_n\sigma_n V_n(k)\overline{V_n}(j)$$
where $\sigma_n$ are the nonnegative eigenvalues, which we assume to be decreasing with respect to $n$, and $V_n$ are the corresponding eigenfunctions.
In general, due to the fast decay of its eigenvalues, we can approximate $H$ by its dominant eigenvectors, that is we may further approximate our intensity functional by using instead of $H$
$$H_{trunc}(k,j)=\sum_{n=1}^{N_0}\sigma_nV_n(k)\overline{V_n}(j)$$
for a suitably low number $N_0$. This leads to a further improvement on the computational efficiency of the aerial intensity. Moreover,
similar to the TCC function, the eigenvalues and eigenvectors are pre-computable when the optical system is fixed. 
Since $N_0$ is usually small, using $H_{trunc}$ instead of $H$ decreases the computational cost of repeated intensive calculations in the inverse problem.  

Finally, the discrete approximated intensity functional $I$ we shall use instead of $\mathcal{I}$ is
$$\frac{I(u)(i)}{\Delta x^2\Delta y^2}=
\sum_{n=1}^{N_0}\sigma_n
\sum_{k,j}\left(u(i-k)V_n(k)\overline{V_n}(j)u(i-j)\right)=\\
\sum_{n=1}^{N_0}\sigma_n(V_n\ast u)(i)\overline{(V_n\ast u)}(i)
$$
for any real valued $N\times N$ matrix $u$.
We notice that $\ast$ denotes the discrete convolution of matrices and we recall that any $N\times N$ matrix is extended to zero on any index $i\in \mathbb{Z}^2$ such that $i\not\in \{1,\ldots,N \}\times \{1,\ldots,N \}$.

In order to use a gradient method for our numerical computations we need to compute differentials with respect to $u$ of functionals applied to $I(u)$. For simplicity, we set for the time being $\Delta x=\Delta y=1$.
Given $u$, a real valued $N\times N$ matrix, and the variation $v$, another real valued $N\times N$ matrix, we have that
$$DI(u)[v]=2\Re\left(
\sum_{n=1}^{N_0}\sigma_n(V_n\ast u)\overline{(V_n\ast v)}
\right).
$$

Let us begin with the following simpler case.
Let
$F(u)=\int f(I(u))=\sum_k f(I(u)(k))\in \mathbb{R}$ for some real valued function $f$.
We need to compute
$\nabla_u F(u)\in \mathbb{R}^{N\times N}$. We remark that, if
$e_{i_1,i_2}$ is the matrix which is one in $(i_1,i_2)$ and zero elsewhere, then
$\nabla_u F(u)(i_1,i_2)=DF(u)[e_{i_1,i_2}]$. Using the same argument as before, we obtain
$$
DF(u)[e_{i_1,i_2}]=
2\Re \left(\sum_{n=1}^{N_0}\sigma_n\left(
\sum_{k}f'(I(u))(k) 
(V_n\ast u)(k)\overline{(V_n\ast e_{i_1,i_2})}(k)
\right)\right).
$$
But $(V_n\ast e_{i_1,i_2})(k)=V_n(k-i)$, with $i=(i_1,i_2)$.
Let us call $W_n(\cdot)=\overline{V_n(-\cdot)}$. Then summing upon $k$ we have
\begin{multline*}
\sum_{k}f'(I(u))(k) 
(V_n\ast u)(k)\overline{(V_n\ast e_{i_1,i_2})}(k)=\\
\sum_{k}f'(I(u))(k)
(V_n\ast u)(k)W_n(i-k)=
\left(W_n \ast \left(
f'(I(u))(V_n\ast u) 
\right)\right)(i).
\end{multline*}
In conclusion
$$\nabla_u F(u)=
2\Re\left(\sum_{n=1}^{N_0}\sigma_n
\left(W_n \ast \left(
f'(I(u))(V_n\ast u) 
\right)\right)
\right).$$
Let us notice that, assuming our computational domain is centered in $0$, we can easily compute $W_n$ through 
the MATLAB command
$W_n=\mbox{\tt conj}(\mbox{\tt flipud}( \mbox{\tt fliplr}(V_n)))$.

Let now
$$\tilde{F}(u)=\int \tilde{f}(\nabla f(I(u)))=
\sum_k \tilde{f}\left(\partial_{x_1}(f(I(u))) (k),\partial_{x_2}(f(I(u))) (k)
\right)
\in \mathbb{R},$$
where $\tilde{f}$ is a real valued function defined on $\mathbb{R}^2$ and $\partial_{x_l}$, $l=1,2$, are partial derivatives with respect to the two variables, in any finite differences sense.
Again we need to compute
$\nabla_u \tilde{F}(u)\in \mathbb{R}^{N\times N}$.
We denote $\tilde{f}_1$ the partial derivative of $\tilde{f}$ with respect to the first variable and  $\tilde{f}_2$ the partial derivative of $\tilde{f}$ with respect to the second variable.
Then similarly we have
$$
\nabla_u \tilde{F}(u)=\\
2\Re\left(\sum_{n=1}^{N_0}\sigma_n\left(W_n\ast A+((Z_n)_1\ast \overline{B_1})+((Z_n)_2\ast \overline{B_2})\right)
\right)
$$
where  $W_n(\cdot)=\overline{V_n(-\cdot)}$,
$(Z_n)_1(\cdot)=(\partial_{x_1}V_n)(-\cdot)$, $(Z_n)_2(\cdot)=(\partial_{x_2}V_n)(-\cdot)$ and
\begin{multline*}
A=
\sum_{l=1}^2\left(  \tilde{f}_l(\nabla f(I(u)))\partial_{x_l}(f'(I(u)))(V_n\ast u)+
\tilde{f}_l(\nabla f(I(u)))f'(I(u))
(\partial_{x_l}V_n\ast u)
\right),\\
B_1=\tilde{f}_1(\nabla f(I(u)))f'(I(u))
(V_n\ast u),\quad B_2=\tilde{f}_2(\nabla f(I(u)))f'(I(u))
(V_n\ast u).
\end{multline*}

Let us now investigate $\partial_{x_1}I(u)$. We have that
$$\partial_{x_1}I(u)=2\Re\left(\sum_{n=1}^{N_0}\sigma_n(\partial_{x_1}V_n\ast u)\overline{(V_n\ast u)}\right)
.$$
Therefore
$$
D(\partial_{x_1}I(u))[v]=
2\Re\left(\sum_{n=1}^{N_0}\sigma_n\left[(\partial_{x_1}V_n\ast u)\overline{(V_n\ast v)}+ (\partial_{x_1}V_n\ast v)\overline{(V_n\ast u)}\right]
\right).
$$

Finally, for a real valued function $g$ defined on $\mathbb{R}^2$, let
$$G(u)=\int g(\|\nabla_{x_1,x_2}I(u)\|^2,I(u))=
\sum_k g\left(\|\nabla_{x_1,x_2}I(u)\|^2(k),I(u)(k)\right)
\in \mathbb{R}.$$
We denote with $g_1$ the partial derivative of $g$ with respect to the first variable and  with $g_2$ the partial derivative of $g$ with respect to the second variable.
A completely similar computation leads to
$$\nabla_u G(u)=2\Re\left(\sum_n\sigma_n(A_1+A_2+\tilde{B}_1+\tilde{B}_2+C)
\right)
$$
where, for $l=1,2$,
$$A_l= W_n\ast (2g_1(\|\nabla_{x_1,x_2}I(u)\|^2,I(u))\partial_{x_l}I(u)(\partial_{x_l}V_n\ast u))
$$
$$\tilde{B}_l=(Z_n)_l\ast (2g_1(\|\nabla_{x_1,x_2}I(u)\|^2,I(u))\partial_{x_l}I(u)\overline{(V_n\ast u)})
$$
$$C=W_n\ast (g_2(\|\nabla_{x_1,x_2}I(u)\|^2,I(u))(V_n\ast u)).
$$

With these results we can compute the gradient for any term of our functional to be minimized. For example, the discretised version of our regularization functional $\mathcal{R}_{\gamma}$ may be expressed as
$$\mathcal{R}_{\gamma}(u)=\int f_{\gamma}\left(\|\nabla_{x_1,x_2}I(u)\|^2-\varphi(I(u))\right)$$
that is
$f_{\gamma}(\|\nabla_{x_1,x_2}I(u)\|^2-\varphi(I(u)))=g(\|\nabla_{x_1,x_2}I(u)\|^2,I(u))$
where
$g(a,b)=f_{\gamma}(a-\varphi(b))
$,
hence
$g_1(a,b)=f'_{\gamma}(a-\varphi(b))$ and $g_2(a,b)=-f'_{\gamma}(a-\varphi(b))\varphi'(b)$.

In our numerical experiments we shall use the following common parameters.
The computational domain is $1600 nm\times 1600 nm$ and, since we take $N=128$,
it is subdivided into $128 \times 128$ squares, each of them with sides of length $
\Delta x=\Delta y=12.5 nm$.

We use the first $N_0=10$ eigenvalues in $H_{trunc}$ and consider the optical system with the following physical parameters
$$\lambda=2\pi/k=193 nm,\quad NA=1,\quad\sigma=0.067.$$
These correspond to the parameters used in \cite{Tuzel09} even if our notation is slightly different. Actually in the numerical computation we use the eigenvalues $\sigma_n$ and eigenfunctions $V_n$ computed in \cite{Tuzel09}.

About the perimeter approximation $\mathcal{P}_{\varepsilon}$ defined in \eqref{modmordef} we choose $p=2$ and
$W(s)=s(1-s)$ for any $0\leq s\leq 1$. Since in our computation we are dropping in 
$\mathcal{P}_{\varepsilon}$ the constant $c_p$, in this section the parameter $b$ actually corresponds to $b/c_p$ in the notation of the previous section.


We choose $f_{\gamma}$ to be an approximation from below of $f$, where $f$ is defined as in \eqref{fdefinition} with $\alpha=1$ and $\varphi$ is as in \eqref{varphidefin}.
Normalizing the threshold $h$ to be equal to $1$ and the length of the pixel $\Delta x=\Delta y$ to be equal to $1$ as well, the term $\mathcal{R}_{\gamma}$ penalizes critical values of the intensity $I(u)$ close to $1$, namely the worst situation is when the triple
$(I(u),\partial_{x_1}I(u),\partial_{x_2}I(u))$ is equal to $(1,0,0)$. If we call
$d$ the distance between $(I(u),\partial_{x_1}I(u),\partial_{x_2}I(u))$ and $(1,0,0)$, in the Euclidean norm, that is
\begin{equation}\label{d_defin}
d=\sqrt{(I(u)-1)^2+(\partial_{x_1}I(u))^2+(\partial_{x_2}I(u))^2}
\quad(\text{normalized }h=1,  \Delta x=\Delta y=1),
\end{equation}
the aim of $\mathcal{R}_{\gamma}$ is not to let $d$ go to zero at any point, actually we wish to avoid the case in which $d$, in this normalized setting, is of the order of $5\%$ or less. Therefore we choose $\mathcal{R}_{\gamma}$ in such a way that it strongly penalizes the case in which $d$ is less than or equal to $5\%$ and has no effect whatsoever when $d$ is above $7 \%$. We keep this property fixed for any $\gamma$ and we let the values of $f_{\gamma}$, thus the value of $\mathcal{R}_{\gamma}$, increase as the positive parameter 
$\gamma$ goes to $0$.

We start with an initial value of $\varepsilon$, $\eta$ and $\gamma$,
namely $\varepsilon_0=0.002$, $\eta_0=0.2$, and $\gamma_0=0.03$, and
its corresponding functional $F_{\varepsilon_0}$, and a suitable
initial guess $u_{initial}$.  By a gradient method, namely a standard
steepest descent, we look for $u_0$, a minimizer of
$F_{\varepsilon_0}$, using $60$ iterations. Then we update the
parameters $\varepsilon$, $\eta$ and $\gamma$, by dividing their
previous values by the corresponding decrease rate given by 
$\mbox{\it rate}_{\varepsilon}$, $\mbox{\it rate}_{\eta}$ and 
$\mbox{\it rate}_{\gamma}$ respectively. We use the computed minimizer $u_0$ of
$F_{\varepsilon_0}$ as the initial guess and try to minimize the
functional $F_{\varepsilon}$ with the updated parameters. We repeat
the procedure after any $60$ iterations. This allows us to have at the
beginning a fast convergence to a reasonably good mask, no matter what
the initial guess is, and a refinement of the optimal mask later
on. The numerical experiments show that in general it is better to
keep these decrease rates rather close to $1$. After we have decreased
the parameters a fixed number of times, we consider the computed
minimizer of the last final functional $F_{\varepsilon}$ as our final
optimal phase-field function $u$.  The final optimal mask is obtained
from this numerical solution of the phase-field variable $u$ by taking
the set where $u>1/2$. Actually, as we shall show in our tests, due to
the presence of the Modica-Mortola functional, on most occasions the
final optimal phase-field function $u$ is already binary taking values
$0$ and $1$ only, therefore it coincides with the final optimal mask.

We recall that, given a phase-field function $u$, which is a function on the computational domain with values in $[0,1]$, 
its outcome pattern is the region where the light intensity is over the threshold value $h$.
The threshold in our tests equals $40 \%$ percent of the maximum value of $I_0$, where $I_0$ is the intensity when the mask is exactly the target pattern, that is, $h = 40\max(I_0)/100$.

\begin{figure}[htb]
\centering
\begin{minipage}[ht]{0.31\linewidth}
\includegraphics[width=1\textwidth]{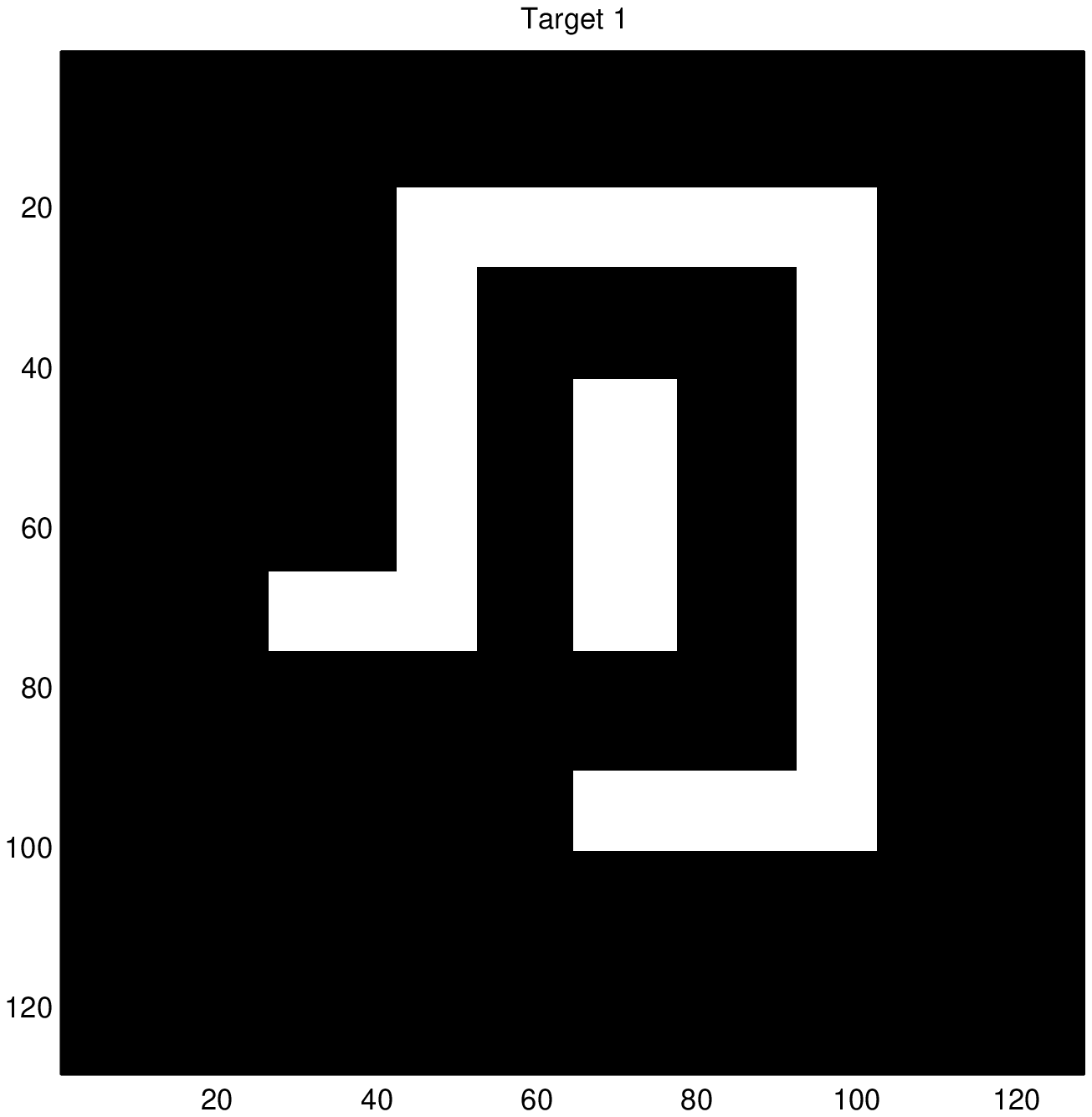}
\end{minipage}
\begin{minipage}[ht]{0.20\linewidth}\phantom{aaaaaaaaaaaaaa}
\end{minipage}
\begin{minipage}[ht]{0.31\linewidth}
\includegraphics[width=1\textwidth]{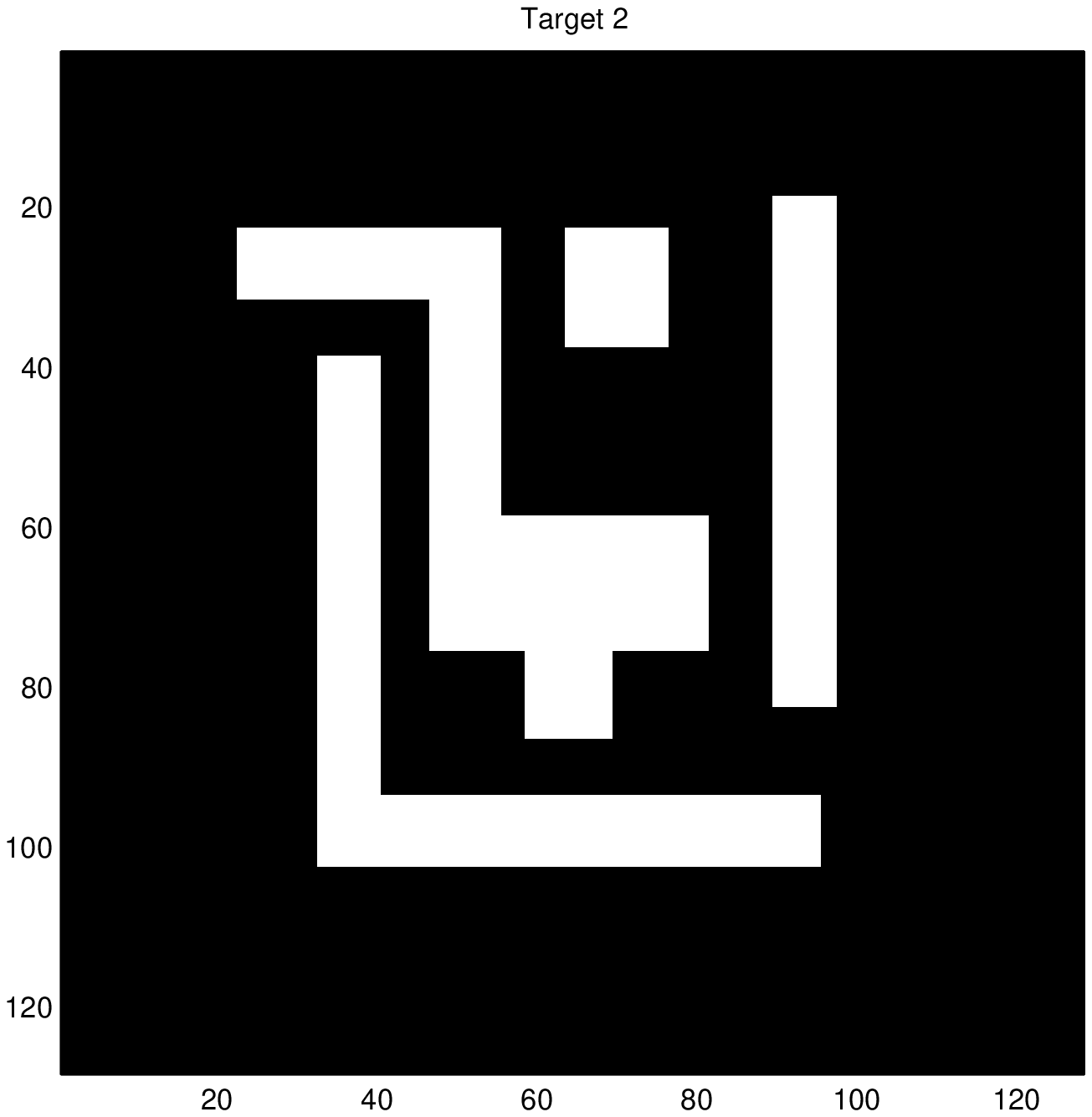}
\end{minipage}
\caption{Test target patterns. Left: Target 1. Right: Target 2.}\label{Fig:targets}
\end{figure}

We now describe the outcome of our numerical tests. We shall use two different types of targets, shown in Figure \ref{Fig:targets}. The first target pattern, Target 1,
is composed of two features.
The smallest width of the outside feature is 10 pixels, the width of the inside vertical bar is 13 pixels, two features are at least 12 pixels apart from each other.
The second target pattern, Target 2, is more complicated, consists of four features, with width as small as 8 pixels and distance between two different features as small as 6 pixels.

We first briefly discuss tests regarding Target $1$, then we move to the more interesting Target $2$.
\newcommand{\rate}{\mbox{\it rate}}
\subsubsection{Target 1}
In these tests we use the following parameters.
The weight of the term containing the difference between the two perimeters, or more precisely the two total variations, in the distance function $d_{st}^2$ is $a = 0$; the weight of the Modica-Mortola term $\mathcal{P}_{\epsilon}$ is $b = 2\times 10^{-4}$. The weight of the regularization term $\mathcal{R}_{\gamma}$ is $c = 0$.
Moreover we set $\rate_{\varepsilon}=1.2$, $\rate_{\eta}=1.2$ and
$\rate_{\gamma}=1.05$ and we perform $1080$ iterations in total, that is we decrease $17$ times our parameters. Correspondingly, at the end we compute the minimizer of the final functional $F_{\varepsilon}$ corresponding to the parameters
$\varepsilon=\varepsilon_0\times \rate_{\varepsilon}^{-17}\approx 9\times 10^{-5}$,
$\eta=\eta_0\times \rate_{\eta}^{-17}\approx 9\times 10^{-3}$,
and
$\gamma=\gamma_0\times \rate_{\gamma}^{-17}\approx 1.3\times 10^{-2}$.

\begin{figure}[htb]
\centering
\begin{minipage}[ht]{0.31\linewidth}
\includegraphics[width=1\textwidth]{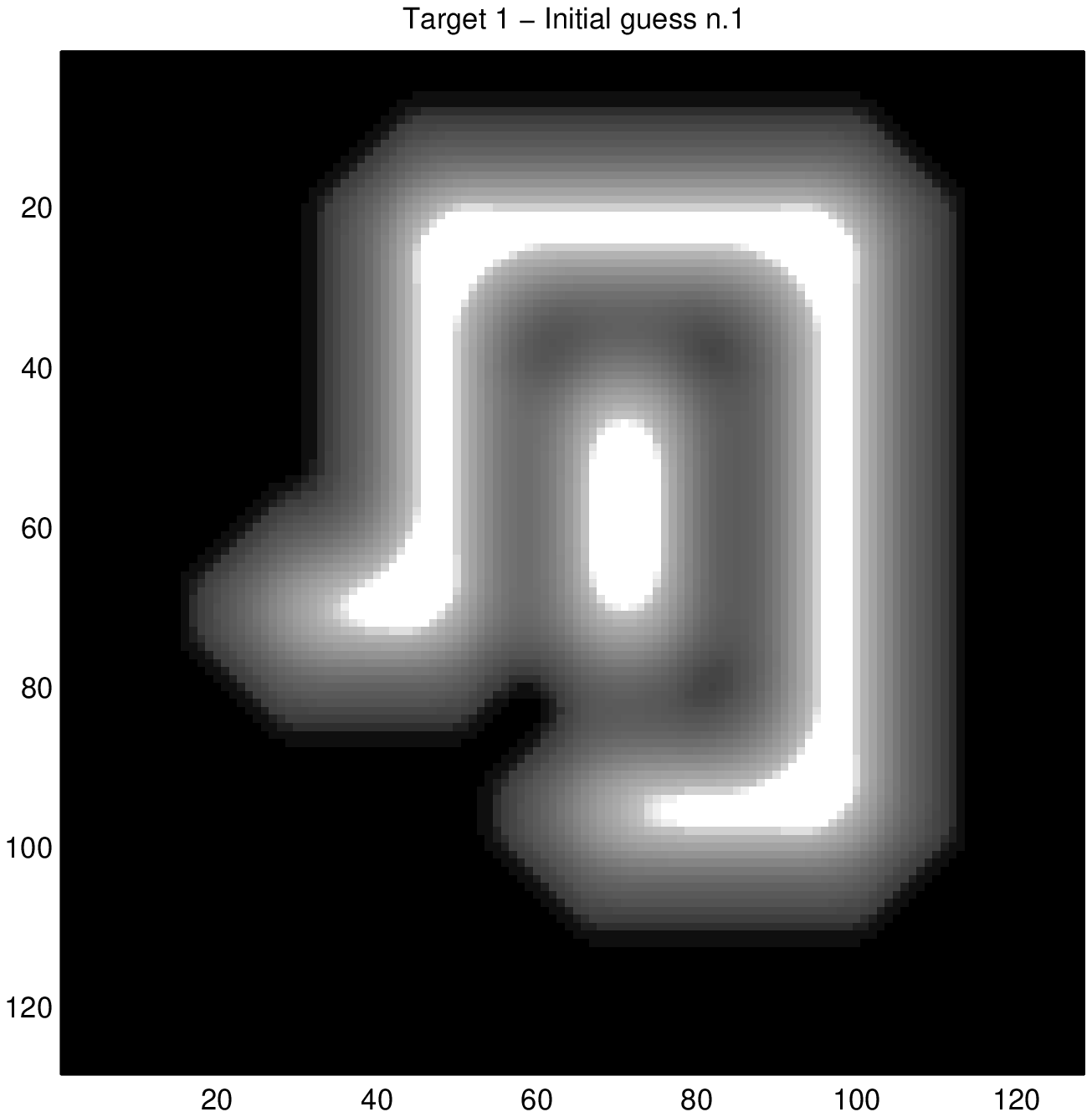}

\includegraphics[width=1\textwidth]{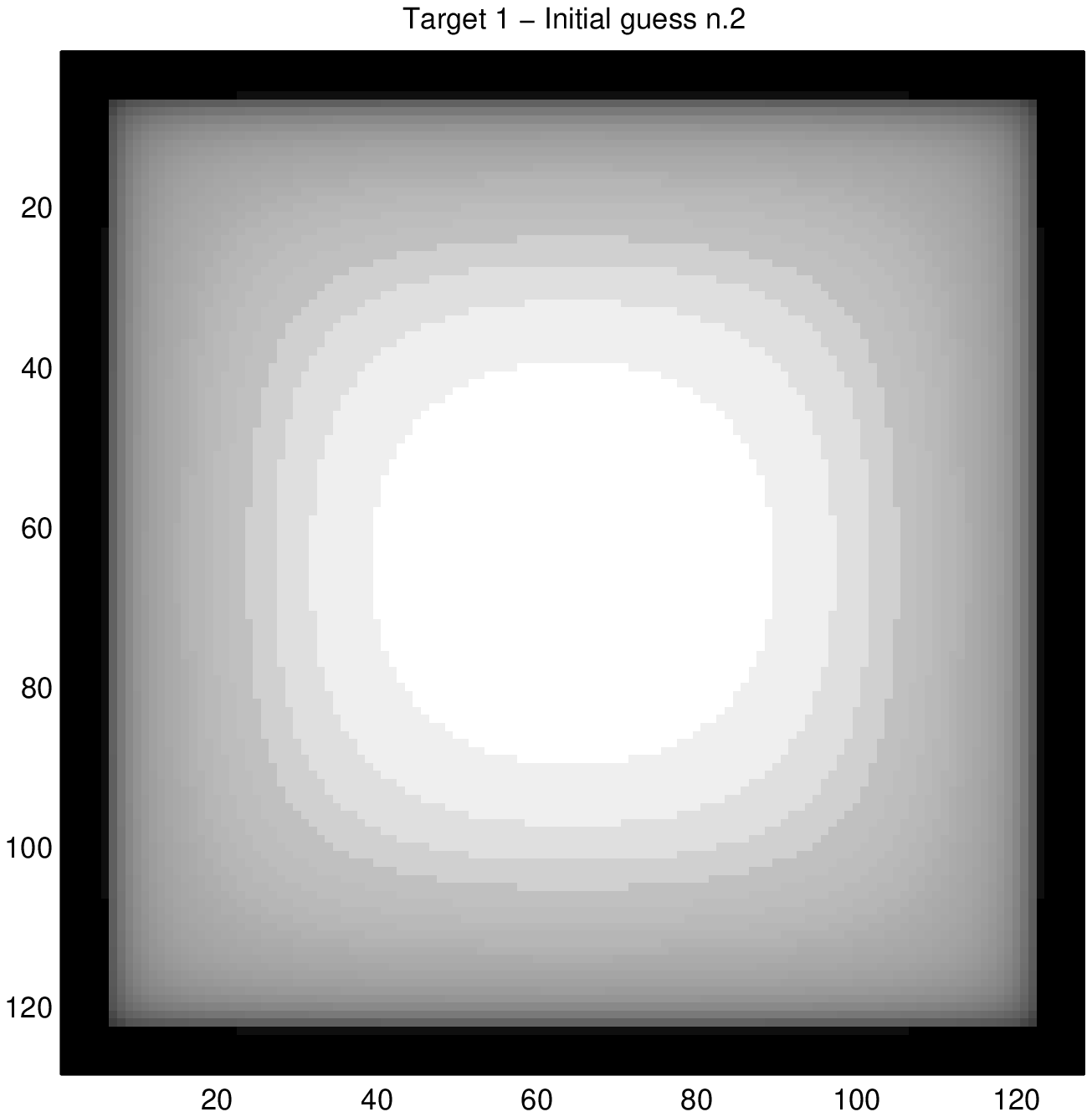}
\end{minipage}
\begin{minipage}[ht]{0.31\linewidth}
\includegraphics[width=1\textwidth]{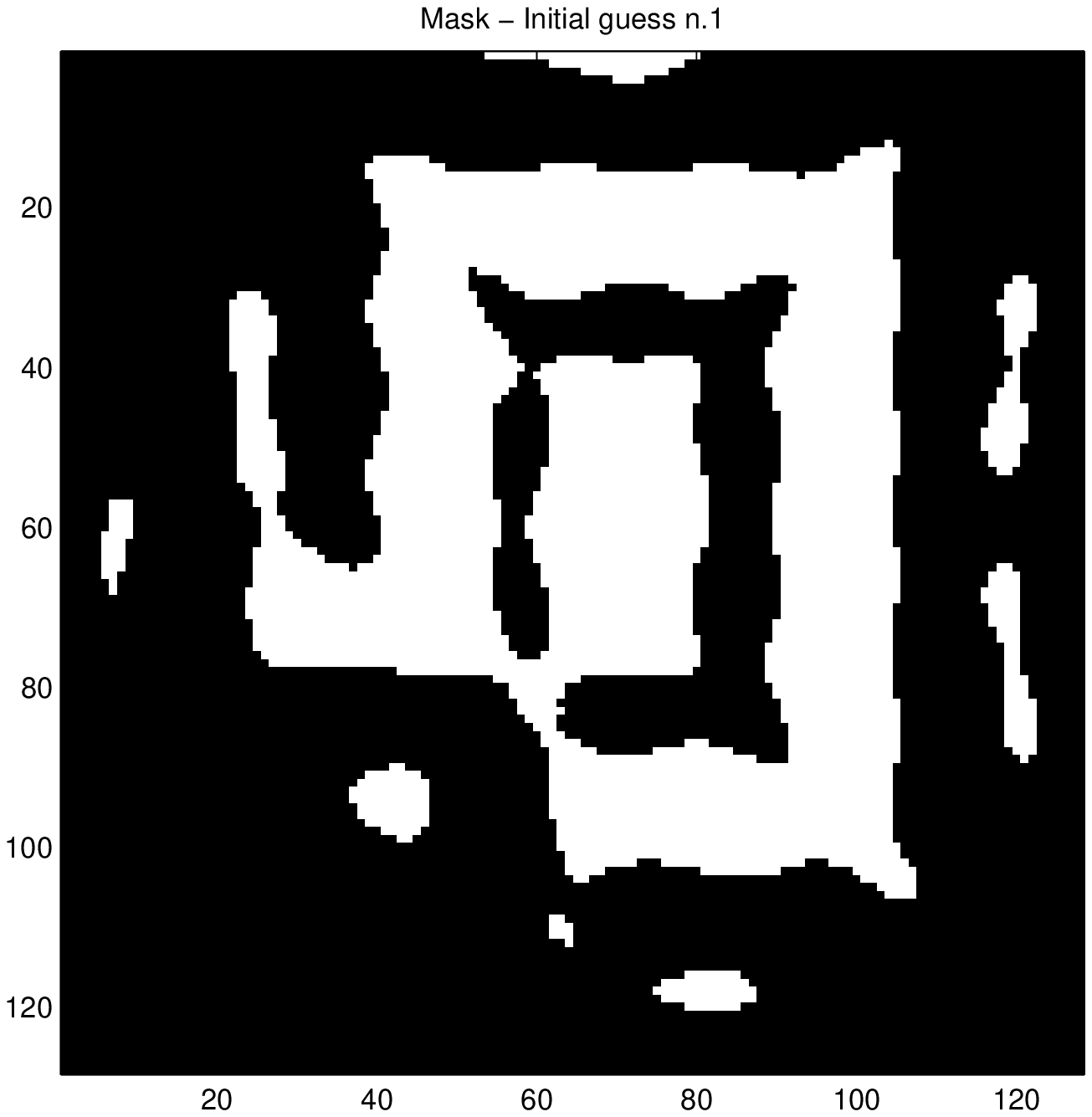}

\includegraphics[width=1\textwidth]{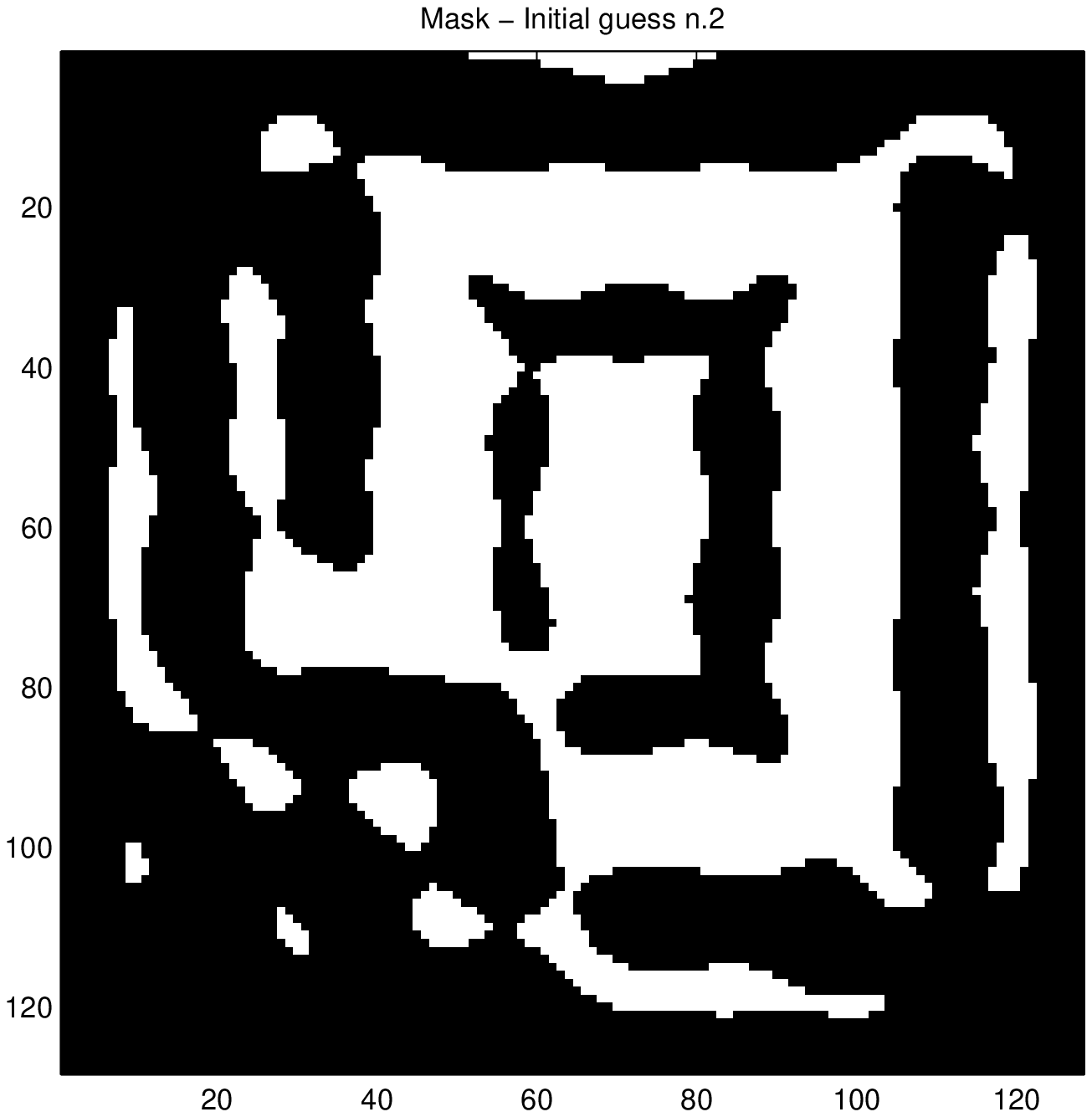}
\end{minipage}
\begin{minipage}[ht]{0.31\linewidth}
\includegraphics[width=1\textwidth]{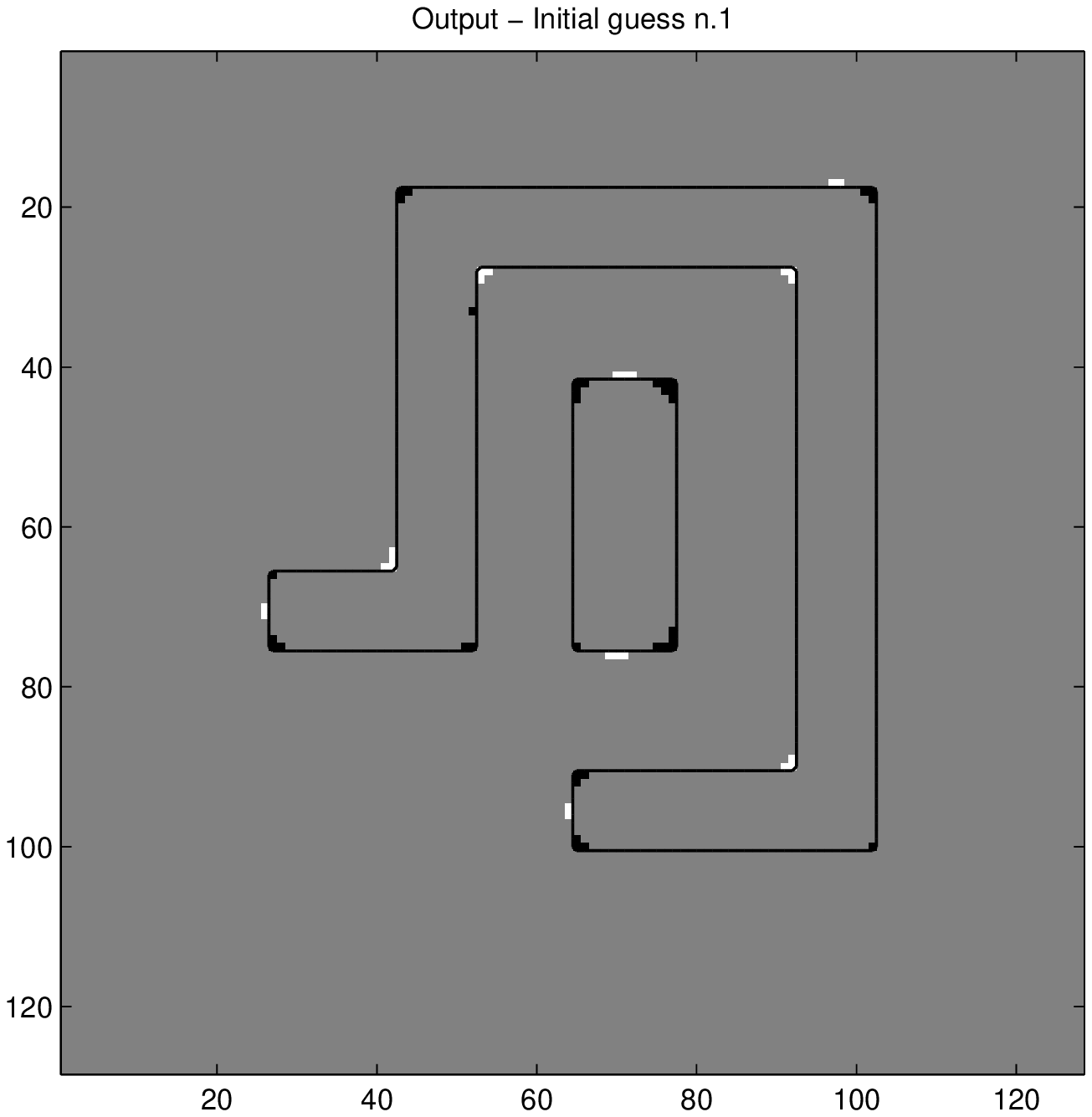}

\includegraphics[width=1\textwidth]{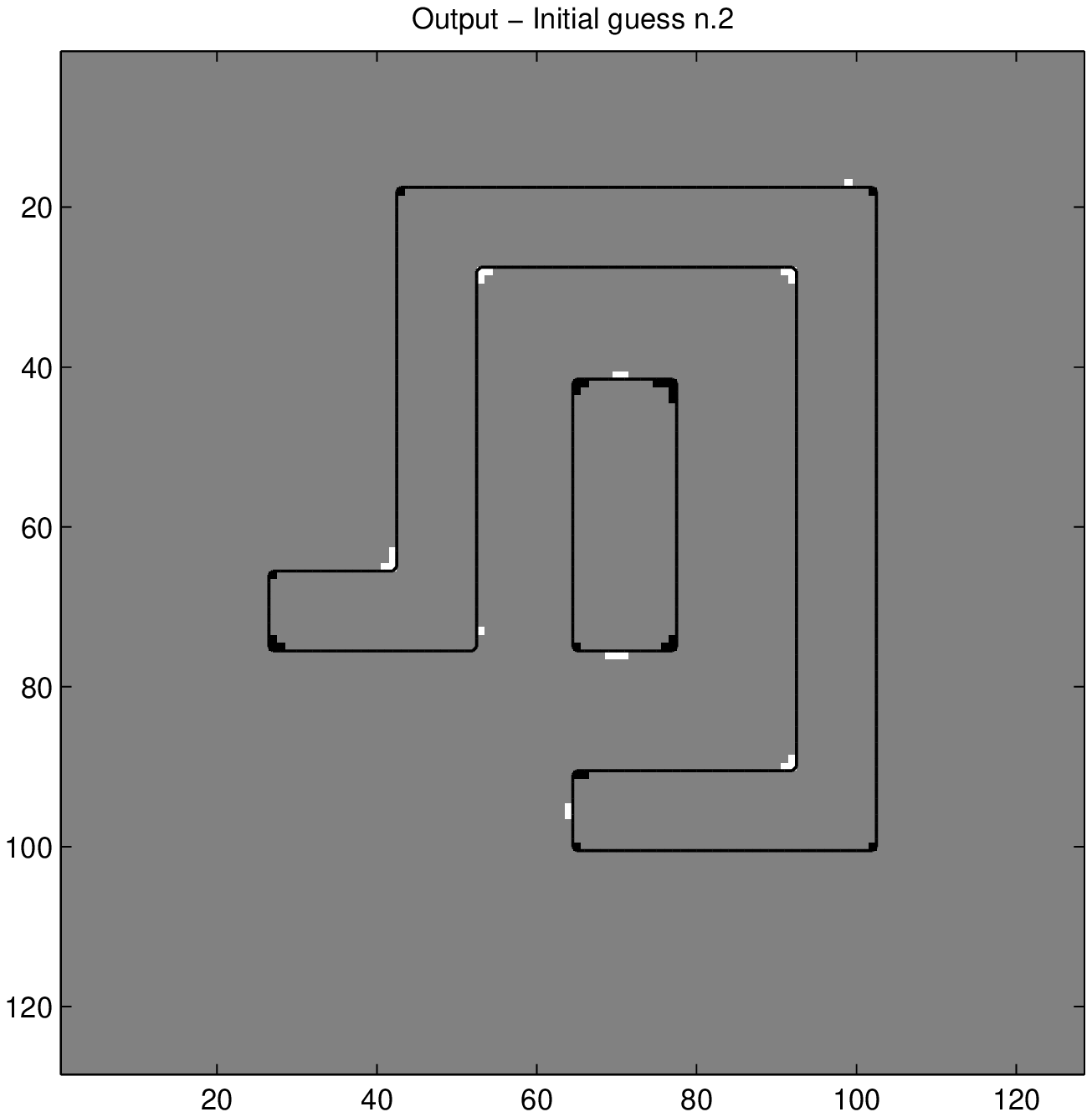}
\end{minipage}
\caption{Top: {\it Test n.1}. Left: initial guess {\it n.1}. Middle: mask. Right: output 
Bottom: {\it Test n.2}. Left: initial guess {\it n.2}. Middle: mask. Right: output.
 }\label{Fig:target1}
\end{figure}

We use two different initial guesses. In the {\it Test n.1} we consider an
initial guess which is a smooth perturbation of the target itself, in
{\it Test n.2} the initial guess is much more diffuse and has nothing to do
with the target itself. The results are presented in
Figure~\ref{Fig:target1}. Let us notice that for initial guesses and masks, the value $0$ is depicted in black, whereas the value $1$ is in white. Concerning the output, we show the difference between the exposed pattern and the target. Namely, in white we have the part of the exposed pattern which is outside the target and in black the part of the target that is not contained in the exposed pattern. The black line is the profile of the target.

First of all we have that in both cases we
converge to a binary function, due to the effect of the Modica-Mortola
functional. The mask so obtained is very diffuse, even with an initial
guess which is not. Actually, the reconstruction is better when the
initial guess is more diffuse. In fact the difference between the
exposed pattern and the target pattern is 61 pixels in {\it Test
n.1} and 44 pixels in {\it Test n.2} and the output is also visibly better.
We also notice that the two masks are rather different in shape, this
may be due to the fact that the original functional $F_0$ may have
several local minima and different initial guesses or different
choices of the parameters may therefore lead to quite different masks.

Since, in both cases, the intensity corresponding to the phase-fields during the iterations has never reached a critical point with value near to the threshold, the result does not change even if we add the regularization term $\mathcal{R}_{\gamma}$ (we have tested it with its coefficient $c$ varying from to $5\times 10^{-4}$ to $2\times 10^{-3}$), in accord to the theory.

\begin{figure}[htb]
\centering
\begin{minipage}[ht]{0.31\linewidth}
\includegraphics[width=1\textwidth]{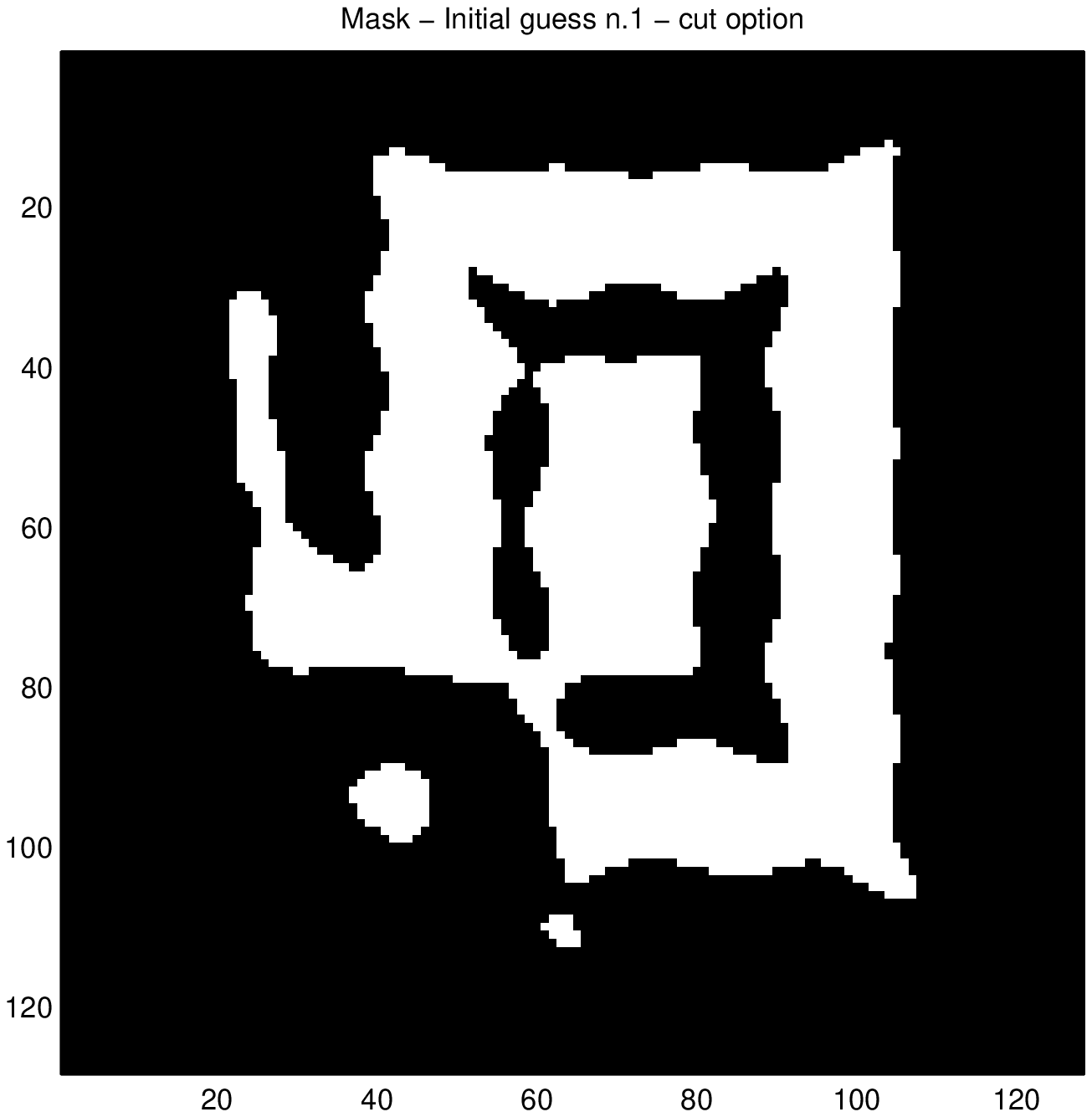}
\end{minipage}
\begin{minipage}[ht]{0.20\linewidth}\phantom{aaaaaaaaaaaaaa}
\end{minipage}
\begin{minipage}[ht]{0.31\linewidth}
\includegraphics[width=1\textwidth]{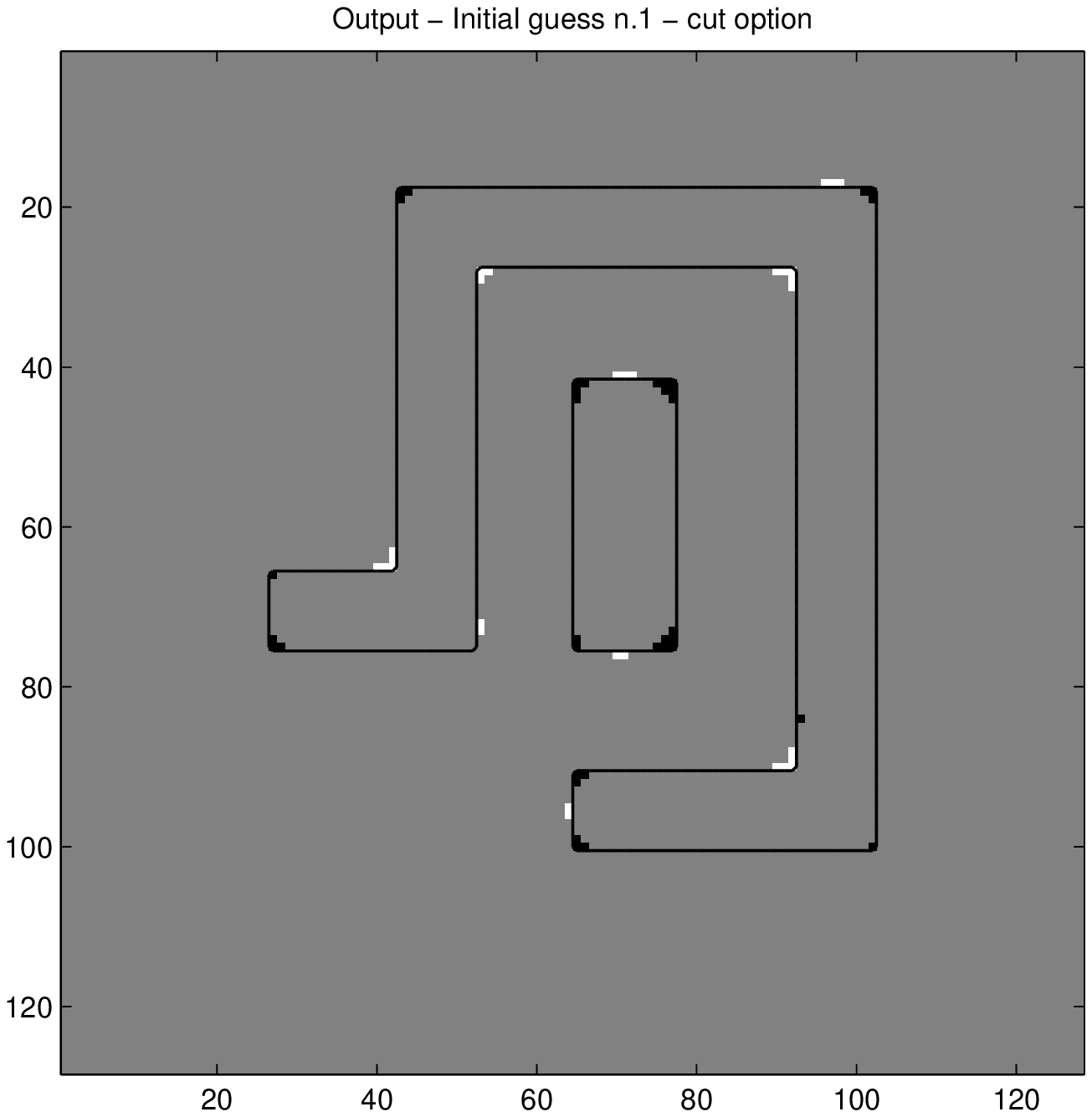}
\end{minipage}
\caption{{\it Test n.3} (cut option). Left: mask. Right: output.}\label{Fig:target1cut}
\end{figure}

In order to verify that having a diffuse mask with lot of assist features is an advantage, we took the initial guess of {\it Test n.1} but we impose our phase-fields during our iterations (and consequently our final mask) to be kept to zero outside a fixed neighbourhood of the target. The outcome is worse, the difference in pixels from the target being 66, see Figure~\ref{Fig:target1cut}.

\subsubsection{Target 2}
In the first two tests we use the same parameters as for Target $1$, namely
the weight of the difference between the perimeters in the distance function is $a = 0$; the weight of the Modica-Mortola term $\mathcal{P}_{\epsilon}$ is $b = 2\times 10^{-4}$. The weight of the regularization term $\mathcal{R}_{\gamma}$ is $c = 0$.
Moreover we set $rate_{\varepsilon}=1.2$, $rate_{\eta}=1.2$ and $rate_{\gamma}=1.05$ and we perform $1080$ iterations in total.

\begin{figure}[htb]
\centering
\begin{minipage}[ht]{0.31\linewidth}
\includegraphics[width=1\textwidth]{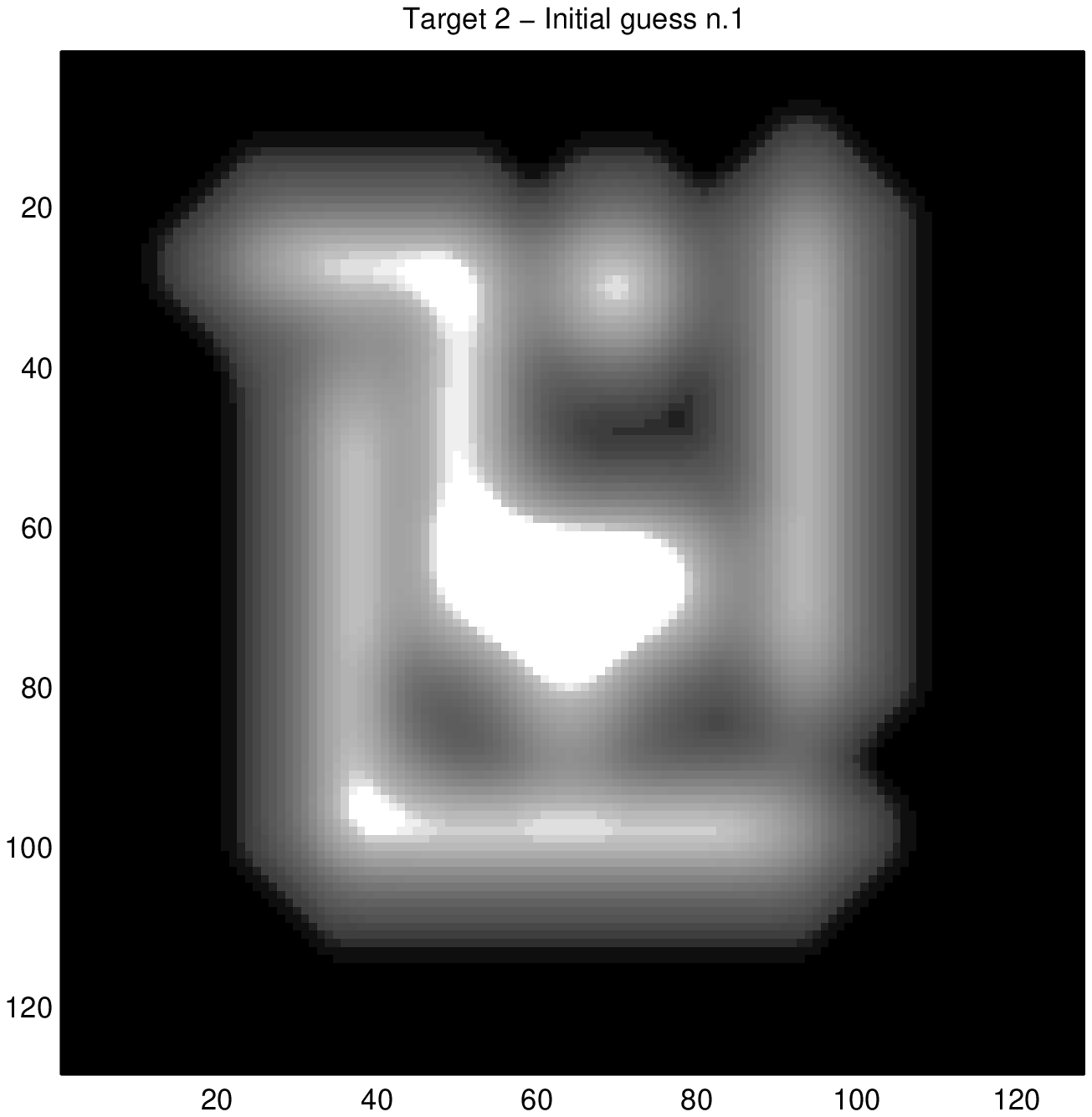}

\includegraphics[width=1\textwidth]{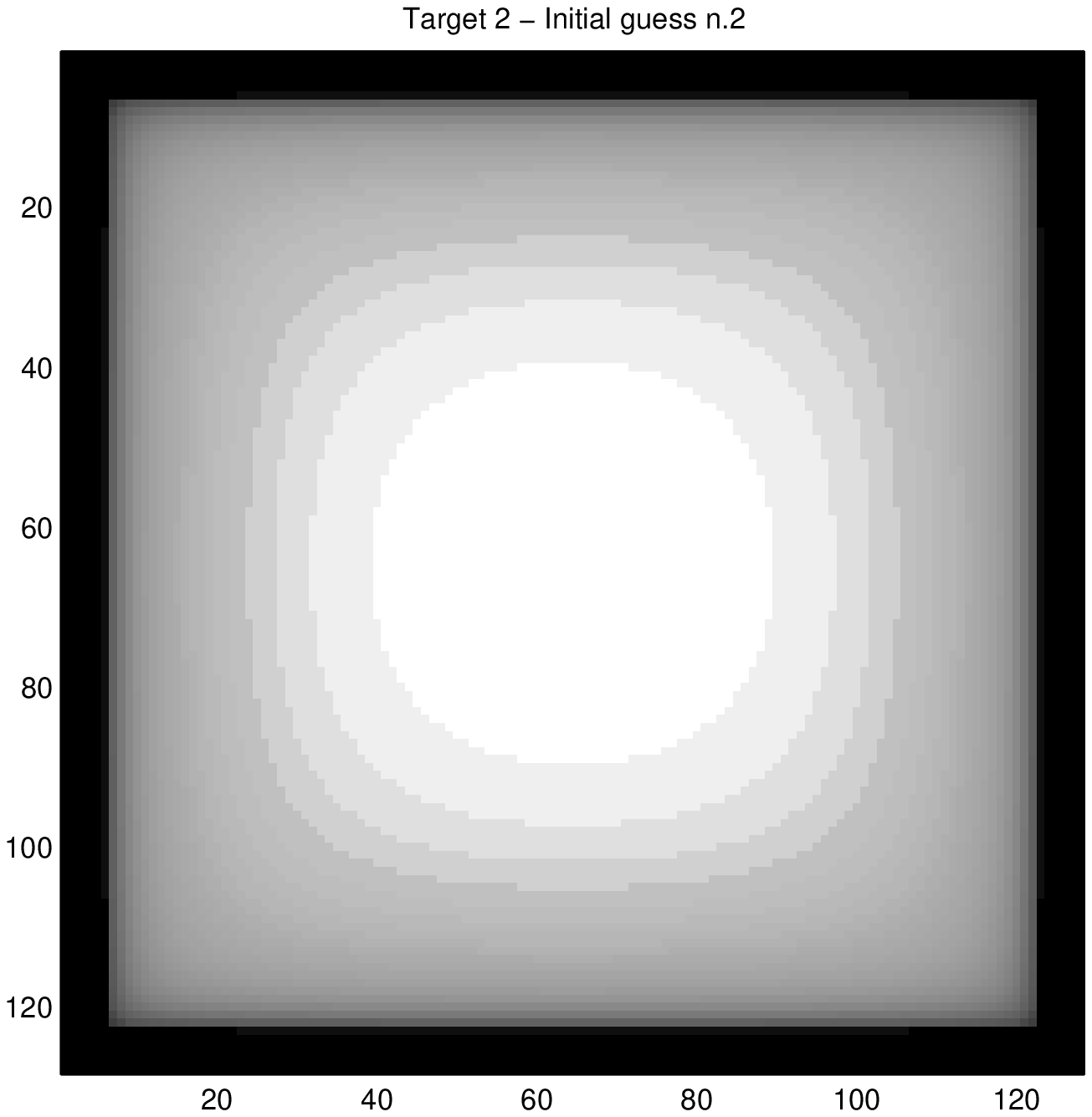}
\end{minipage}
\begin{minipage}[ht]{0.31\linewidth}
\includegraphics[width=1\textwidth]{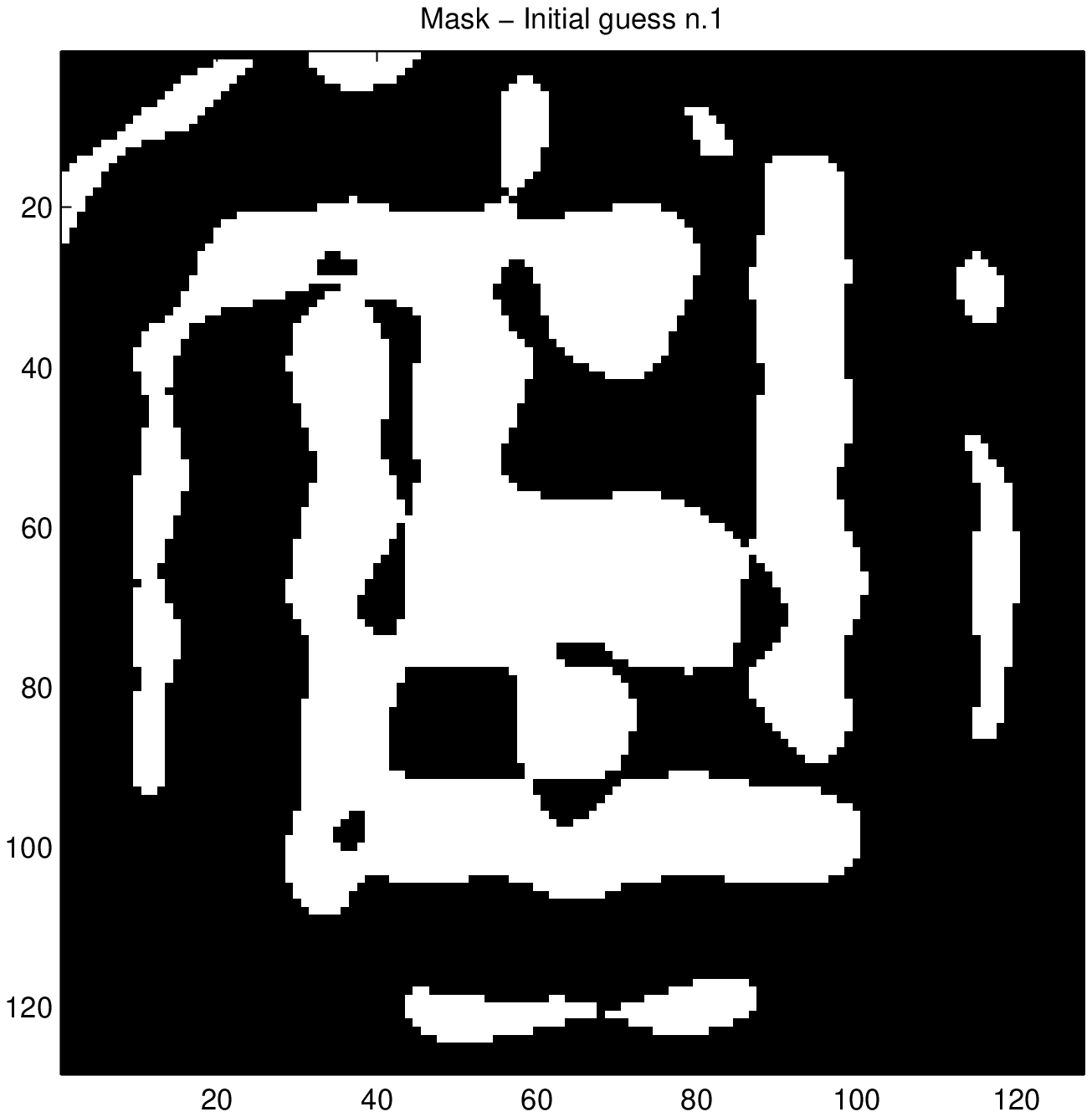}

\includegraphics[width=1\textwidth]{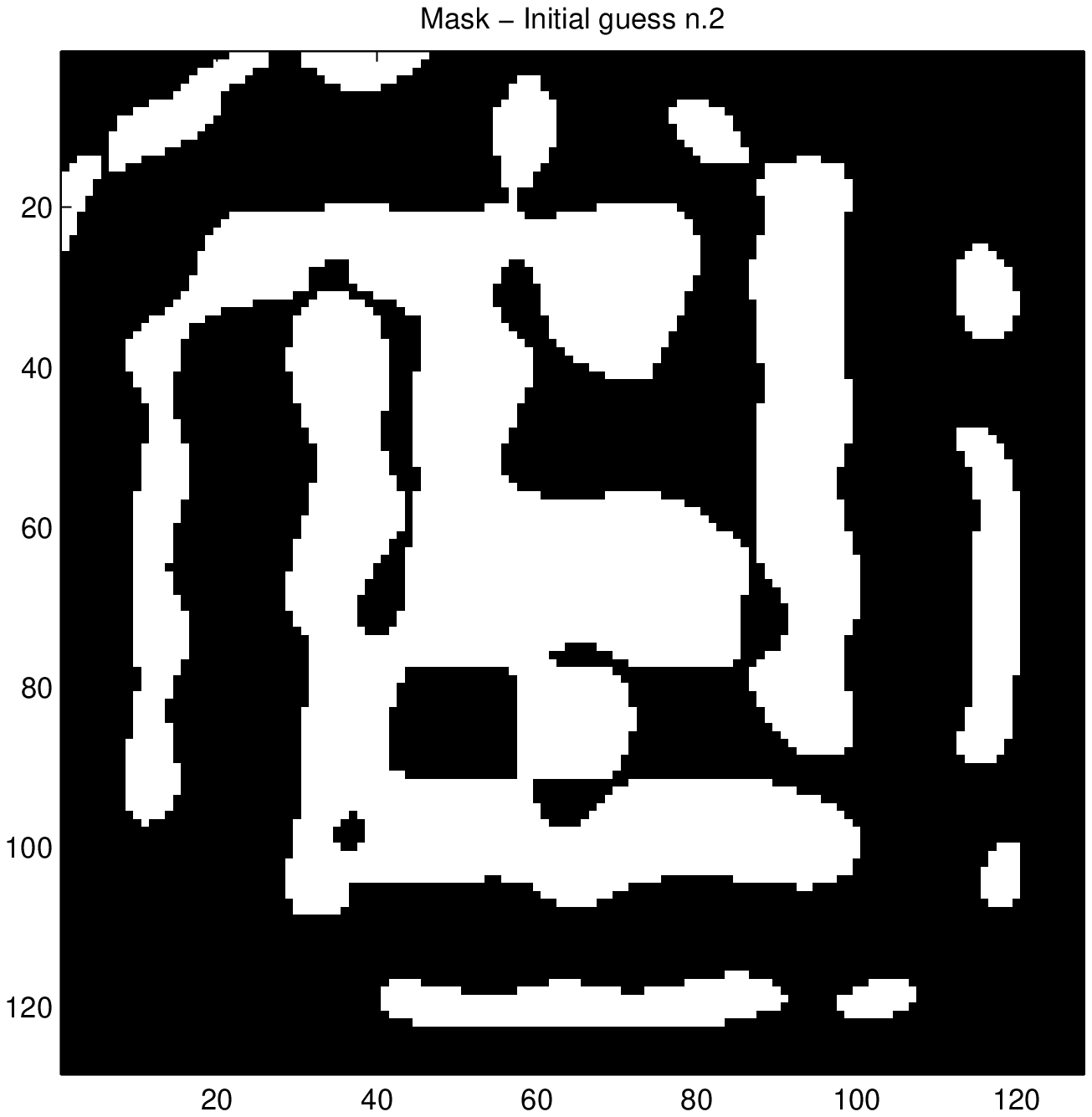}
\end{minipage}
\begin{minipage}[ht]{0.31\linewidth}
\includegraphics[width=1\textwidth]{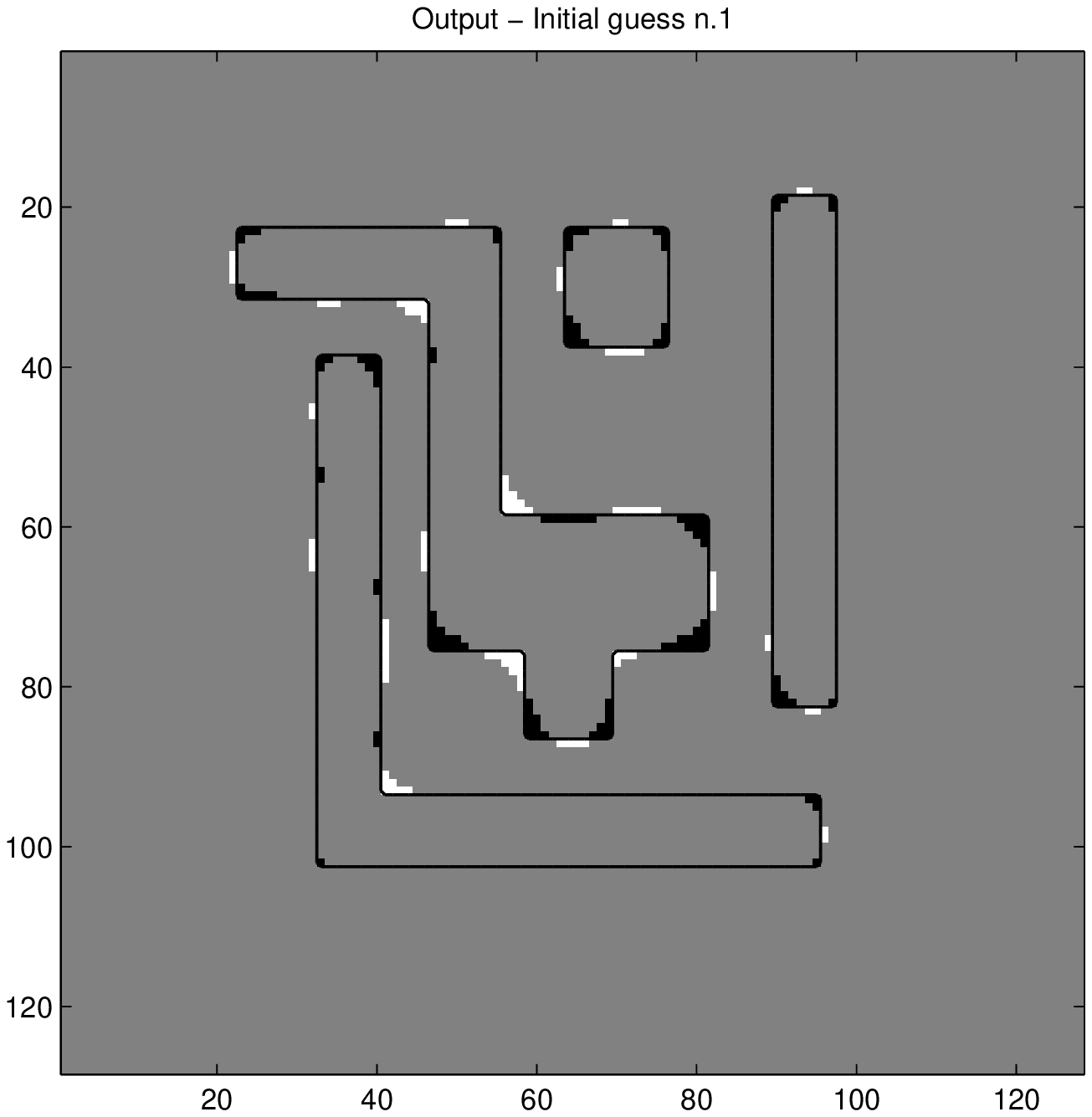}

\includegraphics[width=1\textwidth]{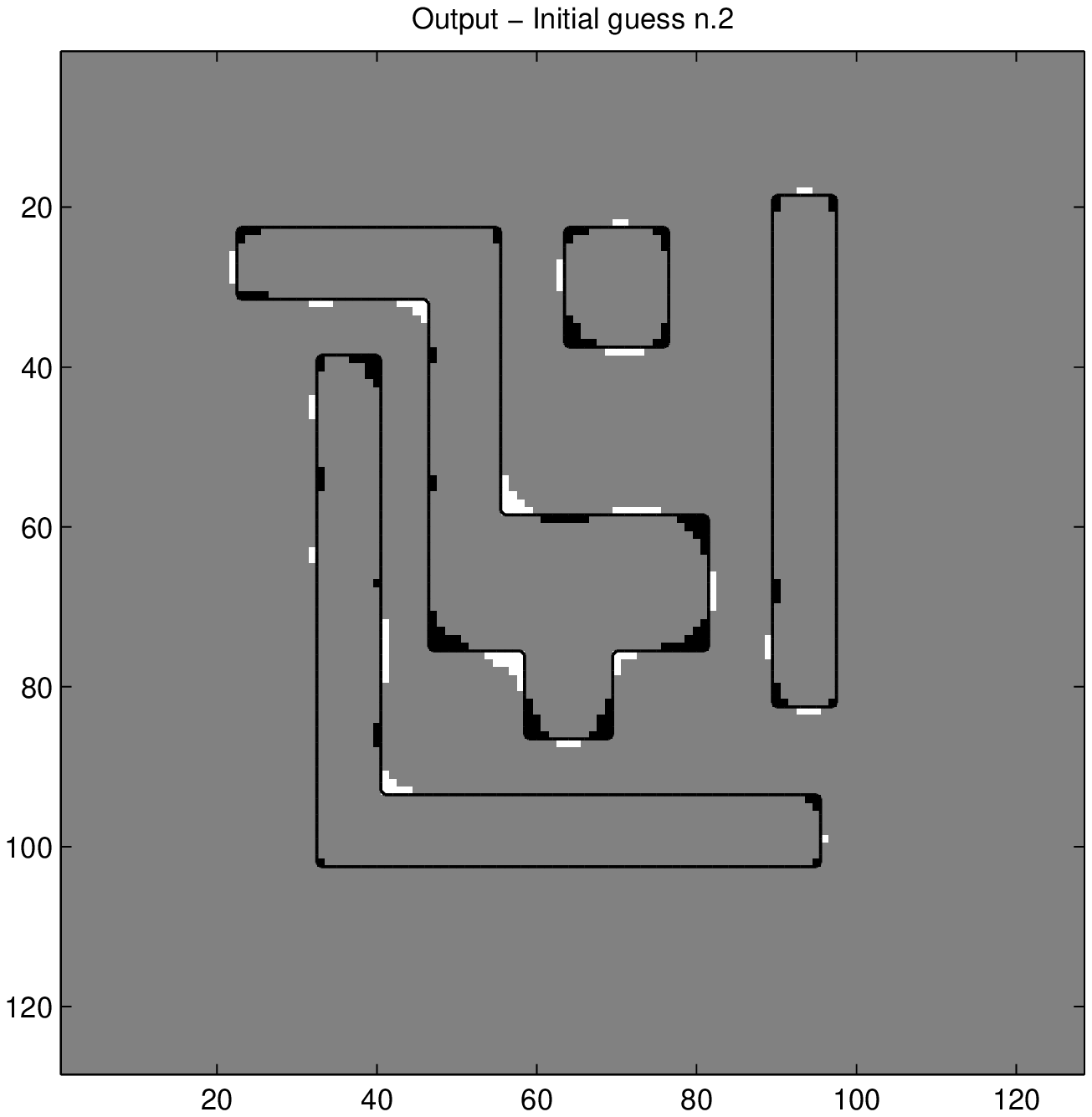}
\end{minipage}
\caption{Top: {\it Test n.1}. Left: initial guess {\it n.1}. Middle: mask. Right: output.
Bottom: {\it Test n.2}. Left: initial guess {\it n.2}. Middle: mask. Right: output.
 }
 \label{Fig:target2}
\end{figure}

We first investigate two tests with different initial guesses. In the {\it Test n.1} we consider an initial guess which is a smooth perturbation of the target itself, in {\it Test n.2} the initial guess is much more diffuse and has nothing to do with the target, it is actually the same as in {\it Test n.2} for Target 1. The results are presented in Figure~\ref{Fig:target2} and the conclusions are similar to those discussed for Target 1. Notice that the difference between exposed pattern and target is 233 for {\it Test n.1} and 227 for {\it Test n.2}. Hence,
we use the diffuse initial guess of {\it Test n.2} in all the following tests.

We shall discuss in detail the effect of the regularization term $R_{\gamma}$, the main theoretical novelty of the paper. Since it penalizes critical points at values close to the threshold value $h$, its effect should be the one to make the reconstruction more stable with respect to perturbations of $h$, especially from a topological point of view.

We consider the following two cases. In the first case we keep the parameters of {\it Test n.2} except the value of the coefficient $c$ of $\mathcal{R}_{\gamma}$. Namely, {\it Test n.3} is equal to {\it Test n.2} ($c=0$) whereas for {\it Test n.4} we set $c=5\times 10^{-4}$ and for {\it {\it Test n.5}} we set $c=2\times 10^{-3}$, that is we steadily increase the coefficient of $\mathcal{R}_{\gamma}$.

%
%
%

Notice that here sometimes the final optimal phase-field function $u$
is not binary, however the number of pixels where $u$ is different
from $0$ and $1$ is very limited. We conjecture that when this happens we
are most likely stuck near a local minimum of the final functional
$F_{\varepsilon}$.

\begin{figure}[htb]
\centering
\begin{minipage}[ht]{0.31\linewidth}
\includegraphics[width=1\textwidth]{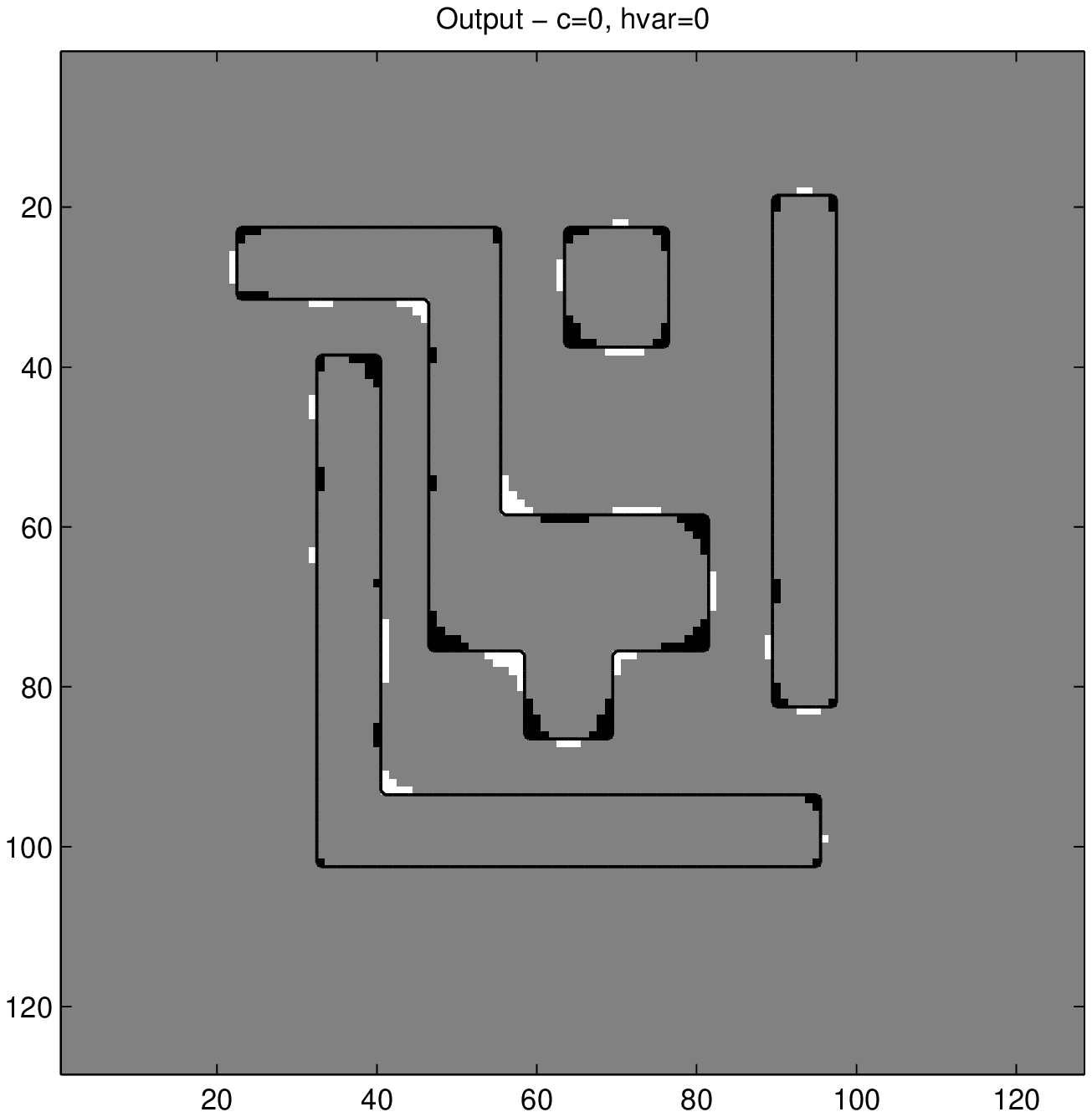}

\includegraphics[width=1\textwidth]{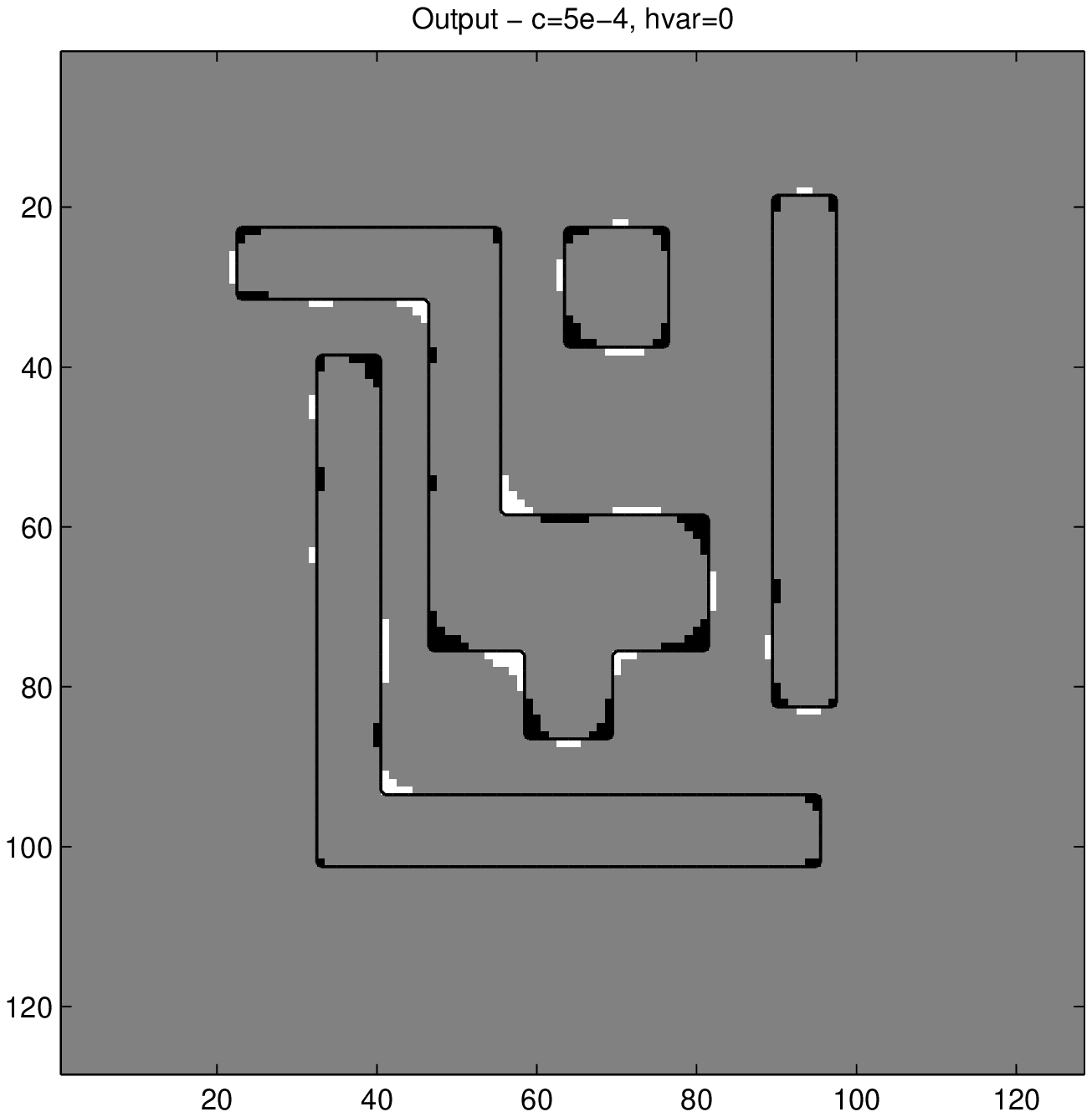}

\includegraphics[width=1\textwidth]{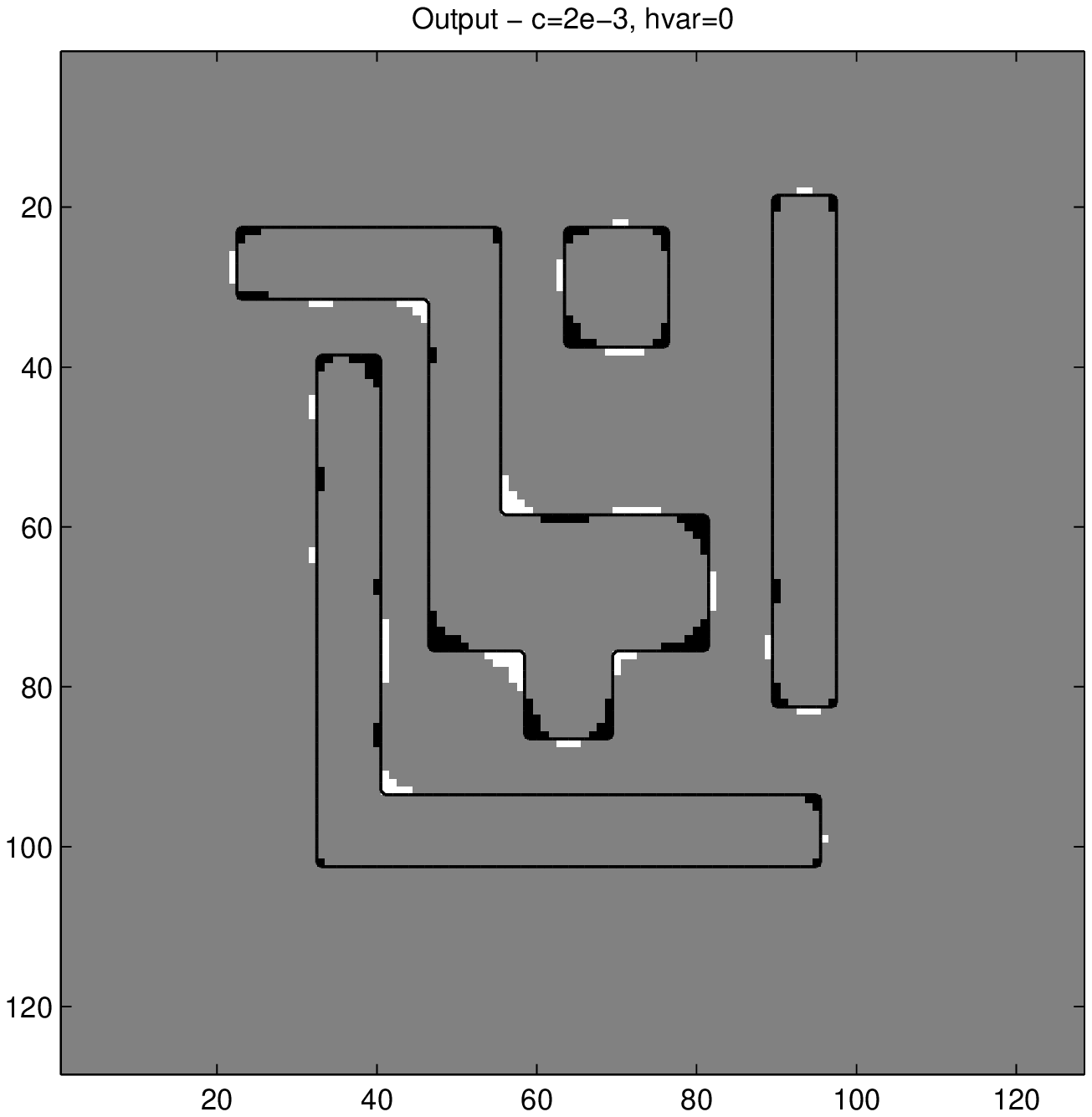}
\end{minipage}
\begin{minipage}[ht]{0.31\linewidth}
\includegraphics[width=1\textwidth]{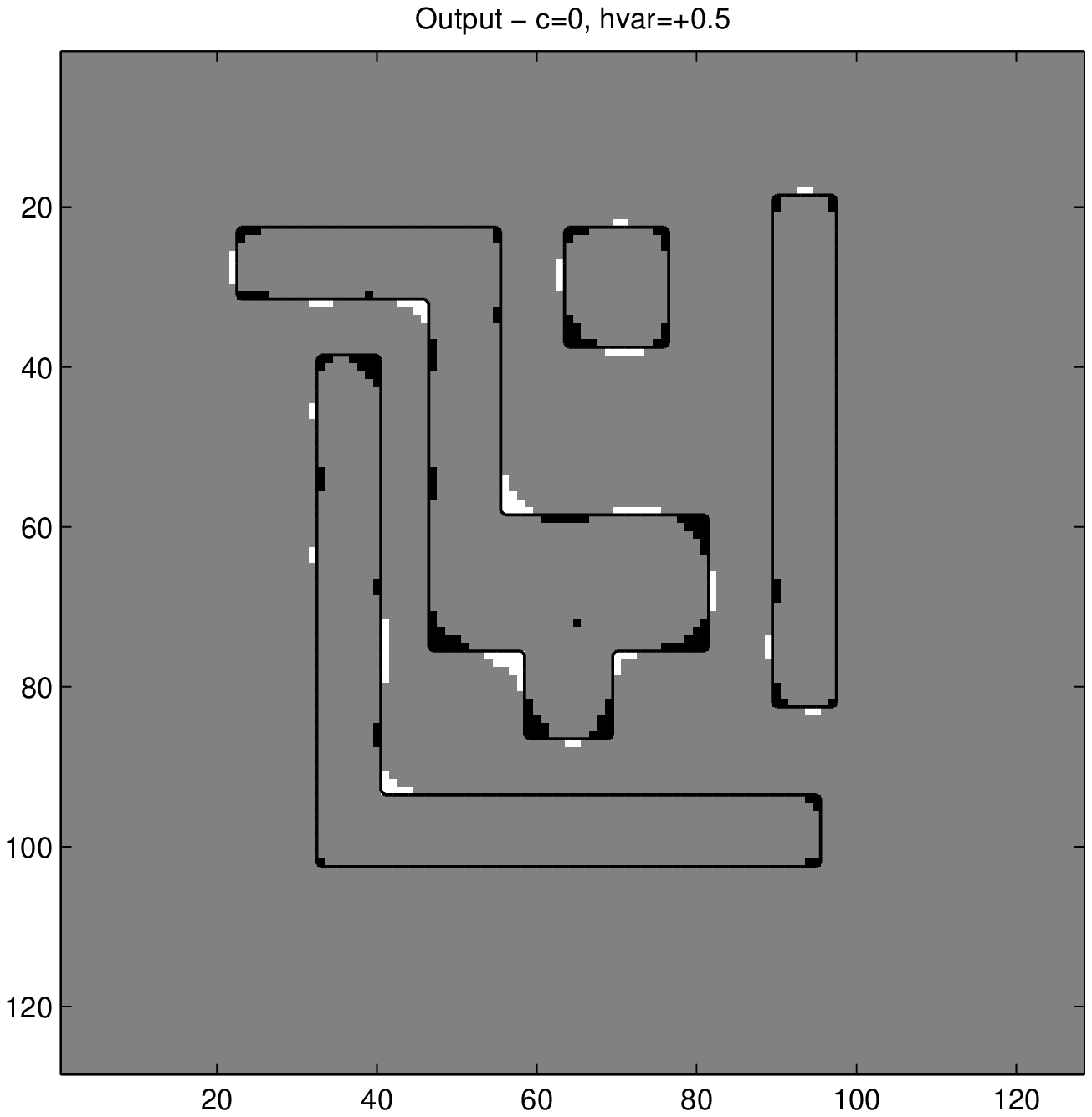}

\includegraphics[width=1\textwidth]{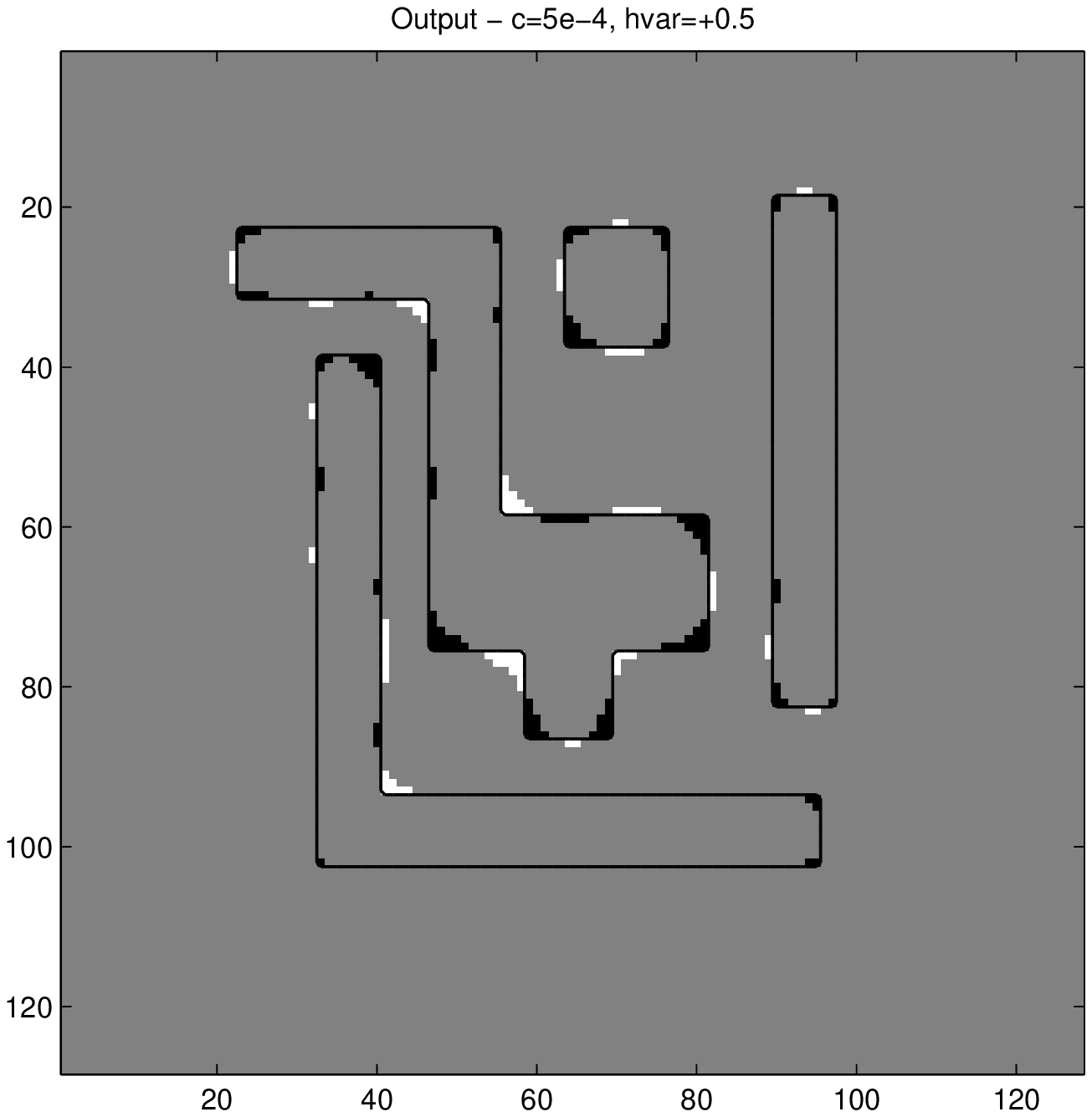}

\includegraphics[width=1\textwidth]{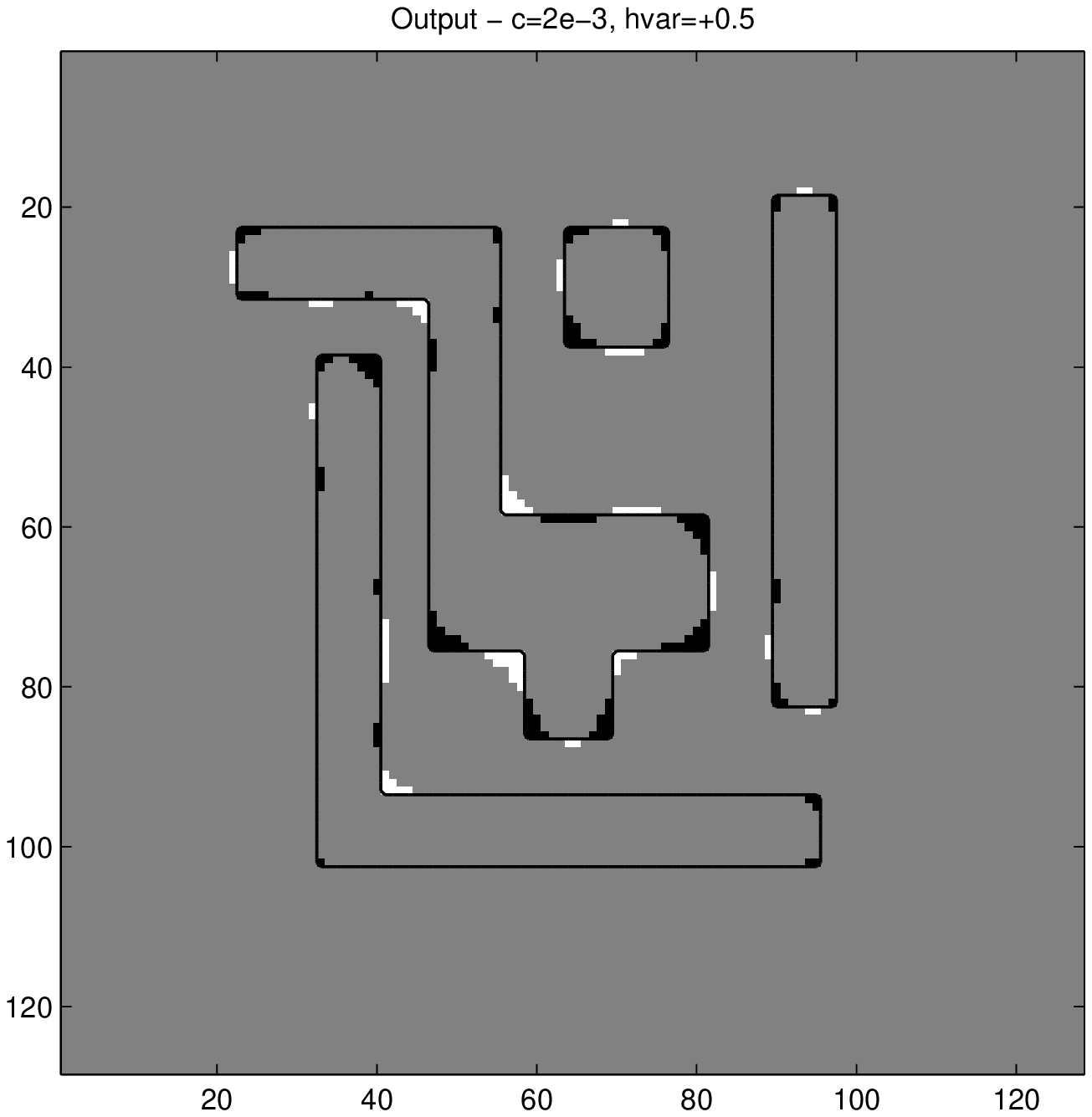}
\end{minipage}
\begin{minipage}[ht]{0.31\linewidth}
\includegraphics[width=1\textwidth]{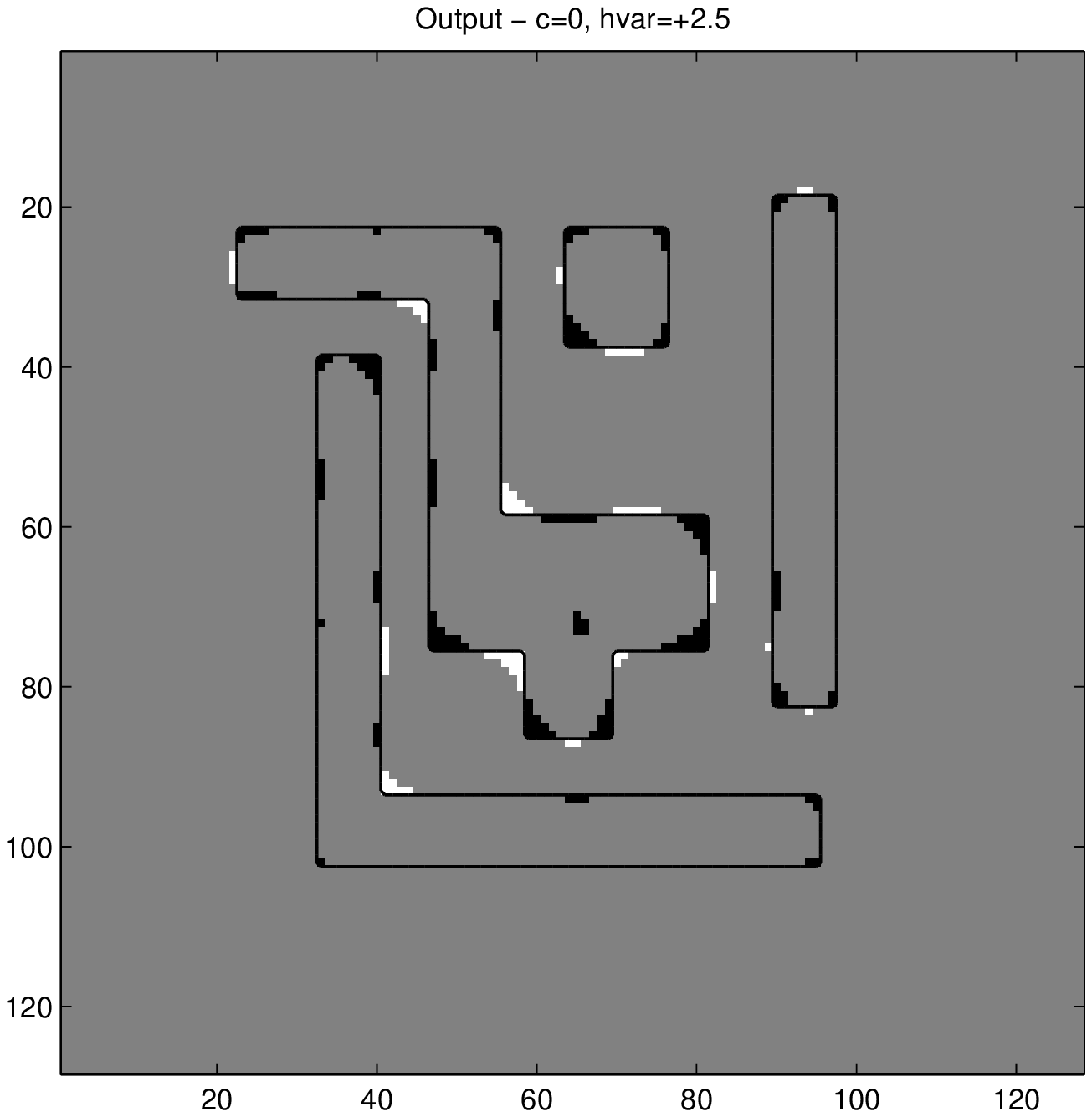}

\includegraphics[width=1\textwidth]{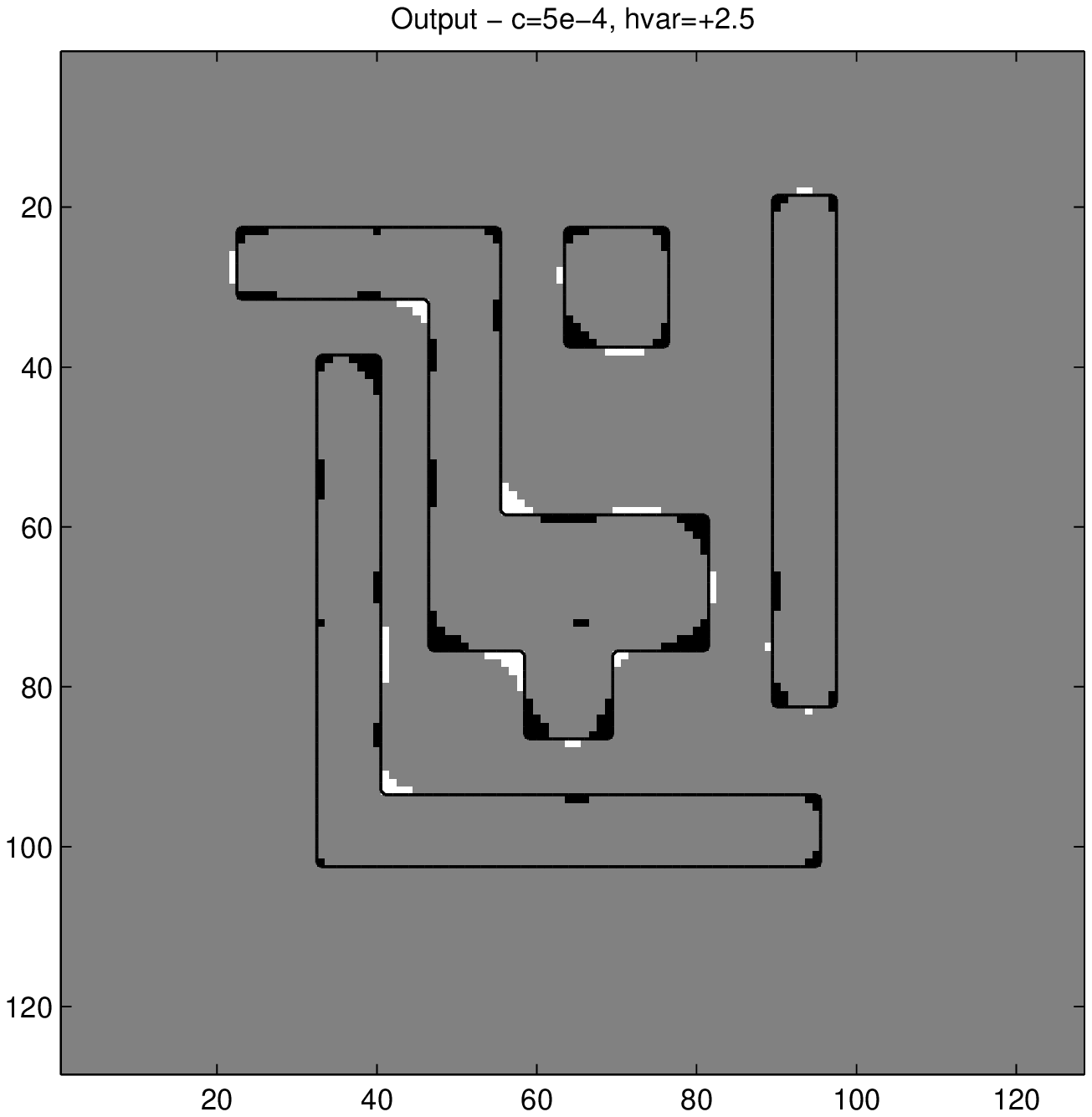}

\includegraphics[width=1\textwidth]{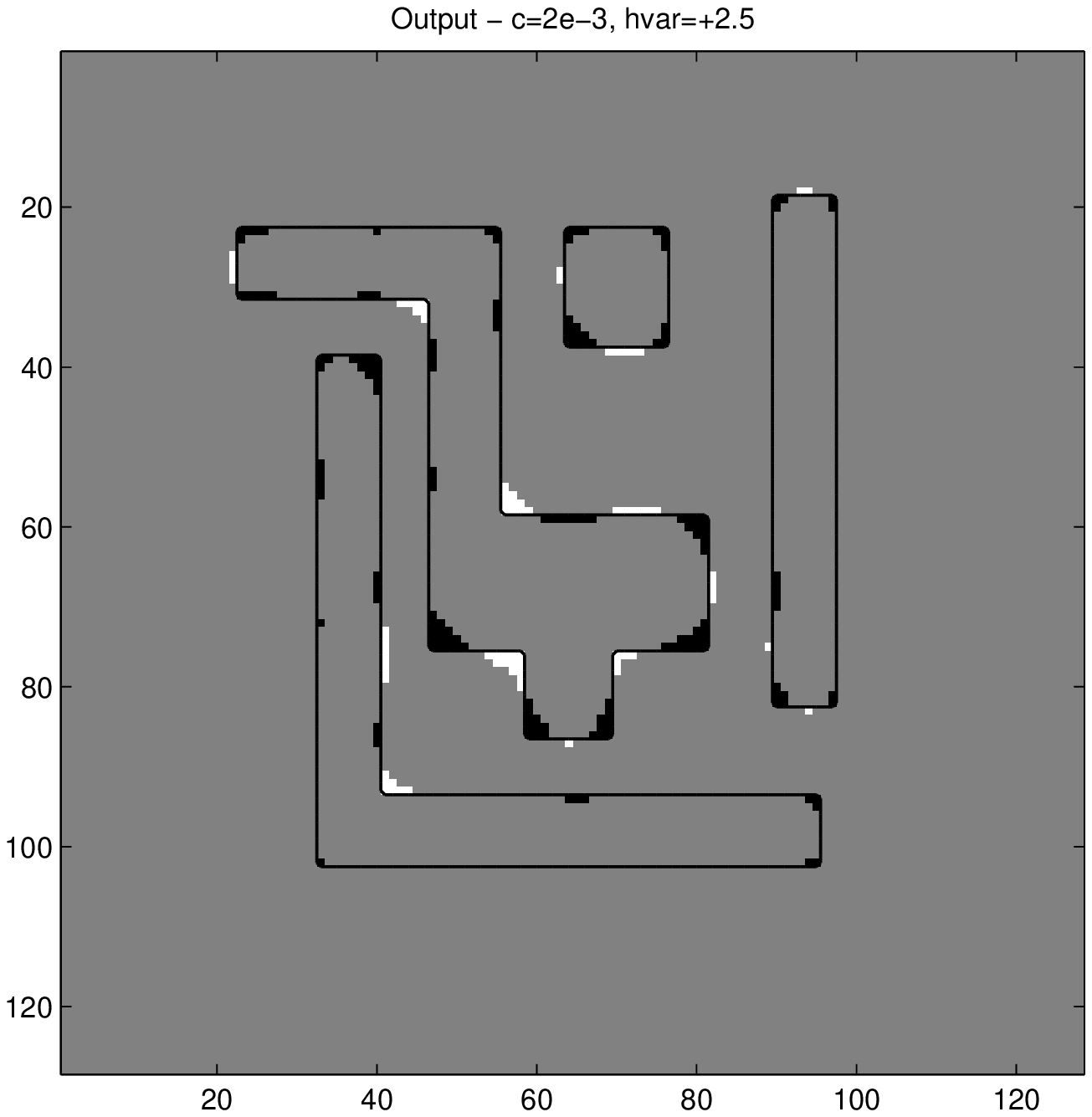}
\end{minipage}\caption{Top: {\it Test n.3} ($c=0$). Center: {\it Test n.4} ($c=5\times 10^{-4}$).
Bottom: {\it Test n.5} ($c=2\times 10^{-3}$).
Left: $hvar=0$. Middle: $hvar=0.5$. Right: $hvar=2.5$.
 }
 \label{Fig:change_c1}
\end{figure}

We remark that there seems to be not much difference in the masks
(which are not shown) and the outputs (the error in pixels is 227 for
{\it Test n.3}, 225 for {\it Test n.4} and 227 again for {\it Test
  n.5}). However, $\mathcal{R}_{\gamma}$ prevents the threshold from
being a critical value. In fact, the minimal value of the function $d$
defined above in \eqref{d_defin} goes from 1.27\% in {\it Test n.3} to
2.24\% in {\it Test n.4} and finally to 4.35\% in {\it Test n.5}.
The benefit of the penalty is 
stability with respect to the changes of the
threshold $h$ as we shall shortly see. We change the value of the
threshold by a percentage value of $hvar$. The outcome is shown in
Figure~\ref{Fig:change_c1}. On the top we have {\it Test n.3} (with
$c=0$), in the center we have {\it Test n.4} ($c=5\times 10^{-4})$ and on
the bottom we have {\it Test n.5} ($c=2\times 10^{-3}$). From left to
right we see how the reconstruction changes if we vary the value of
threshold. On left the threshold is $h$ (corresponding to $hvar=0$),
in the middle it is $(100.5/100)h$ ($hvar=0.5$), and on the right it
is $(102.5/100)h $ ($hvar=2.5$).  Even if the improvement by
increasing the parameter $c$ is not that striking from the point of
view of the error in pixels, from a topological point of view it is
actually remarkable.

\begin{figure}[htb]
\centering
\begin{minipage}[ht]{0.31\linewidth}
\includegraphics[width=1\textwidth]{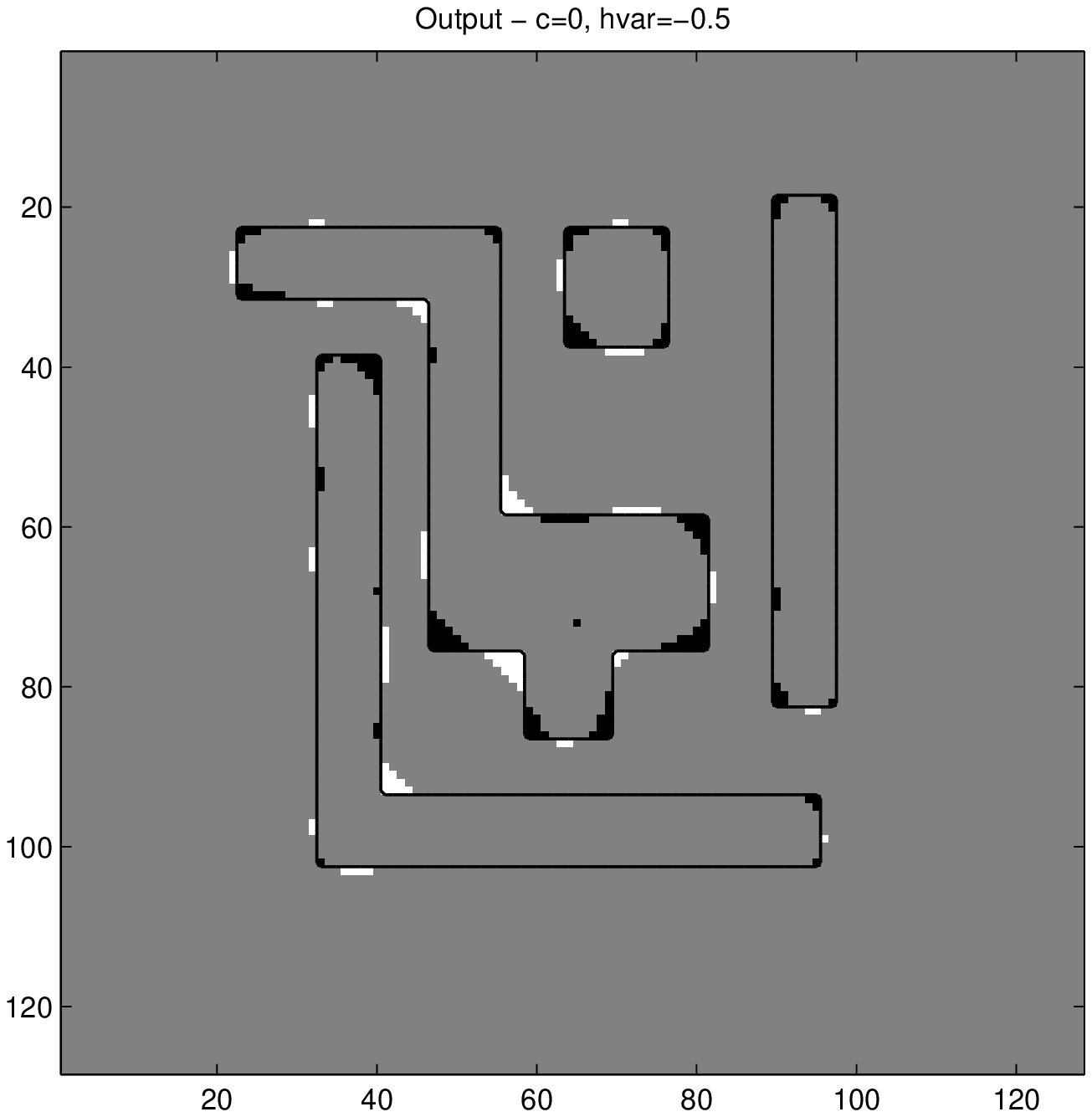}

\includegraphics[width=1\textwidth]{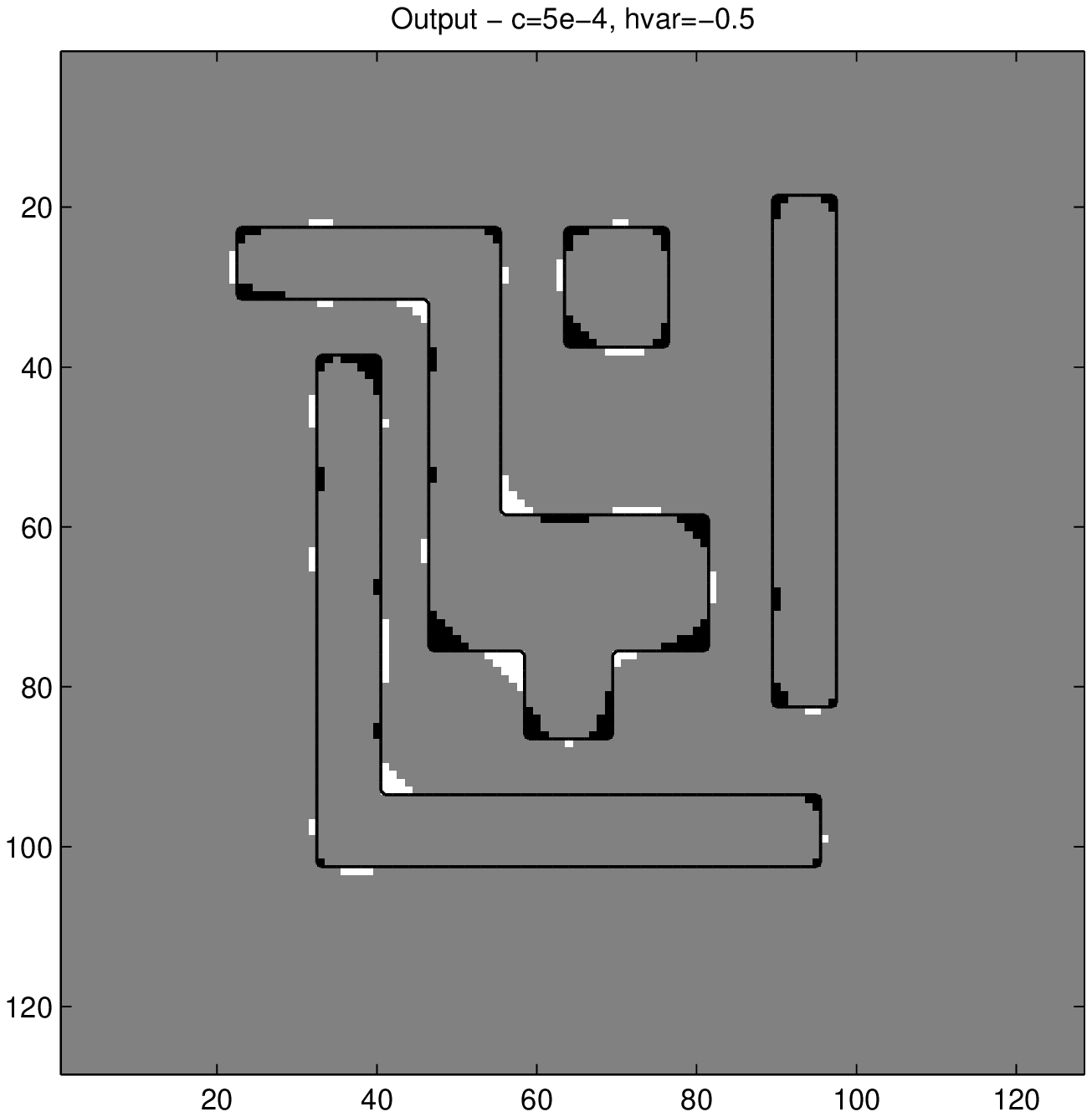}

\includegraphics[width=1\textwidth]{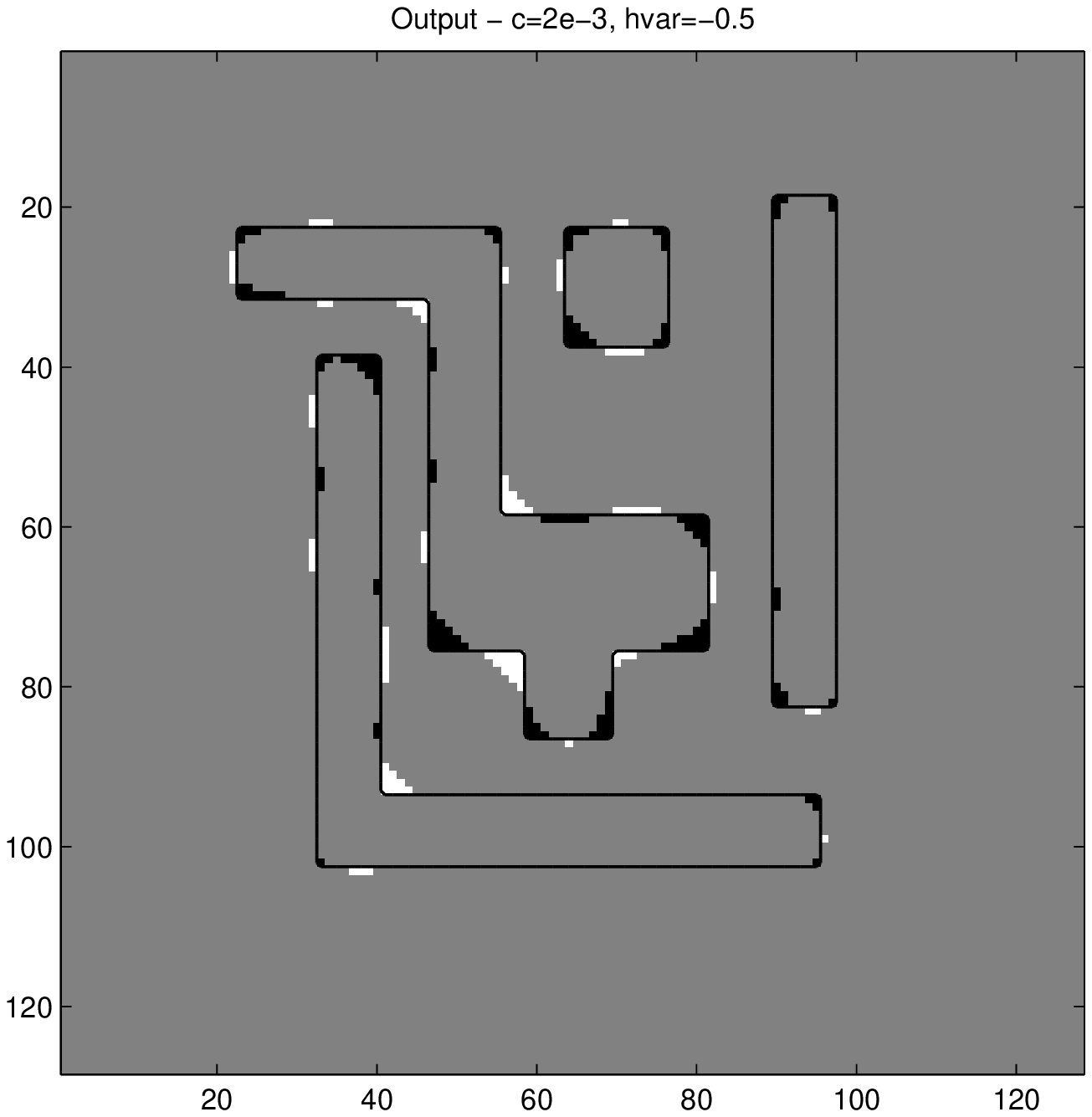}
\end{minipage}
\begin{minipage}[ht]{0.31\linewidth}
\includegraphics[width=1\textwidth]{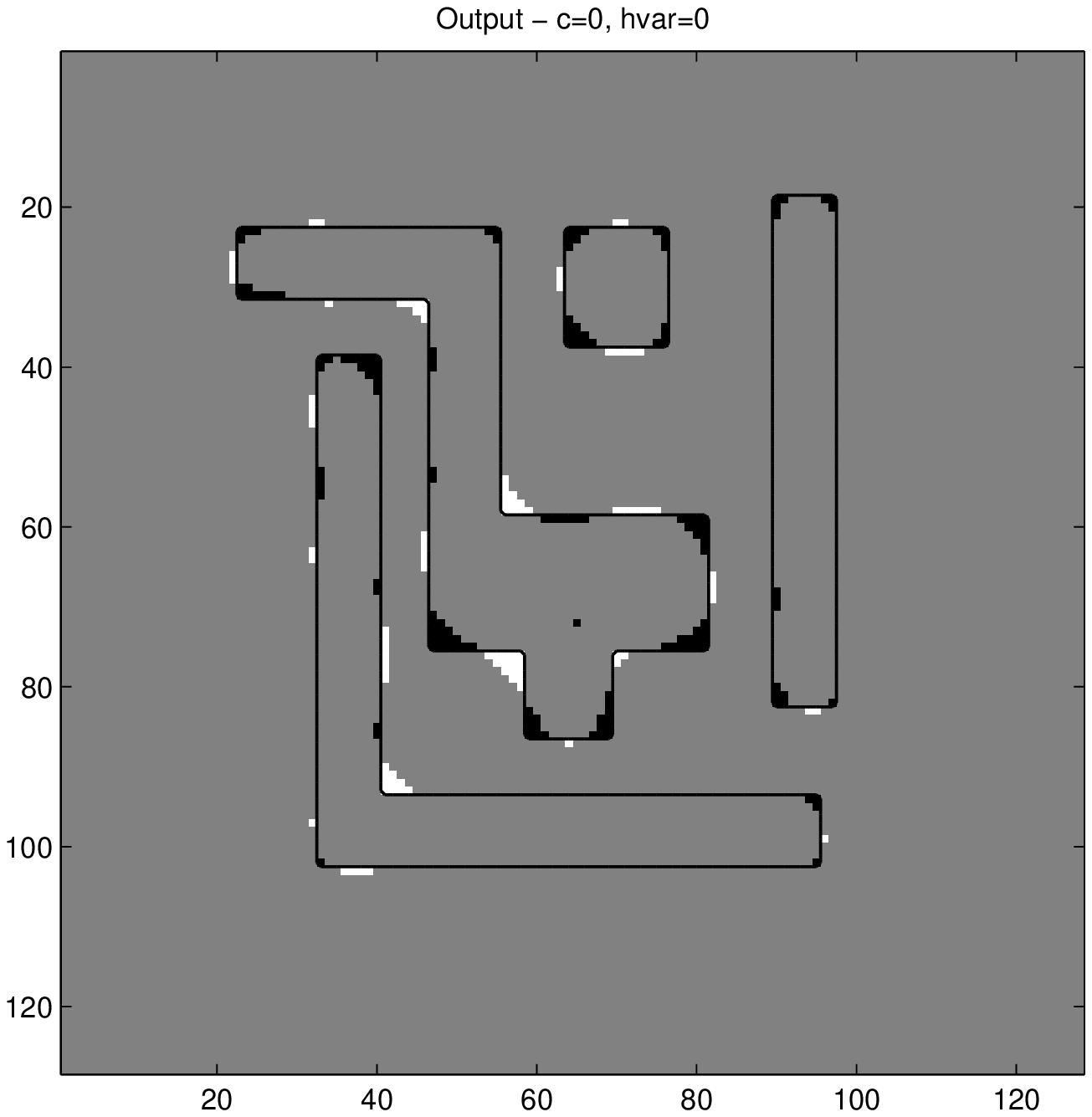}

\includegraphics[width=1\textwidth]{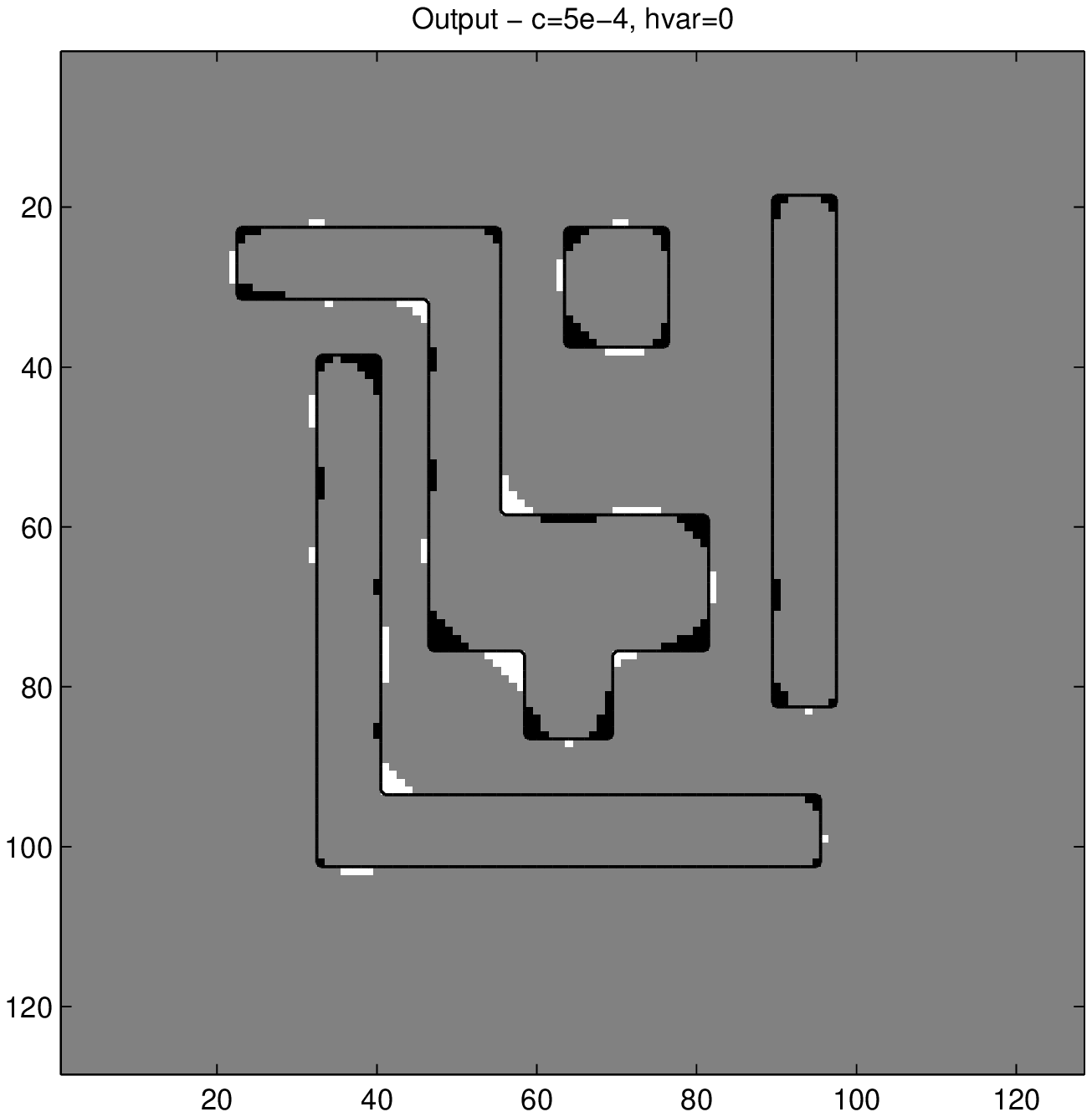}

\includegraphics[width=1\textwidth]{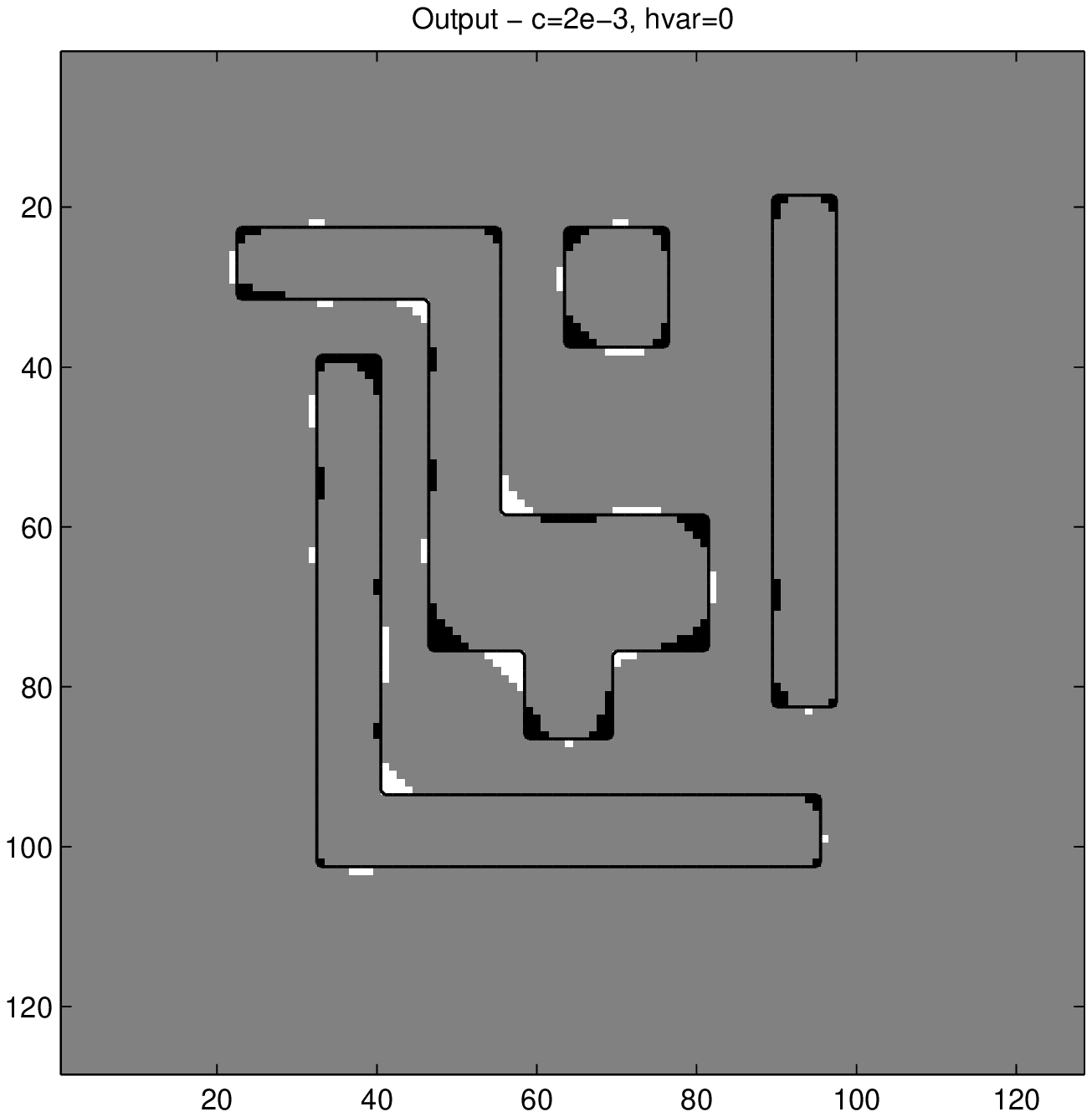}
\end{minipage}
\begin{minipage}[ht]{0.31\linewidth}
\includegraphics[width=1\textwidth]{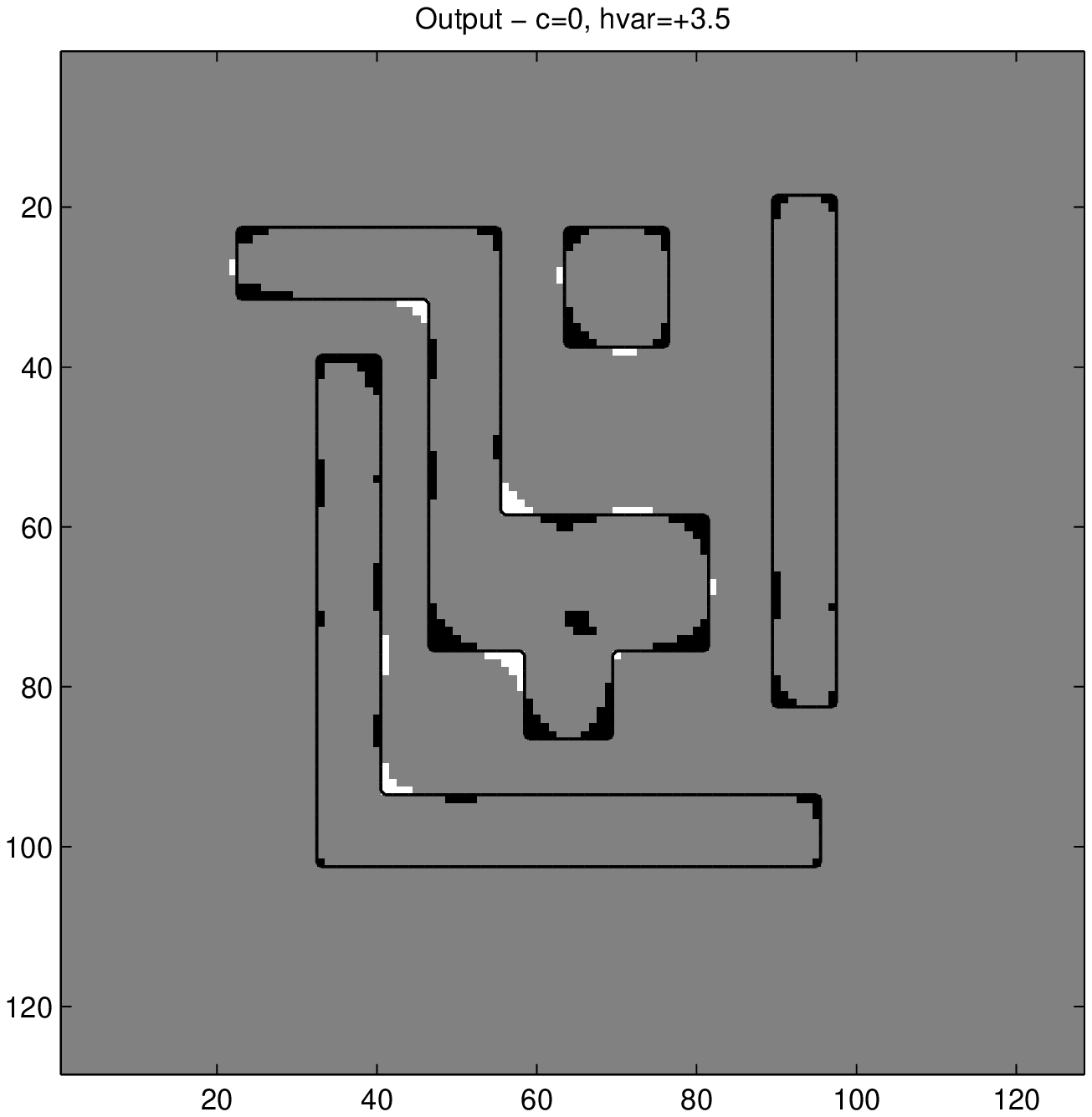}

\includegraphics[width=1\textwidth]{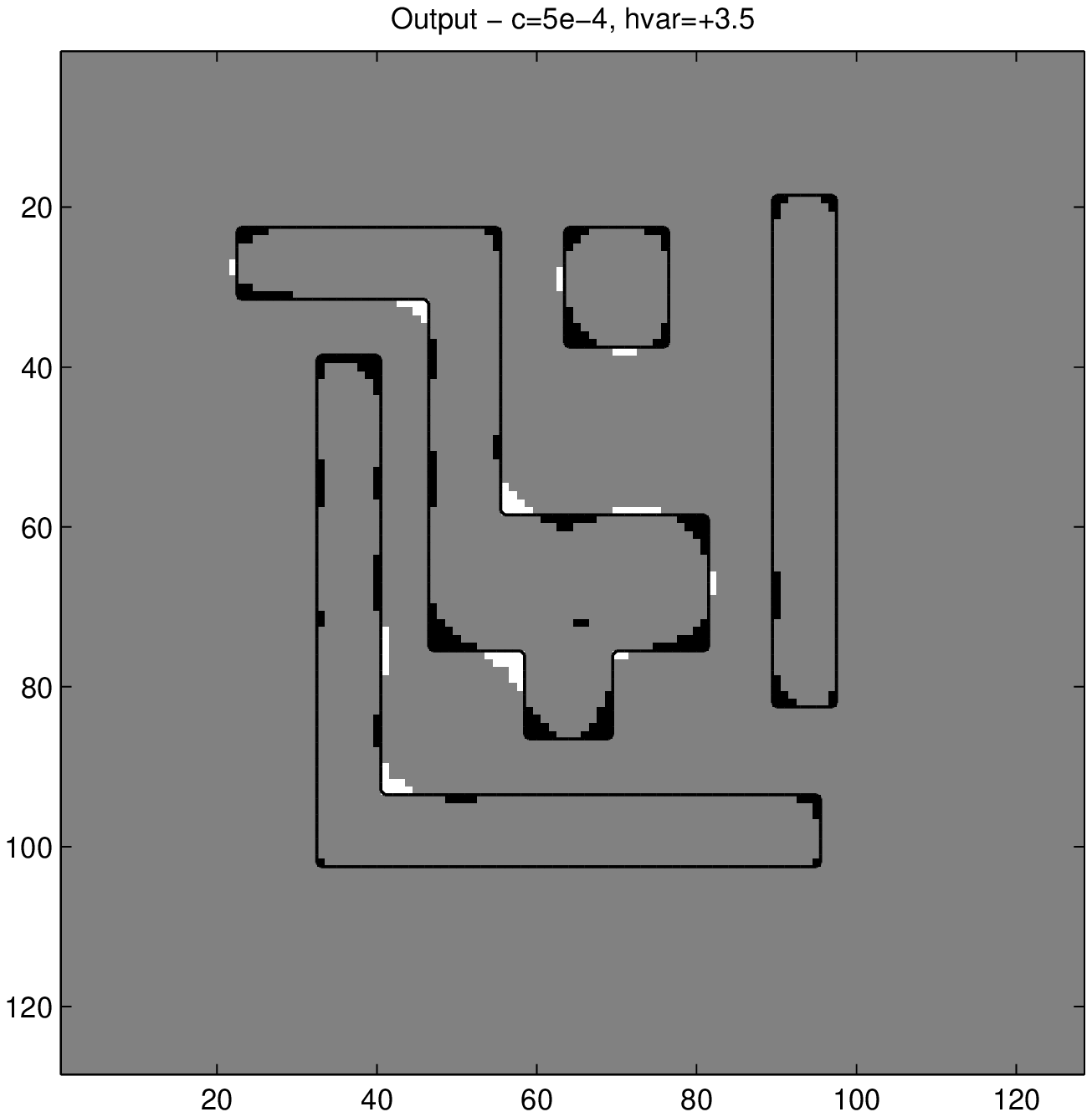}

\includegraphics[width=1\textwidth]{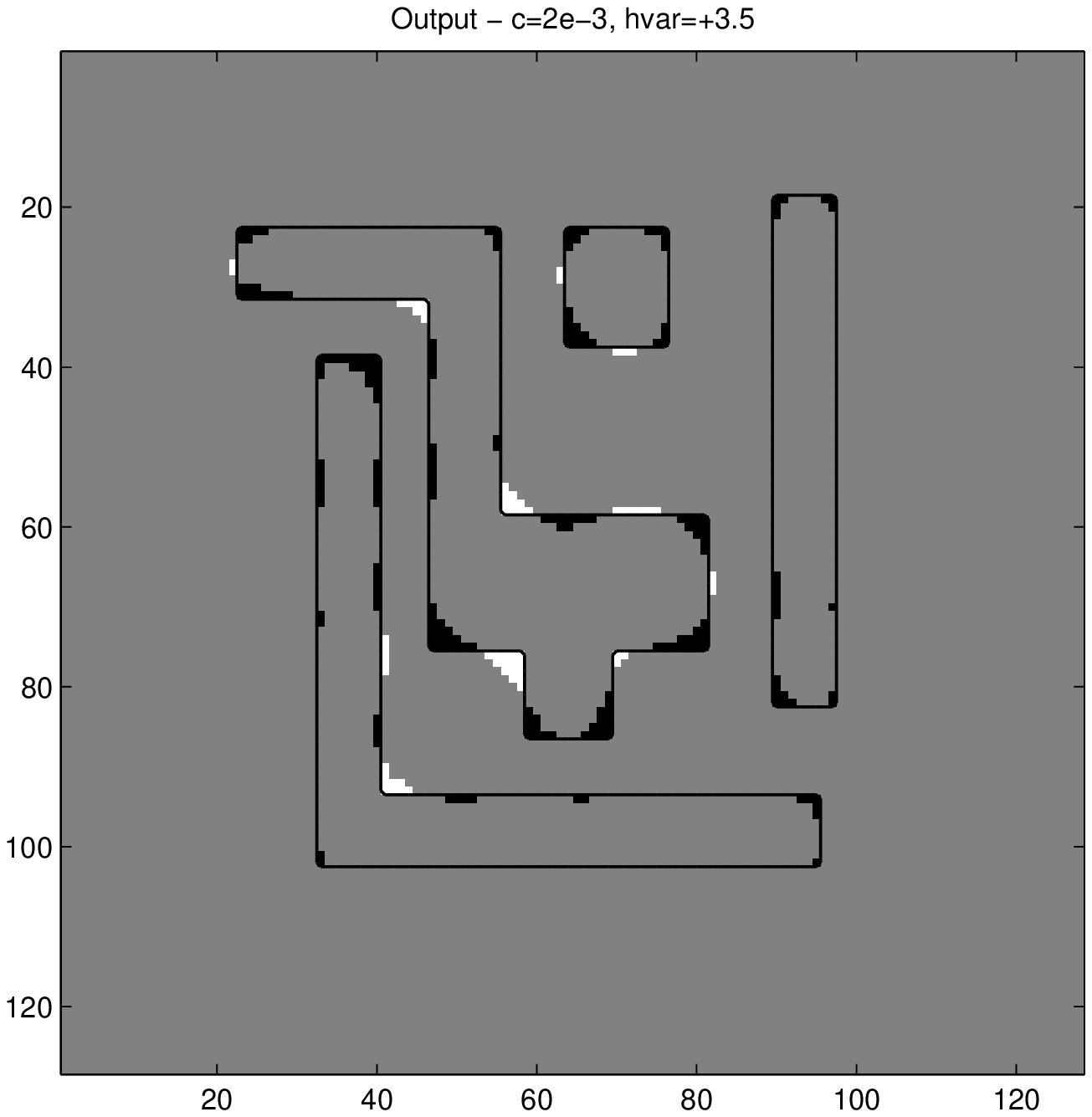}
\end{minipage}
\caption{Top: {\it Test n.6} ($c=0$). Center: {\it Test n.7} ($c=5\times 10^{-4}$).
Bottom: {\it Test n.8} ($c=2\times 10^{-3}$).
Left: $hvar=-0.5$. Middle: $hvar=0$. Right: $hvar=3.5$.
}
 \label{Fig:change_c2}
\end{figure}

This is even more striking in this other example where we decrease
$\varepsilon$, $\eta$ and $\gamma$ faster, by using
$rate_{\varepsilon}=1.5$, $rate_{\eta}=1.5$ and $rate_{\gamma}=1.1$,
and we perform $780$ iterations in total. Keeping all the other
parameters fixed, we call {\it Test n.6} the one with $c=0$, {\it Test
  n.7} the one with $c=5\times 10^{-4}$ and, finally, {\it Test n.8}
the one with $c=2\times 10^{-3}$. The outcome is shown in
Figure~\ref{Fig:change_c2}. On the top we have {\it Test n.6} (with
$c=0$), in the center we have {\it Test n.7} ($c=5\times 10^{-4}$) and
on the bottom we have {\it Test n.8} ($c=2\times 10^{-3}$). From left
to right we see how the reconstruction changes if we vary the value of
threshold. On the left the threshold is $(99.5/100)h$ ($hvar=-0.5$), in
the middle it is $h$ ($hvar=0$), and on the right it is $(103.5/100)h
$ ($hvar=3.5$).

In {\it Test n.6}, without the regularization term $\mathcal{R}_{\gamma}$, the 
hole appears even if we take a threshold lower than $h$. The hole is not present for threshold $h$ if we add $\mathcal{R}_{\gamma}$ with a small coefficient and it is not present for a considerably higher value of the threshold ($+3.5\%$) if the coefficient of $\mathcal{R}_{\gamma}$ is slightly bigger.

\begin{figure}[htb]
\centering
\begin{minipage}[ht]{0.31\linewidth}
\includegraphics[width=1\textwidth]{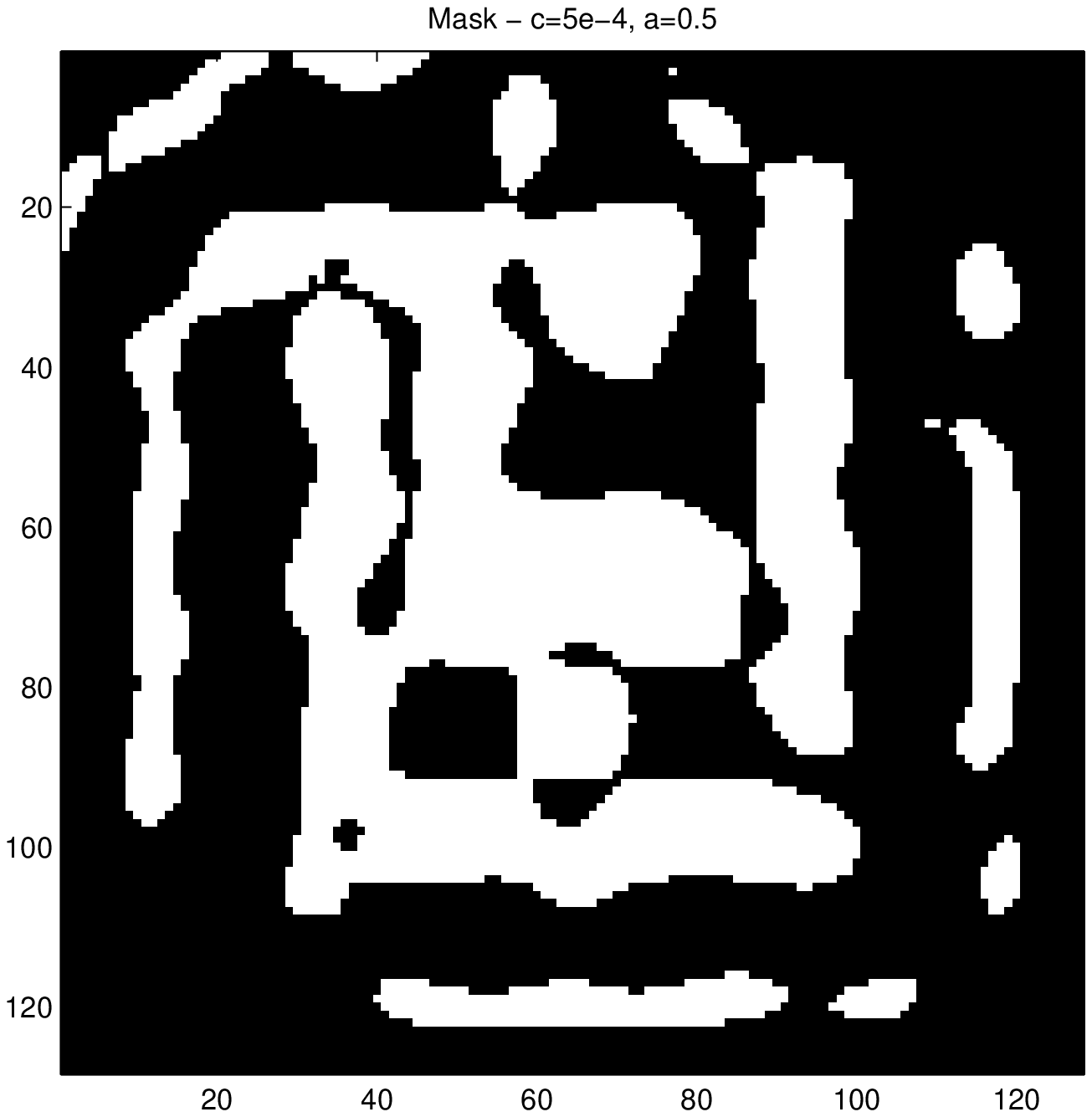}
\end{minipage}
\begin{minipage}[ht]{0.20\linewidth}\phantom{aaaaaaaaaaaaaa}
\end{minipage}
\begin{minipage}[ht]{0.31\linewidth}
\includegraphics[width=1\textwidth]{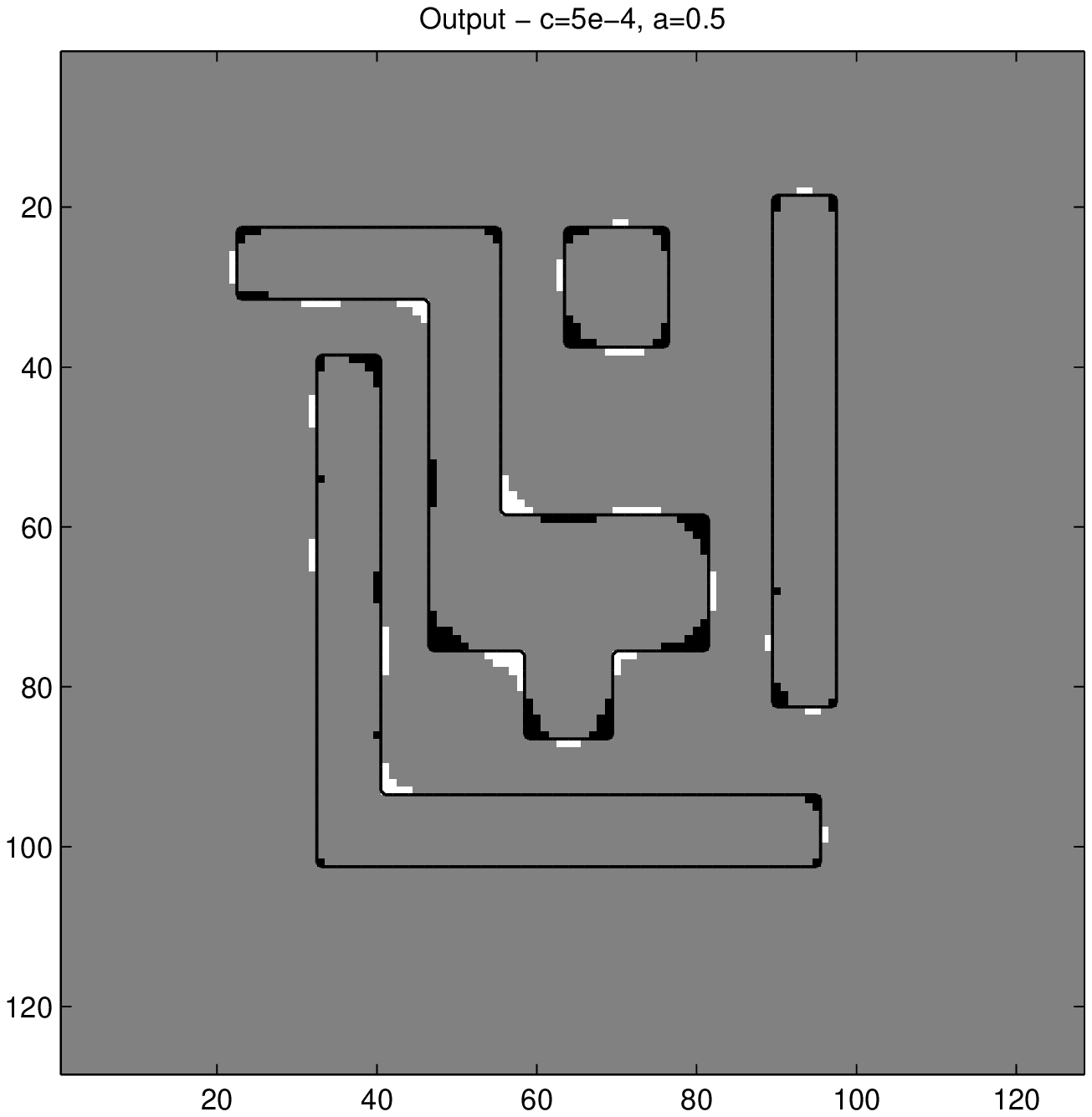}
\end{minipage}
\caption{{\it Test n.9} ($c=5\times 10^{-4}$, $a=0.5$). Left: mask. Right: output.}\label{Fig:perimeter}
\end{figure}

So far we have kept the coefficient $a$ equal to $0$.
In fact, in our experiments we see that the term containing the difference between the perimeters in the definition of the distance function $d_{st}^2$ actually does not play a big role and in general does not improve the reconstruction. However for completeness we show the outcome of an experiment where the full functional is used, namely we modify {\it Test n.4} above by changing the parameter $a$ from $0$ to $0.5$. The error in pixels is 232 in this case and the outcome is illustrated in Figure~\ref{Fig:perimeter}.

\subsubsection*{Conclusions}
From our numerical experiments we can draw the following general
conclusions. 
\begin{enumerate}
\item The outcome mask is very diffuse and is not at all close to target.
While the optimal mask has shapes much more complicated than the target, the
exposed region is close to the target. The main
reason for these two facts is the high nonlocality of the image intensity.
\item The final shape of the mask strongly depends on the initial guess, since
the functional $F_0$ is nonconvex therefore may have several absolute and local minimizers.
\item The presence of several local minima, in this case of the approximated
functional $F_{\varepsilon}$, has also the effect that sometimes we do
not have convergence to a perfectly binary function. However, the
discrepancy with a binary function is limited to very few pixels.
\item We observed that the reconstruction is in general better when the parameters decrease slowly and uniformly.
\item The effect of the term containing the difference between the two perimeters
  in the definition of the distance $d_{st}^2$ does not have a pronounced
influence on the optimal mask.
This can be understood from the fact that this term is
the difference between two numbers which are not very local.
\end{enumerate}

\section{Discussion}\label{conclusionsec}
In this paper we studied the inverse problem of photolithography, which can be viewed
as an optimal shape design problem.  A main novelty of the paper is the regularization 
term $\mathcal{R}$, which has both theoretical and practical value. In solving the inverse
problem, the penalty term has a desirable stabilizing influence.

When the threshold $h$ is not close to a critical value of the
intensity, the penalty term has no effect. This is what happens when
we performed the computation using Target 1. For Target 2, the
intensity has a local minimum inside the biggest feature with a local
minimum value very close to the threshold. This is why a hole may
appear in the reconstruction for small perturbations of the
threshold. In this case the term $d$ defined in \eqref{d_defin} is
very small at this local minimum point, therefore $d_{min}$, the
minimum value of $d$, is very close to $0$.  It happens that the term
$\mathcal{R}_{\gamma}$ raises the value of $d_{min}$, essentially by
pushing away, and actually up, the local minimum value from the
threshold value. As a practical effect, the hole will not show up even
at a higher perturbation of the threshold. Therefore we greatly
improve the topological stability of the reconstruction by adding the
term $\mathcal{R}$.

\subsubsection*{Acknowledgements}
The authors thank Hande T\"uzel for providing the codes to compute the
Hopkins aerial intensity.  L.R. is partly supported by Universit\`a
degli Studi di Trieste through Fondo per la Ricerca di Ateneo --- FRA
2012 and by GNAMPA, INdAM.  The research of F.S. is funded in part by
NSF Award DMS-1211884.  This research was started at
the 
Institute for Mathematics and its Applications (IMA) when Z.W. was a
postdoctoral fellow. The IMA receives funding from the NSF under Award
DMS-0931945.

\end{document}